\documentclass[12pt]{amsart}

\usepackage{amsfonts}
\usepackage{amsmath}
\usepackage{amssymb}
\usepackage{mathabx}

\usepackage{mathrsfs}
\usepackage{amscd}
\usepackage[all]{xy}
\usepackage[pdftex,final]{graphicx}
\usepackage{hyperref}
\usepackage[hyperpageref]{backref}

\newtheorem{lem}{Lemma}[section]
\newtheorem{thm}[lem]{Theorem}
\newtheorem{prop}[lem]{Proposition}
\newtheorem{cor}[lem]{Corollary}
\theoremstyle{definition}

\newtheorem{exa}[lem]{Example}

\newtheorem{rem}[lem]{Remark}

\newcommand{\imp}{\Longrightarrow}

\newcommand{\Q}{\Bbb{Q}}
\newcommand{\F}[1]{\Bbb{F}_{#1}}

\newcommand{\Z}{\Bbb{Z}}

\newcommand{\K}{\mathcal{K}}

\newcommand{\LL}{L^{\tau}}
\renewcommand{\S}{\Sigma}

\newcommand{\wn}[1]{\mathcal{W}_{#1}}% set of units u such that 1-u is a unit 

\newcommand{\oct}{\mathcal{O}}%   octahedra in H_2(Y) for Y 2-dimensional
\newcommand{\cyc}[1]{\left[ #1\right]} %cycle associated to an octahedron
\newcommand{\scyc}[1]{\left[ #1 \right]} %cycle of basic octahedron
\newcommand{\octt}[1]{\left[ #1\right]}% basic 27 encoded cycles of Z/8
\newcommand{\spl}[2]{\operatorname{SL}_{#1}\left(#2\right)}
\newcommand{\pspl}[2]{\mathrm{PSL}_{#1}(#2)}
\newcommand{\gl}[2]{\mathrm{GL}_{#1}(#2)}
\newcommand{\pgl}[2]{\mathrm{PGL}_{#1}(#2)}

\newcommand{\pb}[1]{\mathcal{P}(#1)}%  old pre-Bloch group

\newcommand{\qpb}[1]{\widetilde{\mathcal{P}}(#1)}%  quasi-reducedPre-Bloch Group

\newcommand{\rpb}[1]{\mathcal{RP}(#1)}%old refined pre-Bloch group

\newcommand{\rpbker}[1]{\mathcal{RP}_1(#1)}% old  RP_1(F)
\newcommand{\crpbker}[1]{\widecheck{\mathcal{RP}}_1(#1)}

\newcommand{\crpb}[1]{\widecheck{\mathcal{RP}}(#1)}
\newcommand{\cpb}[1]{\widecheck{\mathcal{P}}(#1)}

\newcommand{\bl}[1]{\mathcal{B}(#1)}% old Bloch group

\newcommand{\rbl}[1]{{\mathcal{RB}}(#1)}%old refined Bloch group
\newcommand{\crbl}[1]{\widecheck{\mathcal{RB}}(#1)}

\newcommand{\gpb}[1]{\left[ #1\right]} %notation for generators of (pre)Bloch

\newcommand{\sus}[1]{\psi(#1) }% "Suslin element"
\newcommand{\suss}[2]{\psi_{#1}\left( #2\right)}
\newcommand{\Suss}[2]{\Psi_{#1}\left( #2\right)}

\newcommand{\ks}[2]{{\K}^{\mbox{\tiny $(#1)$}}_{#2}}
\newcommand{\cks}[2]{\widecheck{{\K}}^{\mbox{\tiny $(#1)$}}_{#2}}

\newcommand{\cjs}[1]{\widecheck{\mathcal{J}}_{#1}}

\newcommand{\bconst}[1]{C_{#1}}
\newcommand{\cconst}[1]{D_{#1}}
\newcommand{\cconstmod}[1]{\mathcal{D}_{#1}}
\newcommand{\ccconstmod}[1]{\widecheck{\mathcal{D}}_{#1}}

\newcommand{\pn}[2]{\mathbb{P}^{#1}(#2)}
\newcommand{\projl}[1]{\pn{1}{#1}}

\newcommand{\card}[1]{\left| #1 \right|}

\newcommand{\id}[1]{\mathrm{Id}_{#1}}

\renewcommand{\ker}[1]{\mathrm{Ker}(#1)}
\newcommand{\image}[1]{\mathrm{Im}(#1)}
\newcommand{\Ker}[1]{\mathrm{Ker}\left( #1\right)}
\newcommand{\coker}[1]{\mathrm{Coker}\left(#1\right)}

\newcommand{\sqcl}[1]{G_{#1}}% group of square classes
\newcommand{\psqcl}[1]{\bar{G}_{#1}}% plus or minus group of square classes

\newcommand{\sgr}[1]{\mathrm{R}_{#1}}% group ring of square classes of field
\newcommand{\rsgr}[1]{\widecheck{\mathrm{R}}_{#1}}% quotient of group ring of square classes of field

\newcommand{\qsgr}[1]{\mathrm{R}^{+}_{#1}} % `quasireduced' group ring = plus group ring of A
\newcommand{\msgr}[1]{\mathrm{R}^{-}_{#1}} % minus group ring of A
\newcommand{\tsgr}[1]{\widecheck{\mathrm{R}}^{+}_{#1}} % `reduced' group ring 
\newcommand{\an}[1]{\left\langle{#1}\right\rangle}
\newcommand{\pf}[1]{\left\langle\!\left\langle{#1}\right\rangle\!\right\rangle}
\newcommand{\aug}[1]{\mathcal{I}_{#1}}
\newcommand{\qaug}[1]{\widetilde{\mathcal{I}}_{#1}}
\newcommand{\asym}[3]{\mathrm{S}^{#1}_{#2}(#3)}

\newcommand{\sgn}[1]{\mathrm{sgn}(#1)}

\newcommand{\supp}[1]{\mathrm{supp}(#1)}

\newcommand{\gab}[1]{{#1}^{\mbox{\tiny{ab}}}}% abelianization of a group

\newcommand{\lconj}[1]{\mathrm{lc}_{#1}}
\newcommand{\rconj}[1]{\mathrm{rc}_{#1}}

\newcommand{\odet}{\overline{\mathrm{det}}}
\newcommand{\oodet}{\widetilde{\mathrm{det}}}
\newcommand{\cls}[2]{\mathrm{Cl}^{#1}(#2)}

\newcommand{\kind}[1]{K^{\mathrm{\small ind}}_3(#1)}

\newcommand{\ho}[3]{\mathrm{H}_{#1}\left( #2,#3 \right)}

\newcommand{\Ind}[3]{\mathrm{Ind}_{#1}^{#2}#3}

\newcommand{\tors}[2]{{#1}[#2]}

\newcommand{\op}[1]{{#1}^{\mathrm{op}}}

\newcommand{\Cor}[2]{\mathrm{cor}_{#1}^{#2}}
\newcommand{\Res}[2]{\mathrm{res}_{#1}^{#2}}

\renewcommand{\sgn}[1]{\mathrm{sgn}(#1)}

\newcommand{\mc}[2]{\Sigma_{#1}(#2)}
\newcommand{\pmc}[2]{\widehat{|\mc{#1}{#2}|}}

\title{Bloch Groups of local rings with residue field $\F{2}$ }
\author{  Kevin Hutchinson}
\address{School of Mathematics and Statistics,
 University College Dublin}
\email{  kevin.hutchinson@ucd.ie}
\date{\today}


\begin{document}
\maketitle
\tableofcontents
\begin{abstract} We give a generic presentation of the pre-Bloch group of a local ring with residue field $\F{2}$. We apply this to the calculation of pre-Bloch groups and Bloch groups of some local rings.
\end{abstract}
\section{Introduction} 
The scissors congruence group, or pre-Bloch group, $\pb{F}$ of a field $F$, and its subgroup  the Bloch group, $\bl{F}$, (see \cite{sah:dupont}) is of interest because it describes the indecomposable $K_3$ of the field, $\kind{F}$,  modulo some known torsion coming from the group $\mu_F$ of roots of unity (\cite{sus:bloch}). 
%(See section \ref{sec:kind} below for a more precise statement.)
%  For any field $F$, $\kind{F}$ is naturally a homomorphic image of the homology group $\ho{3}{\spl{2}{F}}{\Z}$.
 The refined scissors congruence group $\rpb{F}$  and the refined Bloch group $\rbl{F}$ were introduced in \cite{hut:cplx13} to help to understand the kernel of the natural 
surjective homomorphism $\ho{3}{\spl{2}{F}}{\Z}\to\kind{F}$. These groups and their properties were then used to calculate $\ho{3}{\spl{2}{F}}{\Z[\frac{1}{2}]}$ for local fields $F$ as well as well as for $F=\Q$ (\cite{hut:rbl11},\cite{hut:sl2Q}). The definitions and techniques in these results extend readily from fields to local rings with sufficiently large residue fields and thus allow us to describe $\kind{A}$ and $\ho{3}{\spl{2}{A}}{\Z[\frac{1}{2}]}$ for many local rings $A$ (\cite{mirzaii:bwlocal},\cite{hut:slr}, \cite{hut:mokarimirzaii}, \cite{mirzaii:perez25}) . 
However, the given funtorial definition of the (refined) pre- Bloch group does not extend to the fields $\F{2}$ and $\F{3}$ or to local rings with small residue fields.

In \cite{hut:bloch1} the author and  Rodrigo Cuitun Coronado defined functorially the pre-Bloch groups $\rpb{A}$ (and $\pb{A}$)  and the Bloch groups $\rbl{A}$ (and $\bl{A}$) of \emph{any} commutative ring $A$.   Under suitable hypotheses on the ring $A$, there is a natural homomorphism\\
 $\ho{3}{\spl{2}{A}}{\Z}\to \rpb{A}$ whose image is $\rbl{A}$ (or from $\ho{3}{\gl{2}{A}}{\Z}$ to $\pb{A}$ with image $\bl{A}$).  Those suitable hypotheses, which are satisfied by all fields and local rings, concern the low-dimensional homology of a simplicial complex $Y(A)$ associated to the ring $A$, introduced in \cite{hut:ge2}, and whose fundamental group can be described in terms of the algebraic $K$-theory of $A$ (\cite[Theorem 6.9]{hut:ge2}).

When $A$ is a field with at  least four elements or, more generally a local ring whose residue field has at least four elements, then $\pb{A}$ has a natural presentation in terms of generators $\gpb{u}$, $u\in \wn{A}:=\{ a\in A^\times|\  1-a\in A^\times\}$ (the set of exceptional units), and a family of $5$-term dilogarithm relations  $R_{u,v}$ where $u,v\in\wn{A}$ also satisfy $u/v\in \wn{A}$. Likewise, $\rpb{A}$ has an analogous presentation as an $\sgr{A}:=\Z[A^\times/(A^\times)^2]$-module. These presentations have proved useful in calculating $K$-theory groups and the third homology of $\spl{2}{A}$ (see, in particular, \cite{hut:sl2Q}).

 However, when $A$ is a ring mapping to $\F{2}=\Z/2$ the set $\wn{A}$ is necesssarily empty and the proposed generators $\gpb{u}$ do not exist. When $A$ admits a homomorphism to $\F{3}=\Z/3$ and if $u,v\in\wn{A}$ then $u/v\not\in \wn{A}$. Thus when $A$ is a local ring with residue field $\F{3}$, one can show that $\rpb{A}$ is generated  by the set $\{ \gpb{u}\ |\ u\in\wn{A}\}$ but that the above $5$-term relations no longer exist.

Our main objective in the present paper is to give a presentation of $\rpb{A}$ (and of $\pb{A}$ which is a quotient of $\rpb{A}$) in the case of a local ring $A$ with residue field $\F{2}$, and to do illustrate this result by doing some particular calculations. The presentation is achieved  in Theorem \ref{thm:rpaf2} below. The space $Y(A)$ mentioned above plays a crucial role here. Along with some standard familiess of elements which exist  in $\rpb{A}$ for any ring $A$ the additional exceptional generators correspond to certain canonical generators of $H_2(|Y(A)|)$. When $A$ is a local ring with residue field $\F{2}$, $Y(A)$ is homotopy equivalent to a pointed sum of $2$-spheres. In fact the generating $2$-spheres are represented by the surfaces of octahedra in the simplicial complex $Y(A)$ and the corresponding homology classes lift naturally to elements of $\rpb{A}$ which we call \emph{standard octahedra}. Then these octahedra, together with the standard families of elements mentioned, give a complete explicit generating set of $\rpb{A}$ (and of $\pb{A}$); see Corollary \ref{cor:localf2}.

In this article we also apply the general results to calculate Bloch groups, and the homology of $\spl{2}{A}$, for the local rings $A=\Z/2^n$ for $n=1,2,3$ and also for the ring $\Z_{\an{2}}$ of rational numbers with odd denominator.

In a subsequent paper, we will give a presentation of $\rpb{A}$ for local rings $A$ with residue field $\F{3}$  (\cite{hut:blochf3}).

In section \ref{sec:ya} we review the definition and basic properties of the simplicial complex $Y=Y(A)$. In particular $\rpb{A}:=Z_2(L_\bullet)_{\spl{2}{A}}$ where  $L_\bullet=L_\bullet(Y)$ is the complex of \emph{non-degenerate} ordered simplices. When $A$ is a local ring $Y(A)$ is the clique complex of complete multipartite graph.
In appendix \ref{app:clique} we show that the clique complex of a complete multipartite graph with $q+1$ parts  is homotopy equivalent to  a pointed sum of $q$-spheres and show how to write down generators for its top-dimensional homology. We apply this in section \ref{sec:ya} to write down an explicit basis  of $H_q(|Y(A)|)$ when $A$ is a local ring whose residue field has $q$ elements  (Proposition \ref{prop:localhom}).

In section \ref{sec:lowdim} we develop some quite general results about simplicial complexes. We introduce the complex $\Sigma_\bullet$ of \emph{partially symmetric chains} of a simplicial complex $Y$ and prove that when $Y$ is $2$-dimensional and has the property that every $1$-simplex is contained in a $2$-simplex then there is a natural short exact sequence 
$0\to Z_2(\Sigma_\bullet)\to Z_2(L_\bullet)\to H_2(|Y|)\to 0$. The group $Z_2(\Sigma_\bullet)$ has a simple presentation  (Corollary \ref{cor:z2pres}).

In section \ref{sec:wna=0} we apply the results of section \ref{sec:lowdim} to the special case that $Y=Y(A)$ where $A$ is a ring for which $\wn{A}=\emptyset$ (which includes the case of local rings with residue field $\F{2}$).  In this case we deduce an exact sequence $\crpb{A}\to\rpb{A}\to H_2(|Y(A)|)_{\spl{2}{A}}\to 0$ and we prove a presentation of the $\sgr{A}$-module
$\crpb{A}$. This module is generated by the elements $\suss{i}{u}$, $i\in \{ 1,2\}$, $u\in A^\times$,  and the constant $\cconst{A}$ whcih are subject to certain universal families of relations. 

In section \ref{sec:crpba} we give an explicit description of the structure of $\crpb{A}=\ho{0}{\spl{2}{A}}{Z_2(\Sigma_\bullet)}$ as an abelian group and as an $\sgr{A}$-module
 (Proposition \ref{prop:crpbstruct}).  More generally, we describe some  
calculations of the groups
 $\ho{n}{\spl{2}{A}}{Z_2(\Sigma_\bullet)}$ for small $n$. These are relevant to certain spectral sequence calculations in applications.

In section \ref{sec:f2} we apply the previous results to prove Theorem \ref{thm:rpaf2} giving a presentation of the pre-Bloch group $\rpb{A}$ when $A=(A,{M},k=\F{2})$ is a local ring with residue field $\F{2}$. Along with the standard families of elements $\suss{i}{u}$  and $\cconst{A}$, the generators include the standard octahedra which can be represented by symbols $\scyc{\frac{1}{m}\ n\ u}$,
$m,n\in {M}$, $u\in A^\times$, which are subject to certain algebraic identities which we enumerate.

Finally, in section \ref{sec:exa}, we apply the general results to the calculation, in whole or in part, of the Bloch groups of  specific local rings with residue field $\F{2}$.

\begin{rem}
No AI tools or technologies were used in the production of this article.
\end{rem}
\subsection{Notation}  In this article all  rings are commutative and have a unit.

For a ring $A$, $A^\times$ denotes the group of units. $\sgr{A}$ denotes the group ring
$\Z[A^\times/(A^\times)^2]$ of the group of square classes of units. For a unit $u\in A^\times$, the image in $\sgr{A}$ will be denoted $\an{u}$. The augmentation ideal $\aug{A}$ of the ring $\sgr{A}$ is the additive subgroup generated by the elements $\pf{u}:=\an{u}-1$, $u\in A^\times$.  

 Let $A$ be a ring. Composing the homomorphism $\det:\gl{2}{A}\to A^\times$ with the quotient map $A^\times\to A^\times/(A^\times)^2$ gives a homomorphism
$\odet:\gl{2}{A}\to A^\times\to A^\times/(A^\times)^2$, which factors through $\pgl{2}{A}$ since the determinant of scalar matrix $u\cdot I$ is $u^2$.

%For $m\geq 0$, we will say that $Y$ \emph{satisfies (the condition) $(C)_m$}  if $Y^\alpha$ is path connected for all $\alpha\in Y_m$.
%%%%%%%%%%%%%%%%%%%%%%%%%%%%%%%%%%%%%%%%%%%%%%%%%%%%%%%%%%%%%%%%%%%%%%%%%%%%%%%%
\section{The simplicial complex $Y(A)$}\label{sec:ya}
%%%%%%%%%%%%%%%%%%%%%%%%%%%%%%%%%%%%%%%%%%%%%%%%%%%%%%%%%%%%%%%%%%%%%%%%%%%%%%
\subsection{Preliminaries on simplicial complexes}
Let $Y$ be a simplicial complex.  For $n\geq 0$, $Y_n$ will denote the set of $n$-simplices of $Y$. Let $C_\bullet=C_\bullet(Y)$  be the complex whose $n$-th term $C_n$ is the free abelian group on the set 
\[
E_n=E_n(Y):=\{ (x_0,\ldots,x_n)\ |\ \{ x_0,\ldots,x_n\} \mbox{ is an $m$-simplex  for some }m\leq n\},
\] 
equipped with the usual boundary homomorphism 
\[
d_n (x_0,\ldots,x_n):= \sum_{i=0}^n(-1)^i(x_0,\ldots, \widehat{x_i},\cdots,x_n).
\]
Thus $H_n(C_\bullet)\cong H_n(|Y|)$ for all $n\geq 0$, where $|Y|$ is the geometric realization of $Y$.
Let $L_\bullet$ be the subcomplex whose $n$-the term is $\Z[X_n]$ where $X_n$ is the set of (non-degenerate) ordered $n$-simplices:
\[
X_n=X_n(Y):=\{ (x_0,\ldots,x_n)\ |\ \{ x_0,\ldots,x_n\} \mbox{ is an $n$-simplex}\}.
\]
Let $\bar{C}_\bullet$ denote the quotient complex $C_\bullet/L_\bullet$. Thus, as an abelian group, $\bar{C}_n$ is isomorphic to  the free abelian group $D_n=D_n(Y)$ on the  set 
$E_n\setminus X_n$ of (ordered) degenerate simplices and $C_n=L_n\oplus D_n \cong L_n\oplus \bar{C}_n$ as abelian groups for all $n$.

If $Y$ is a simplicial complex and if $\alpha\in Y_m$ we will denote by $Y^\alpha$ the subcomplex defined by $Y^\alpha_n:=\{ \beta\in Y_n\ |\ \alpha\cup\beta\in Y_{n+m+1}\}$. 
%A path or loop in $Y^\alpha$ will also be called an $\alpha$-path or $\alpha$-loop in $Y$. Thus $a_0,a_1,\ldots,a_n$  ($a_i\in Y_0$) is an $\alpha$-path in $Y$ if $\{ \alpha_{i-1},\alpha_i\}\cup\alpha\in Y_{m+2}$ for $i=1,\ldots,n$.  

\subsection{Definition and basic properties}\label{sec:yabasic}
Let $A$ be a commutative ring.  Let $Y_0(A)$ be the set of unimodular rows $u\in A^2$ up to multiplication by an element of $A^\times$. Furthermore, for any $a\in A$, we will also denote the unimodular row $(a,1)\in Y_0(A)$ by $a$. For any $u\in A^\times$, the symbol $\infty$ will denote the class of $(u,0)=(1,0)$ in $Y_0(A)$.
% We also set $\infty:=(1,0)\in  Y_0(A)$.
 Let $Y:=Y(A)$ be the simplicial complex whose set $Y_n(A)$ of $n$-simplices consists of all sets $\{ [u_0],\ldots,[u_n]\}$, $[u_i]\in Y_0(A)$ with the property that the matrices
\[
M_{ij}:=
\left[
\begin{array}{c}
u_i\\
u_j\\
\end{array}
\right]
\]
lie in $\gl{2}{A}$. The group $\pgl{2}{A}$, and therefore also $\gl{2}{A}$, acts on $Y(A)$ by right multiplication. The associated complex $L_\bullet=L_\bullet(Y(A))$ is naturally a complex of $\pgl{2}{A}$-modules.

 Our main concern below is to understand the structure of $Z_2(L_\bullet)$ as a $\Z[\pgl{2}{A}]$-module.
% This is a $\gl{2}{A}$-module, and hence also a $\spl{2}{A}$-module. 
In particular,  wish to calculate the groups   $\rpb{A}:=(Z_2(L_\bullet))_{\pspl{2}{A}}$ and $\pb{A}:=(Z_2(L_\bullet))_{\pgl{2}{A}}$ (respectively, the \emph{refined pre-Bloch group} and the \emph{pre-Bloch group} of the ring $A$).

 Because of the exact sequence
 \[
1\to\pspl{2}{A}\to\pgl{2}{A}\to A^\times/(A^\times)^2\to 1,
\]
 the group $\rpb{A}$ is naturally an $\sgr{A}$-module (where $\sgr{A}:= \Z[A^\times/(A^\times)^2]$): for all $z\in Z_2(L_\bullet)$,  $\an{a}\cdot [z] = [z\cdot M]$ for any $M\in \gl{2}{A}$ with determinant $a$, and where 
$[z]$ denotes the class of $z$ in  $\rpb{A}$.  When $M\in\pgl{2}{A}$, $\an{\mathrm{det}(M)}$ will denote the resulting square class in $A^\times/(A^\times)^2$.  If $N$ is a (right) 
$\pgl{2}{A}$-module, then the group of coinvariants $M_{\pspl{2}{A}}$ is thus naturally an $\sgr{A}$-module via
\[
\an{\mathrm{det}(M)}\cdot\bar{x}=\overline{x\cdot M}
\]
for all $x\in N$, $M\in \pgl{2}{A}$.

Note  that $\{ 0,\infty,1\}\in Y_2(A)$ for any ring $A$. Virtually by definition, the group $\pgl{2}{A}$ acts transitively on $X_1(A)$, and hence on the set of $1$-simplices $Y_1(A)$; i.e., if $a=(a_1,a_2), b=(b_1,b_2)$ are the rows of an invertible matrix, then the resulting matrix 
\[
T_{a,b}:=\left[
\begin{array}{cc}
a_1&a_2\\
b_1&b_2\\
\end{array}
\right]\in \pgl{2}{A}
\]
has the property that $(\infty,0)\cdot T_{a,b}=(a,b)$.
 Furthemore, $\gl{2}{A}$ acts transitively on $X_2(A)$, and hence also on $Y_2(A)$:
Let $\{  a=(a_1,a_2),b=(b_1,b_2),c=(c_1,c_2)\}\in X_2(A)$. Let
\[
M_{a,b,c}:= 
\left[
\begin{array}{cc}
d(b,c)&0\\
0&d(c,a)\\
\end{array}
\right]\cdot T_{a,b}=
\left[
\begin{array}{cc}
a_1d(b,c)&a_2d(b,c)\\
b_1d(c,a)&b_2d(c,a)\\
\end{array}
\right]
\in \pgl{2}{A}
\]
where, whenever $x=(x_1,x_2), y=(y_1,y_2)$ are the rows of a $2\times 2$ invertible matrix, we define
\[
d(x,y):= \mathrm{det}
\left[
\begin{array}{cc}
x_1&x_2\\
y_1&y_2\\
\end{array}
\right].
\]
Then $(\infty,0,1)\cdot M_{a,b,c}=(a,b,c)$. 

Observe that for any $(a,b,c)\in X_2(A)$, we have $\an{\mathrm{det}(M_{a,b,c})}=\an{d(a,b)d(b,c)d(c,a)}\in A^\times/(A^\times)^2$.

Note that if $\{ a,b,c\}\in Y_2$ then
\[
T_{a,b}^{-1}= \left[
\begin{array}{cc}
b_2&-a_2\\
-b_1&a_1\\
\end{array}
\right]
\mbox{ and }
M_{a,b,c}^{-1}=T_{a,b}^{-1}\cdot \left[
\begin{array}{cc}
d(c,a)&0\\
0&d(b,c)\\
\end{array}
\right]\mbox{ in } \pgl{2}{A}.
\]

The following facts are easily verified:
\begin{lem}  \label{lem:cr}
\begin{enumerate}
\item If $\{ a,b\}\in Y_1$ and if $c\in Y^{\{ a,b\}}_0$ then $c\cdot T_{a,b}^{-1}=(d(c,b),d(a,c))\in Y_0$.
\item If $\{a,b,c\}\in Y_2$ and if $d\in Y^{\{ a,b\}}_0$ then 
\[
u(a,b,c,d):=d\cdot M_{a,b,c}^{-1}=\left( d(c,a)d(d,b),d(a,d)d(b,c)\right)=\frac{ d(a,c)d(b,d)}{d(a,d)d(b,c)}:=\an{d(a,b,c)}\in A^\times.
\]
\end{enumerate}
\end{lem}

Since also $d(q,p)=-d(p,q)$ for $\{ p,q\}\in Y_1$, we deduce:
\begin{cor}\label{cor:cr}
If $\{x,y\}\in Y_1$ and if $a,b,c\in Y^{\{ x,y\}}_0$ then 
\[
u(x,y,a,b)\cdot u(x,y,b,c)=u(x,y,a,c)   \mbox{ in }A^\times.
\]
\end{cor}

We can express certain topological properties of the  simplicial complex $Y(A)$ in direct ring-theoretic terms. For a ring $A$ we let $\wn{A}:=\{ a\in A\ |\ a(1-a)\in A^\times\}\subset A^\times$. 
\begin{lem}\label{lem:ya} Let $A$ be a ring.
\begin{enumerate}
\item Every $0$ -simplex of $Y(A)$ is contained in a $1$-simplex and every $1$-simplex is contained in a $2$-simplex.
\item $Y_3(A)\not=\emptyset$ if and only if $\wn{A}\not=\emptyset$ and, if so, every $2$-simplex is contained in a $3$-simplex.
%\item The simplicial complex $Y(A)$ satisfies the condition $(C)_1$ if and only if $\wn{A}$ generates $A^\times$ as a group.
\end{enumerate}
\end{lem}

\begin{proof} 
\begin{enumerate}
\item Every $0$-simplex is contained in a $1$-simplex by the definition of \emph{unimodular}. If $\{a,b\}$ is a $1$-simplex, then it is contained in the $2$-simplex $\{ a,b,c\}$ where 
$c=(1 1)\cdot T_{a,b}^{-1}$. 
\item Note that $Y^{\{ \infty,0\}}_0=\{ (u,1)\ |\ u\in A^\times\}=A^\times$. Now one easily verifies that if $u\in A^\times$, then $\{ \infty,0,1,u\}\in Y_3$ if and only if $1-u\in A^\times$; i.e., 
if and only if $u\in \wn{A}$. Conversely, if $\{ a,b,c,d\}$ is any $3$-simplex, then so is $\{\infty,0,1,u\}$ where $u=d\cdot M_{a,b,c}^{-1}$. 

Conversely,  if $u\in \wn{A}$ and if $\{ a,b,c\}\in Y_2(A)$, then $\{a,b,c,d\}\in Y_3(A)$ when $d=u\cdot M_{a,b,c}$. 
%\item Since $\gl{2}{A}$ acts transitively on edges of $Y(A)$, it follows that for any edge $\{ x,y\}\in Y_1(A)$, the simplicial complex $Y^{\{ x,y\}}$ is path-connected if and only if 
%$Y^{\{ \infty,0\}}$ is path-connected. Recall that $Y^{\{ \infty,0\}}_0=A^\times$. Furthermore, if $u,v\in Y^{\{ \infty,0\}}_0$, then $\{ u,v\}\in Y^{\{ \infty,0\}}_1$ if and only if 
%$\{ 1,u^{-1}v\}\in Y^{\{ \infty,0\}}_1$ (since multiplication by $\mathrm{diag}(u^{-1},1)$ is an automorphism of $Y^{\{ \infty,0\}}$) if and only if $\{ \infty,0,1,u^{-1}v\}\in Y_3(A)$ if and only if $u^{-1}v\in \wn{A}$. Thus an $\{ \infty,0\}$-path from $u$ to $v$ in $Y(A)$ is a sequence of units $u=u_0,u_1,\cdots, u_{n-1},u_n=v$ satisfying $w_i:=u_{i-1}^{-1}u_i\in \wn{A}$ for 
%$i=1,\ldots,n$. It follows that $u^{-1}v= w_1\cdots w_n$. 

%Conversely, given $u,v\in A^\times$, if there exist $w_1,\ldots,w_n\in\wn{A}$ such that $u^{-1}v=w_1w_2\cdots w_n$, then $u, uw_1,uw_1w_2,\cdots,uw_1\cdots w_n=v$ is an $\{ \infty,0\}$-path 
%in $Y(A)$ from $u$ to $v$. 
\end{enumerate}
\end{proof}

%\begin{rem} Similarly, it can also be easily seen that $Y(A)$ has property $(C)_0$ if and only if $A^\times$ generates $A$ as an additive group, but we will not need this property 
%below.
%\begin{proof} $Y(A)$ satisfies $(C)_0$ if and only if  Y^\infty$ is path-connected. Now $Y_0^\infty=A=\{ (a,1)\ |\ a\in A\}\subset $Y_0(A)$. If $a,b\in A$ then $\{ a,b\}\in Y^\infty_1$
%if and only if $\{\infty, a,b\}\in Y_2(A)$ if and only if $d(a,b)=a-b\in A^\times$.
%\end{proof}
%\end{rem}
%%%%%%%%%%%%%%%%%%%%%%%%%%%%%%%%%%%%%%%%%%%%%%%%%%%%%%%%%%%%%%%%%%%%%
\subsection{The connection with the low-dimensional homology of $\spl{2}{A}$}
%%%%%%%%%%%%%%%%%%%%%%%%%%%%%%%%%%%%%%%%%%%%%%%%%%%%%%%%%%%%%%%%%%%%%%%
Let $\LL_\bullet:=\LL_\bullet(A)$ be the following truncation of the complex $L_\bullet$:
\[
\LL_n(A)=\left\{
\begin{array}{ll}
L_n(A),& n\leq 2\\
Z_2(L_\bullet):=\ker{L_2(A)\to L_1(A)},& n=3\\
0,& n>3
\end{array}
\right.
\]
(where the map $\LL_3(A)=Z_2(L_\bullet)\to \LL_2(A)=L_2(A)$ is the inclusion homomorphism).
 Let $\epsilon:L_0(A)=Z[X_0]\to \Z$ be the homomorphism sending each $x\in  X_0$ to $1$. If we regard $\Z$ as a complex of (trivial) $\pgl{2}{A}$-modules concentrated in dimension $0$, then $\epsilon$ induces maps of complexes $\epsilon:\LL_\bullet \to \Z$ and $\epsilon:L_\bullet\to \Z$. It can be shown (\cite{hut:ge2}) that this map induces an isomorphism on homology for many rings $A$,
including when $A$ is a field or a local ring or when $A=\Z$. It follows that for such rings $A$, and for any group $G$ admitting a homomorphism to $\pgl{2}{A}$, there is a hyperhomology spectral sequence of the form 
\[
E^1_{p,q}=\ho{p}{G}{L_q}\Longrightarrow \ho{p+q}{G}{\LL_\bullet}\cong\ho{p+q}{G}{\Z}.
\]

Thus, when $G=\spl{2}{A}$, we have $E^1_{0,3}=\ho{0}{\spl{2}{A}}{Z_2(L_\bullet)}=\rpb{A}$ and $E^1_{0,2}=\ho{0}{\spl{2}{A}}{L_2}\cong\sgr{A}$ via $(a,b,c)\mapsto \an{d(a,b,c)}$.

 We let $\lambda_1$ be the map 
\[
d^1_{0,3}:\rpb{A}=E^1_{0,3}\to E^1_{0,2}=\sgr{A}
\]
and we let $\rpbker{A}:=\ker{\lambda_1}=E^2_{0,3}$. Now $E^2_{1,1}=0$ always (\cite[Corollary 4.4]{hut:bloch1}). We let $\lambda_2$ be the map $d^3_{0,3}:\rpbker{A}\to E^3_{2,0}$. Then $E^4_{0,3}=E^\infty_{0,3}=\ker{\lambda_2:\rpbker{A}\to E^3_{2,0}}=\rbl{A}$ is the \emph{(refined) Bloch group} of the ring $A$. Thus when $A$ is a field or a local ring or $\Z$, there is a natural surjective edge homomorphism $\ho{3}{\spl{2}{A}}{\Z}\to \rbl{A}\subset\rpb{A}$. (For more details about this spectral sequence, see \cite[Sections 4,5]{hut:bloch1}
%%%%%%%%%%%%%%%%%%%%%%%%%%%%%%%%%%%%%%%%%%%%%%%%%%%%%%%%%%%%%%%%%%%%%%%%%%%%%%%%
%%%%%%%%%%%%%%%%%%%%%%%%%%%%%%%%%%%%%%%%%%%%%%%%%%%%%%%%%%%%%%%%%%%%%%%%%%%%%
\subsection{ The space $Y(A)$ for local rings $A$}
%%%%%%%%%%%%%%%%%%%%%%%%%%%%%%%%%%%%%%%%%%%%%%%%%%%%%%%%%%
For a ring $A$ let $\Gamma(A)$ be the $1$-skeleton of the simplicial complex $Y(A)$. Thus $\Gamma(A)$ is a graph whose vertices are classes of unimodular rows and whose edges are pairs 
$\{ [u_0], [u_1]\}$ with the property that the matrix\\
 $T_{u_0,u_1}=\left[
\begin{array}{c}
u_0\\
u_1\\
\end{array}
\right]$ is invertible. Then, by definition, the simplicial complex $Y(A)$ is just the clique complex of the graph $\Gamma(A)$; i.e., it is the complex whose $r$-simplices are the vertices of $(r+1)$-cliques (complete subgraphs on $r+1$ vertices)  of $\Gamma(A)$.

Let $A$ be a local ring with maximal ideal ${M}$ and residue field $k=A/M$. Then $(a_1,a_2)$ is unimodular if and only if either $a_1\in A^\times$ or $a_2\in A^\times$. If $a_2\in A^\times$, 
then $(a_1,a_2)=(a,1)=a\in A$ where $a=a_1a_2^{-1}$. If $a_1\in A^\times$, then $(a_1,a_2)=(1,b)$ where $b=a_1^{-1}a_2\in A$. For convenience, we will use the notation $1/b$ for such classes $(1,b)$ of unimodular rows.  Note that $\infty:=1/0$ in this notation. Also,  when $a\in A^\times$ we have $1/a=(1,a)=(a^{-1},1)=a^{-1}$. It follows that $1/A\cap A=A^\times$ and the vertices of $\Gamma(A)$ 
are $A\cup 1/A=A\sqcup 1/M$. 

For $a,b\in A\subset \Gamma(A)_0$, $\mathrm{det}(T_{a,b})=a-b$ and hence $\{ a,b\}$ is an edge if and only if $a-b\in A^\times$ if and only if $a\not\equiv b\pmod{M}$.
For any $a\in A, m\in M$, $\mathrm{det}(T_{a,1/m})=am-1\in A^\times$ and hence $\{ a,1/m\}$ is always an edge. Finally, for 
any $m,n\in M$, $\mathrm{det}(T_{1/m,1/n})=n-m\in M$ and hence $\{ 1/m,1/n\}$ is never an edge of $\Gamma(A)$.

Suppose now that the residue field $k$ is finite of cardinality $q$.
% Fix a set of representatives $\{ x_0,x_1,\ldots,x_{q-1}\}\subset A$ for the cosets of $M$ in $A$.  Furthermore set $x_q:=\infty=1/0$. For $0\leq i\leq q$, let $T_i\subset %Y(A)_0=\Gamma(A)_0$ be defined by
For $i\in \projl{k}$ define $T_i\subset Y(A)_0$ by
\[
T_i=\left\{
\begin{array}{ll}
\mbox{ the corresponding coset of $M$ in $A$},& i\in k\\
1/M,& i=\infty.\\
\end{array}
\right.
\]  
Note that each $T_i$ is naturally in bijective correspondence with $M$ for all $i\in\projl{k}$. In particular, if $A$ is finite, then $|T_i|=|M|=|A|/q$ for all $i\in\projl{k}$.

We can describe the
 vertices, $\Gamma(A))_0$, and edges, $\Gamma(A)_1$, of the graph $\Gamma(A)$ as follows:\\
 $\Gamma(A)_0=\bigsqcup_{i\in\projl{k}}T_i$ and 
$\{ t,s\}\in \Gamma(A)_1$ if and only if $t\in T_{i}$ and $s\in T_{j}$ with $i\not=j$; in other words, $\Gamma(A)$ is the \emph{complete multipartite 
graph}, $K(\{ T_i\}_{i\in\projl{k}})$  on the collection of sets $\{ T_i\}_{i}$. It follows that $Y(A)=Y(\{ T_i\}_{i})$ is the clique complex of this complete multipartite 
graph. In detail, $\{ t_0,\ldots,t_r\}\subset \Gamma(A)_0$ is an $r$-simplex of $Y(A)$ if and only if $t_k\in T_{i_k}$, for $0\leq k\leq r$ with $i_k\not=i_{k'}$ when $k\not=k'$. 
To review:

\begin{prop}\label{prop:localya} Let $A$ be a local ring with maximal ideal $M$ and finite residue field $K$ of size $q$. Then $\Gamma(A)$ is the complete multipartite graph 
$K(\{ T_i\}_{i\in\projl{k}})$ where the $T_i$
for $i\in k$ is the corresponding  coset of $M$  in $A$ and $T_\infty:= 1/M\subset Y(A)_0$. Thus the simplicial complex $Y(A)$ is the clique complex of this graph.
\end{prop}

%\begin{rem}
%Observe that the index set $I\cup \{ \infty\}$ is naturally in bijective correspondence with $k\cup \{ \infty\}=\projl{k}$.
%\end{rem}

For example, for any prime $p$, $\Gamma(\Z/p^n)$ is the complete bipartite graph $K(\{T_i\}_{i=0}^p)$ on $p+1$ sets $T_i$, each of cardinality $|M|=p^{n-1}$; it is the complete 
multipartite graph 
\[
\Gamma(\Z/p^n)\cong K_{\underbrace{p^{n-1},\ldots,p^{n-1}}_{p+1}}
\]
and hence $Y(\Z/p^n)$ is the corresponding clique complex
\[
Y(\Z/p^n)\cong Y_{\underbrace{p^{n-1},\ldots,p^{n-1}}_{p+1}}.
\]

In Appendix \ref{app:clique} below we determine the homotopy type and the homology of spaces of the form $Y(T_0,\ldots,T_q)$.    By Proposition \ref{prop:mc2} and Theorem \ref{thm:spheresum}, we have:

\begin{prop}\label{prop:localhom} Let $A$ be a local ring with maximal ideal $M$ and a finite residue field $k=A/M$ of size $q$. 
For convenience, fix a bijection between $\projl{k}$ and $[0,q]:=\{ i\in \Z\ | \ 0\leq i\leq q\}$. For $i\in [0,q]$, let $T_i$  be as above.
For each $i\in [0,q]$, fix $x_i\in T_i$.
\begin{enumerate} 
\item For any $s=(s_0,\ldots,s_q)\in \prod_{i=0}^q (T_i\setminus\{ x_i\})$ let $Y_s$ be the simplicial complex\\
 $Y(\{ x_0,s_0\};\cdots ;\{ x_q,s_q\})$; i.e., the clique complex of the complete multipartite graph on this collection of  $q+1$ sets. Then $|Y_s|$ is homeomorphic to a $q$-sphere and the cycle
 \[
\gamma^s_q:=\sum_{y=(y_0,\ldots,y_q)\in\prod_{i=0}^q\{ x_i,s_i\}}(-1)^{\delta(y,s)}y\in  L_q(Y_s)\subset C_q(Y_s)
\]
(where $\delta(y,s):=\card{\{ i\ |\ y_i\not=s_i\}}=\card{\{ i\ |\ y_i=x_i\}}$) represents a generator of $H_q(|Y_s|)\cong \Z$.
\item The inclusion maps of simplicial complexes $Y_s\to Y(A)$, as $s$ varies over\\
 $S:=\prod_{i=0}^q (T_i\setminus\{ x_i\})$,  induce a map of simplicial complexes $\bigvee_s Y_s\to Y(A)$ 
(where the $Y_s$ are pointed with basepoint $x_0$) which induces
a homotopy equivalence of spaces
\[
\bigvee_{s\in S}|Y_s|\cong |Y(A)|.
\]
In particular, $H_n(|Y(A)|)=0$ for $n\not=0,q$ and $H_q(|Y(A)|)\cong \oplus_{s\in S} \Z$ with basis $\{ \gamma_q^s\ |\ s\in S=\prod_{i=0}^q (T_i\setminus\{ x_i\})\}$.
\end{enumerate}
\end{prop}
 Thus in this case $Y(A)$ is homotopy equivalent to a pointed sum of $q$-spheres and its $q$-dimensional homology is a free $\Z$-module with an explicit basis. For example, by Proposition \ref{prop:localhom}, $Y(\Z/p^n)$ has the homotopy type of a 
pointed sum of $(p^{n-1}-1)^{p+1}$ $p$-spheres. 

We will be particularly concerned below with the case when $A$ is a local ring 
with residue field of order $q=2$ or $q=3$, in which cases, as we will see below,  the corresponding homology groups $H_q(Y(A))$ play a role in the description of $\rpb{A}$.

\begin{rem} It follows from Proposition \ref{prop:localhom} that, for any local ring $A$, $\pi_1(Y(A),\infty)=0$. By \cite[Theorem 6.9]{hut:ge2}, $\pi(Y(A),\infty)$ is isomorphic to the quotient of the algebraic $K$-theory group $K_2(2,A)$ by the subgroup $C(2,A)$ generated by symbols $c(u,v)$. Thus we derive an alternative proof of the known fact that $K_2(2,A)$ is generated by symbols when $A$ is a local ring. 
\end{rem}

%%%%%%%%%%%%%%%%%%%%%%%%%%%%%%%%%%%%%%%%%%%%%%%%%%%%%%%%%%%%%%%%%%%%%%%%%%%%%%%
\section{Low-dimensional homology of simplicial complexes}\label{sec:lowdim}
%%%%%%%%%%%%%%%%%%%%%%%%%%%%%%%%%%%%%%%%%%%%%%%%%%%%%%%%%%%%%%%%%%%%%%%%%%%%%%%

%%%%%%%%%%%%%%%%%%%%%%%%%%%%%%%%%%%%%%%%%%%%%%%%%%%%%%%%%%%%%%%%%%%%%%%%%%%
%%%%%%%%%%%%%%%%%%%%%%%%%%%%%%%%%%%%%%%%%%%%%%%%%%%%%%%%%%%%%%%%%%%%%%%%%%%%%%
\subsection{Elementary results in low-dimensions}
\begin{lem}  For any simplicial complex $Y$ we have $H_1(\bar{C}_\bullet)=0$.
\end{lem}
\begin{proof} $\bar{C}_1$ is generated by the terms $(x,x)$, $x$ a $0$-simplex. Since $(x,x)=d_2(x,x,x)$, the results follows.
\end{proof}

\begin{lem}\label{lem:c2} Suppose that every $1$-simplex of  $Y$  is contained in a $2$-simplex. Then $H_2(\bar{C}_\bullet)=0$.
\end{lem}

\begin{proof}  $\bar{C}_2$ is generated by the terms $(x,x,x)$, $x$  a $0$-simplex and the terms\\
 $(x,x,y), (x,y,x), (y,x,x)$, $\{ x,y\}$ a $1$-simplex.  Now, for any edge $\{ x,y\}$,
$
d_3(x,x,x,y)=(x,x,y)-(x,x,x)$ and hence 
\begin{eqnarray}\label{eqn:1}
(x,x,y)\equiv (x,x,x)\mod{d_3(\bar{C}_3)}.
\end{eqnarray}
Similarly, $d_3(y,x,x,x)=(x,x,x)-(y,x,x)$ implies
\begin{eqnarray}\label{eqn:2}
(y,x,x)\equiv (x,x,x)\mod{d_3(\bar{C}_3)}.
\end{eqnarray}

By assumption, given any $1$-simplex $\{ x,y\}$, by hypothesis there exists a $0$-simplex $z$ such that $\{ x,y,z\}$ is a $2$-simplex. Then 
\begin{eqnarray*}
d_3(z,x,y,x)&=& (x,y,x)-(z,y,x)+(z,x,x)-(z,x,y)= (x,y,x)+(z,x,x) \mbox{ in } \bar{C}_2.
\end{eqnarray*}
and hence 
\begin{eqnarray}\label{eqn:3}
(x,y,x)\equiv-(z,x,x)\equiv - (x,x,x)\mod{d_3(\bar{C}_3)}.
\end{eqnarray}

By (\ref{eqn:1}), (\ref{eqn:2}) and (\ref{eqn:3}), $\bar{C}_2$  is generated modulo $d_3(\bar{C}_3)$ by the terms of the form $(x,x,x)$. Thus, if $Z\in \bar{C}_2$ is a cycle, we may assume $Z$ is of the form $\sum_x n_x(x,x,x)$. Hence $0=d_2(Z)=\sum_x n_x(x,x)$ and this forces $n_x=0$ for all $x$, and hence $Z=0$.
\end{proof}

From the long exact homology sequence associated to the short exact sequence of complexes $0\to L_\bullet\to C_\bullet\to \bar{C}_\bullet\to 0$ we deduce:
\begin{cor}\label{cor:exact} If every $1$-simplex of $Y$ is contained in a $2$-simplex then
 \[
H_1(L_\bullet)\cong H_1(C_\bullet)\cong H_1(|Y|)
\]
and there is an exact sequence 
\[
H_3(\bar{C}_\bullet)\to H_2(L_\bullet)\to H_2(C_\bullet)\cong H_2(|Y|)\to 0.
\]
\end{cor}

\begin{cor}\label{cor:ses} Let $Y$ be a simplicial complex with the property that every $1$-simplex is contained in a $2$-simplex. 
%Let $Z_2:=Z_2(L_\bullet)=\ker{d_2:L_2\to L_1}$. 
Then there is a short exact sequence
\[
0\to \left[d_3(D_3)\cap L_2\right]+d_3(L_3)\to Z_2(L_\bullet)\to H_2(|Y|)\to 0.
\]
\end{cor}

\begin{proof} Let $\delta$ denote the connecting homomorphism $H_3(\bar{C}_\bullet)\to H_2(L_\bullet)$.
By Corollary \ref{cor:exact}, there is a short exact sequence $0\to\image{\delta}\to Z_2(L_\bullet)/d_3(L_3)\to H_2(|Y|)\to 0$. 

%Let $\bar{Z}_3:= \ker{d_3:\bar{C}_3\to\bar{C}_2}$. 
Under the isomorphism $D_3\cong \bar{C}_3$, the subgroup $D_3^0:=\{ x\in D_3\ |\ d_3(x)\in L_2\}$ maps isomorphically to 
${Z}_3(\bar{C}_\bullet)$.  Now the definition of $\delta$ tells us that its image  is $\left(d_3(D_3^0)+d_3(L_3)\right)/d_3(L_3)$. Finally, observe that $d_3(D_3^0)=d_3(D_3)\cap L_2$.
\end{proof}

If $Y$ is $2$-dimensional then $Y_3=\emptyset$ and hence $L_3=0$. Thus we have 
\begin{cor}\label{cor:2dim} If $Y$ is $2$-dimensional then there is a short exact sequence
\[
0\to \left[d_3(D_3)\cap L_2\right]\to Z_2\to H_2(|Y|)\to 0.
\]
\end{cor}
Our immediate goal is to describe $d_3(D_3)\cap L_2$ explicitly. In fact we will derive an explicit  presentation for this group (see Proposition \ref{prop:pres}).

%%%%%%%%%%%%%%%%%%%%%%%%%%%%%%%%%%%%%%%%%%%%%%%%%%%%%%%%%%%%%%%%

\subsection{The subcomplex of partially symmetric chains}
For $\sigma\in S_{n+1}$ and $\{ x_1,\ldots,x_{n+1}\}\in Y_n$ we denote by $R_{\sigma}(x_1,\ldots,x_n)$ the element $(x_1,\ldots,x_n)-\sgn{\sigma}(x_{\sigma(1)},\ldots,x_{\sigma(n+1)})$ 
in $L_n$.  For $n\geq 0$, let $\Sigma_n\subset L_n$ be the subgroup generated by the elements $R_{\sigma}(x_1,\ldots,x_{n+1})$, $\{x_1,\ldots,x_{n+1}\}\in Y_n$, $\sigma\in \S_{n+1}$.

\begin{lem}\label{lem:gen}  For any $n\geq 1$, $\Sigma_n$ is generated by elements of the form
 \[
R_{(i\  i+1)}(x_1,\ldots,x_{n+1})=(x_1,\ldots, x_i,x_{i+1},\ldots ,x_{n+1})+(x_1,\ldots, x_{i+1},x_i\ldots ,x_{n+1})
\]
with $1\leq i\leq n$, $\{ x_1,\ldots,x_{n+1}\}\in Y_n$.
\end{lem}

\begin{proof}  This follows from the identity
\[
R_{\sigma\tau}(x_1,\ldots,x_{n+1})=R_{\tau}(x_1,\ldots,x_{n+1})+\sgn{\tau}R_{\sigma}(x_{\tau(1)},\ldots,x_{\tau(n+1)})
\]
together with the fact that the symmetric group  $S_{n+1}$ is generated by the set\\
 $\{ (i\ i+1)\ |\ 1\leq i\leq n\}$.
\end{proof}

If  $\{ x,y\}\in Y_1$ and if $\{ x,y\}\cup \{ a_1,\ldots,a_{n-1}\}\in Y_n$  we will let 
 \begin{eqnarray*}
R^{\{ x,y\}}_i(a_1,\ldots,a_{n-1}):&=& R_{(i\ i+1)}(a_1,\ldots,a_{i-1},x,y,a_i,\ldots,a_{n-1})\\
:&=& (a_1,\ldots,a_{i-1},x,y,a_i,\ldots,a_{n-1})+(a_1,\ldots,a_{i-1},y,x,a_i,\ldots,a_{n-1}).\\
\end{eqnarray*}
Thus $\Sigma_n$ is generated by the elements of the form $R^{\{ x,y\}}_i(a_1,\ldots,a_{n-1})$.

The following lemma can be verified directly from the definitions:
\begin{lem}\label{lem:d} If $\{x,y\}\cup\{ a_1,\ldots,a_{n-1}\}\in Y_n$ then for all $1\leq i\leq n$ we have
\[
d_n(R^{\{ x,y\}}_i(a_1,\ldots,a_{n-1}))=\sum_{j=1}^{i-1}(-1)^{j-1}R^{\{x,y\}}_{i-1}(a_1,\ldots,\hat{a_j}\ldots,a_{n-1})+\sum_{j=i}^{n-1}(-1)^{j-1}R^{\{x,y\}}_{i}(a_1,\ldots,\hat{a_j}\ldots,a_{n-1})
\]
(where the first term of the right is the empty sum  - i.e., equals $0$ --  if $i=1$, and the second term is $0$ if $i=n$).
\end{lem}

In view of Lemmas \ref{lem:gen} and \ref{lem:d} we immediately deduce:

\begin{cor}\label{cor:subcplx}  For all $n\geq 1$, $d_n(\Sigma_n)\subset \Sigma_{n-1}$ and hence $\Sigma_\bullet:=\{ \Sigma_n,d_n\}_{n\geq 0}$ is a subcomplex of $L_\bullet$.
\end{cor}

We may thus set $Z_n(\Sigma_\bullet):= \ker{d_n:\Sigma_n\to \Sigma_{n-1}}$ for $n\geq 1$. Of course, $Z_n(\Sigma_\bullet)\subset Z_n=Z_n(L_\bullet)$ for all $n\geq 1$.

\begin{prop}\label{prop:symm} 
 For all $n\geq 1$, $d_{n+1}(D_{n+1})\cap L_n\subset Z_n(\Sigma_\bullet)$.
\end{prop}

\begin{proof} 
%We begin by proving $d_{n+1}(D_{n+1})\cap L_n\subset Z_n(\Sigma_\bullet)$. 
Since $d_{n+1}(D_{n+1})\cap L_n\subset Z_n$ and since $Z_n(\Sigma_\bullet)=\Sigma_n\cap Z_n$, it is enough to prove that $d_{n+1}(D_{n+1})\cap L_n\subset \Sigma_n$. 

Let $\pi_L:C_n=L_n\oplus D_n\to L_n$ be the projection on the first factor. Clearly, given $w\in C_n$ then $w\in L_n$ if and only if $w=\pi_L(w)$. In particular, 
$d_{n+1}(D_{n+1})\cap L_n\subset \pi_L(d_{n+1}(D_{n+1})$. Hence it is enough to prove that $\pi_L(d_{n+1}(D_{n+1})\subset \Sigma_n$. In turn, it suffices to  prove this for a set of generators of $D_{n+1}$ which consists of elements $w=(x_1,\ldots,x_{n+2})$ where $x_i=x_j$ for some $i\not=j$. Now if three or more of the $x_i$ coincide then clearly 
$\pi_L(d_{n+1}(w))=0\in \Sigma_n$. So we may assume there is only one pair of indices $i<j$ such that $x_i=x_j$. In this case can write $w=(x_1,\ldots,x_i,\ldots, x_{j-1},x_i,x_j\ldots,x_{n+1})$ with $\{x_1,\ldots,x_{n+1}\}\in Y_n$ and hence 
\begin{eqnarray*}
\pi_L(d_{n+1}(w))&=& (-1)^{i-1}(x_1,\ldots,\hat{x_i},\ldots, x_{j-1},x_i,x_j\ldots ,x_{n+1})\\
&&+ (-1)^{j-1}(x_1,\ldots, x_i,\ldots, x_{j-1},x_j,\ldots,x_{n+1})\\
&=&(-1)^{j-1}\left[(x_1,\ldots, x_{n+1})+(-1)^{i+j}(x_1,\ldots,x_i,x_{i+1},\ldots, x_{j-1},x_i,x_j,\ldots x_{n+1}) \right]\\
&=&(-1)^{j-1}R_\sigma(x_1,\ldots,x_{n+1})\in \Sigma_n 
\end{eqnarray*}
as required, where $\sigma=(i\ i+1\cdots j-1)$ and hence $\sgn{\sigma}=(-1)^{(j-1)-i}=-(-1)^{i+j}$. 

%We prove the converse inclusion $Z_n(\Sigma_\bullet)\subset d_{n+1}(D_{n+1})$. By Corollary \ref{cor:decomp}, it is enough to prove this for $Z_n^{\{ x,y\}}$ for any given 
%$\{ x,y\}\in Y_1$:  Let $w\in \Sigma_n^{\{ x,y\}}$ with $d_n(w)=0$. Now consider $w':=(x,w)\in D_{n+1}$. Then $d_{n+1}(w')=w-(x,d_n(w))=w$ as required.
\end{proof}

%%%%%%%%%%%%%%%%%%%%%%%%%%%%%%%%%%%%%%%%%%%%%%%%%%%%%%%%%%%%%%%%%%%%%

\subsection{$2$-dimensional complexes}
%%%%%%%%%%%%%%%%%%%%%%%%%%%%%%%%%%%%%%%%%%%%%%%%%%%%%%%%%%%%%%%%%%%%%%%%%%%%%%
%\subsubsection{Generators and Relations for $Z_2(\Sigma_\bullet)$}

We introduce some families of elements in $Z_2(\Sigma_\bullet)$:  First observe that whenever $\{ x,y,a\}\in Y_2$ we have \[
d_2(R_1^{\{ x,y\}}(a))=d_2(R_2^{\{ x,y\}}(a))=R_1^{\{ x,y\}}:= (x,y)+(y,x).
\]
Whenever $\{x,y,a\},\{x,y,b\}\in Y_2$, we have the following elements of  $Z_2(\Sigma_\bullet)$:
\begin{eqnarray*}
D_1^{\{ x,y\}}(a,b)&:=& R_1^{\{ x,y\}}(b)-R_1^{\{ x,y\}}(a)\\
D_2^{\{ x,y\}}(a,b)&:=& R_2^{\{ x,y\}}(b)-R_2^{\{ x,y\}}(a).\\
\end{eqnarray*}
Furthermore, for any $\{ x,y,a\}\in Y_2$, let $C^{\{ x,y\}}(a):=R_1^{\{ x,y\}}(a)-R_2^{\{ x,y\}}(a)\in Z_2(\Sigma_\bullet)$. 

%Given a $1$-simplex $\{ x,y\}$, let 
%\[
%C_{x,y}:= (x,y)+(y,x)\in L_1.
%\]

%If $\{ x,y,z\}\in Y_2$, let $L_{x,y}(z):=(z,C_{x,y}):= (z,x,y)+(z,y,x)\in L_2$ and let $R_{x,y}(z):=(C_{x,y},z)$. 

%Let $K_\ell$ be the subgroup of $L_2$ generated by the union of the set $\{ L_{x,y}(z)\ |\ \{ x,y,z\}\in X_2\}$ and let $K_r$ be the subgroup generated by 
%$\{ R_{x,y}(z)\ |\ \{ x,y,z\}\in X_2\}$. Let $K=K_\ell+K_r$. 

%Given $\{ x,y\}\in Y_1$, let $S_{x,y}$ denote the set $\{ z \ |\ \{x,y,z\}\in Y_1\}$. For $\{ x,y\}\in Y_1$ and $a,b\in S_{x,y}$, let $L_{x,y}(a,b):= L_{x,y}(a)-L_{x,y}(b)$,
%$R_{x,y}(a,b):=R_{x,y}(a)-R_{x,y}(b)$, $C_{x,y}(a):=R_{x,y}(a)-L_{x,y}(b)$.

%Let $Z_c\subset K $ be the subgroup of $Z_2$ generated by the cycles $\{ C_{x,y}(a)\ |\ \{ x,y,a\}\in Y_2\}$. \\
%Let $Z_\ell\subset K_\ell$ be the subgroup of $Z_2$ generated by the cycles $\{ (L_{x,y}(a,b)\ | \{ x,y\}\in Y_1, a,b\in S_{x,y}\}$. 
%Let $Z_r\subset K_r$ be the subgroup of $Z_2$ generated by the cycles $\{ R_{x,y}(a,b)\ | \{ x,y\}\in Y_1, a,b\in S_{x,y}\}$.\\
%Let $Z:=Z_c+Z_\ell+Z_r$.

Now $\Sigma_2=\Sigma_2^1+\Sigma_2^2$ where $\Sigma_2^i$ is the group generated by the terms $R_i^{\{ x,y\}}(a)$ for $\{ x,y,a\}\in Y_2$ and $i=1,2$.

For $i=1,2$, let $Z^i\subset \Sigma_2^i$ be the subgroup generated by the terms $D_i^{\{ x,y\}}(a,b)$ and let $Z^c$ be the subgroup generated by the terms $C^{\{ x,y\}}(a)$.

\begin{lem}\label{lem:z}
 Let $Y$ be a simplicial complex. Then $Z_2(\Sigma_\bullet)$ is generated by the elements $D_i^{\{ x,y\}}(a,b)$, $i=1,2$ and $C^{\{ x,y\}}(a)$: i.e., $Z_2(\Sigma_\bullet)=Z^1+Z^2+Z^c$.
\end{lem}

\begin{proof} Let $W_1\in \Sigma_2^1$. By definition, we have 
\[
W_1=\sum_{\{ x,y\}}\sum_z n_{x,y,z}R_1^{\{ x,y\}}(z),
\]
where the first sum is over all $1$-simplices $\{ x,y\}$ which are contained in a $2$-simplex and the second is over the set $\{ z\ |\ \{ x,y,z\} \mbox{ is a $2$-simplex}\}$.
For each such $1$-simplex $\{ x,y\}$, fix a $0$-simplex $z_{x,y}$ such that $\{ x,y,z_{x,y}\}$ is a $2$-simplex. For each $\{ x,y\}$, let $N_{x,y}:=
\sum_z n_{x,y,z}$. Then
\begin{eqnarray*}
W_1&=& \sum_{\{ x,y\}}N_{x,y}R_1^{\{ x,y\}}(z_{x,y})+\sum_{\{ x,y\}}\sum_z n_{x,y,z}\left(R_1^{\{ x,y\}}(z)-R_1^{\{ x,y\}}(z_{x,y})\right)\\
&=& \sum_{\{ x,y\}}N_{x,y}R_1^{\{ x,y\}}(z_{x,y}) +z_1
\end{eqnarray*}
where $z_1=\sum_{\{ x,y\}}\sum_z n_{x,y,z}D_1^{\{ x,y\}}(z_{x,y},z)\in Z^1$.

Similarly, if $W_2=\sum_{\{ x,y\}}\sum_z m_{x,y,z}R_2^{\{ x,y\}}(z)\in \Sigma_2^2$ then we can write 
\[
W_2= \sum_{\{ x,y\}}M_{x,y}R_2^{\{ x,y\}}(z_{x,y})+z_2
\]
where $z_2\in Z^2$.  Hence $W:=W_1+W_2\in \Sigma$ can be written 
\[
 \sum_{\{ x,y\}}N_{x,y}R_1^{\{ x,y\}}(z_{x,y})+ \sum_{\{ x,y\}}M_{x,y}R_2^{\{ x,y\}}(z_{x,y})+z_1+z_2
\]
with $z_i\in Z^i$.

Now suppose $W\in Z_2(\Sigma_\bullet)$. Then
 \begin{eqnarray*}
0=d_2(W)&=&  \sum_{\{ x,y\}}N_{x,y}R_1^{\{ x,y\}}+ \sum_{\{ x,y\}}M_{x,y}R_1^{\{ x,y\}}\\
&=& \sum_{\{ x,y\}}(N_{x,y}+M_{x,y})R_1^{\{ x,y\}},
\end{eqnarray*}
which implies that $M_{x,y}=-N_{x,y}$ for all $\{ x,y\}$, since the set $\{ R_1^{\{ x,y\}}\}$ is linearly independent. Thus
\[
W= \sum_{\{ x,y\}}N_{x,y}\left(R_1^{\{ x,y\}}(z_{x,y})-R_2^{\{ x,y\}}(z_{x,y}))\right)+z_1+z_2=z_c+z_1+z_2
\]
with $z_c= \sum_{\{ x,y\}}N_{x,y}C^{\{ x,y\}}(z_{x,y})\in Z^c$, as required.
\end{proof}

\begin{lem}\label{lem:d3} Let $Y$ be a simplicial complex. Then $d_3(D_3)\cap L_2=Z_2(\Sigma_\bullet)$.
\end{lem}

\begin{proof}
Consider any one of the generators $C$  of the groups $Z^1,Z^2$ or $Z^c$ listed above. Then there exists some $\{ x,y\}\in Y_1$ such that $C=R_i^{\{ x,y\}}(a)$ or $C^{\{ x,y\}}(a)$. Then 
$E:=(x,C)\in D_3$ and $d_3(E)=C-(x,d_2(C))=C$ since $d_2(C)=0$. Thus $C\in d_3(D_3)$ as required. 
\end{proof}
Combining Lemma \ref{lem:d3} with Corollary \ref{cor:ses} we obtain:
\begin{prop} \label{prop:z2}
Let $Y$ be a simplicial complex for which every $1$-simplex is contained in a $2$-simplex. There is a short exact sequence
\[
0\to Z_2(\Sigma_\bullet)+d_3(L_3)\to Z_2(L_\bullet)\to H_2(|Y|)\to 0.
\]
\end{prop}

If $Y$ is $2$-dimensional, then $L_3(Y)=0$ and $Z_2(L_\bullet)=H_2(L_\bullet)$. We deduce:
\begin{cor} \label{cor:2dim}
If the simplicial complex $Y$ is $2$-dimensional and if each $1$-simplex is contained in a $2$-simplex then there is a short exact sequence 
\[
0\to Z_2(\Sigma_\bullet)\to Z_2(L_\bullet)\to H_2(|Y|)\to 0.
\]
If furthermore $H_2(|Y|)=0$ then $Z_2(\Sigma_\bullet)=Z_2(L_\bullet)$.
\end{cor}

%Given a $1$-simplex $\{x,y\}$, let $S_{x,y}:=\{ a\ | \{x,y,a\} \mbox{ is a $2$-simplex}\}$.  Now given a $1$-simplex $\{ x,y\}$, for any $a\in S_{x,y}$ let $L_{x,y}(a):=(a,C_{x,y})$ and let %$R_{x,y}(a):=(C_{x,y},a)$.
Given $\{ x,y\}\in Y_1$, for $i=1,2$,  let $\Sigma_2^i(x,y)$ denote the (free abelian) ) subgroup of $L_2$ generated by $\{ R_i^{\{ x,y\}}(a)\ |\   \{ x,y,a\}\in Y_2\}$. Then clearly 
$
\Sigma_2^i=\bigoplus_{\{ x,y\}\in Y_1}\Sigma_2^i(x,y)\ \mbox{ for }i=1,2.
$

  For any $2$-simplex $\{ x,y,z\}$ let $\Sigma(\{ x,y,z\}):= \sum_{\sigma}(\sigma x,\sigma y, \sigma z)\in L_2$, where $\sigma$ runs over the six permutations of $\{ x,y,z\}$. Observe that
\[
\Sigma(\{ x,y,z\})= R_1^{\{ x,y\}}(z)+R_1^{\{ x,z\}}(y)+R_1^{\{ y,z\}}(x) =  R_2^{\{ x,y\}}(z)+R_2^{\{ x,z\}}(y)+R_2^{\{ y,z\}}(x)\mbox{ in }\Sigma_2^1\cap\Sigma_2^2.
\]

 Recall that for $\{ x,y\}\in Y_1$,  $Y^{\{ x,y\}}_0=\{ z\in Y_0\ |\ \{ x,y,z\}\in Y_2\}$.

\begin{lem}\label{lem:K} We have 
\[
\Sigma_2^1\cap\Sigma_2^2 =\bigoplus_{\{ x,y,z\}\in Y_2}\Z\cdot \Sigma(\{ x,y,z\}).
\]
\end{lem} 

\begin{proof}  Let $W=\sum_{(x,y,z)\in Y_2}n_{x,y,z}(x,y,z)\in \Sigma_2^1\cap\Sigma_2^2 $.  Since $W\in \Sigma_2^1$ it can be written 
\[
W=\sum_{\{ x,y\}\in Y_1}\sum_{z\in {Y^{\{ x,y\}}_0} }m_{\{x,y\},z}R_{x,y}(z)=\sum_{\{ x,y\}\in Y_1}\sum_{z\in Y^{\{ x,y\}}_0}m_{\{x,y\},z}\left( (x,y,z)+(y,x,z)\right)
\]
Similarly, since $W\in \Sigma_2^2$ we have
\[
W=\sum_{\{ y,z\}\in Y_1}\sum_{x\in Y^{\{ y,z\}}_0}t_{x,\{ y,z\}}\left( (x,y,z)+(x,z,y)\right).
\]
Thus given $\{x,y,z\}\in Y_2$ we must have $n_{x,y,z}=m_{\{ x,y\},z}=m_{\{ y,x\},z}=n_{y,x,z}$ and $n_{x,y,z}=t_{x,\{y,z\}}=t_{x,\{z,y\}}=n_{x,z,y}$. Thus 
$n_{x,y,z}=n_{y,x,z}=t_{y, \{ x,z\}}=t_{y, \{ z,x\}}=n_{y,z,x}$ etc and we deduce that $n_{x,y,z}=n_{\sigma x,\sigma y,\sigma z}$ for all permutations $\sigma$ of $\{ x,y,z\}$, as required.
\end{proof}

%Given $\{x,y\}\in Y_1$ and $a,b\in S_{x,y}$, let $L_{x,y}(a,b):=L_{x,y}(a)-L_{x,y}(b)\in Z_\ell$ and let $R_{x,y}(a,b):= R_{x,y}(a)-R_{x,y}(b)\in Z_r$. Let $C_{x,y}(a):=R_{x,y}(a)-%L_{x,y}(a)\in Z_c$.

We will need the following well-known fact:
\begin{lem}\label{lem:exact} For any set nonempty set $S$, the complex
\[
\xymatrix{
\cdots\ar^{d_{n+1}}[r]&\Z[S^{\times n}]\ar^{d_n}[r]&\Z[S^{\times (n-1)}]\ar[r]&\cdots\ar[r]&\Z[S]\ar^{d_1}[r]&\Z\ar[r]&0,
}
\]
where $d_n(s_1,\ldots,s_n):=\sum_{i=1}^n(-1)^{i+1}(s_1,\ldots,\widehat{s_i},\ldots,s_n)$, is exact.
\end{lem}
\begin{proof} Fix $x\in S$ and, for all $n\geq 1$, let $\sigma_x:\Z[S^{\times (n-1)}]\to\Z[S^{\times n}]$ be the homomorphism $\sigma_x(s_1,\ldots,s_{n-1}):=(x,s_1,\ldots,s_n)$. Then for 
all $w\in \Z[[S^{\times (n-1)}]$, $d_n(\sigma_x(w))=w-\sigma_x(d_{n-1}(w))$ . Thus if $d_{n-1}(w)=0$, then $w=d_n(\sigma_x(w))$.
\end{proof}
\begin{prop}\label{prop:pres}
The abelian group $Z_2(\Sigma_\bullet)=d_3(D_3)\cap L_2$ has the following presentation: 

It is generated by the elements  $D_1^{\{x,y\}}(a,b)$, $D_2^{\{x,y\}}(a,b)$, $C^{\{ x,y\}}(a)$ with $\{ x,y\}\in Y_1$, $a,b\in Y^{\{ x,y\}}_0$, subject to the relations
\begin{enumerate}
\item $D_1^{\{ x,y\}}(a,b)+D_1^{\{ x,y\}}(b,c)=D_1^{\{ x,y\}}(a,c)$.
\item  $D_2^{\{ x,y\}}(a,b)+D_2^{\{ x,y\}}(b,c)=D_2^{\{ x,y\}}(a,c)$.
\item $C^{\{ x,y\}}(b)-C^{\{ x,y\}}(a)=D_1^{\{ x,y\}}(a,b)-D_2^{\{ x,y\}}(a,b)$.
\item For all $\{ x,y,z\}\in Y_2$, $C^{\{ x,y\}}(z)+C^{\{ x,z\}}(y)+C^{\{ y,z\}}(x)=0$.
\end{enumerate}
In (1), (2) we allow $a=b$, $b=c$ or $a=c$ and define $D_1^{\{ x,y\}}(z,z):=0$.
\end{prop}

\begin{proof}  Let $\widetilde{\Sigma_2}:=\Sigma_2^1\oplus \Sigma_2^2$ with generators $\tilde{R}^{\{ x,y\}}_1(a):=(R^{\{ x,y\}}_1(a),0)$ and 
$\tilde{R}^{\{ x,y\}}_2(a):=(0,R^{\{ x,y\}}_2(a))$. Let $\pi$ be the natural surjection $\widetilde{\Sigma_2}\to \Sigma_2=\Sigma_2^1+\Sigma_2^2$. 
By Lemma \ref{lem:K} 
\[
K:=\ker{\pi}=\oplus_{\{x,y,z\}\in X_2}\Z\cdot \tilde{C}(\{ x,y,z\})
\]
where $\tilde{C}(\{ x,y,z\}):=(\Sigma(\{ x,y,z\}),-\Sigma(\{ x,y,z\}))\in \widetilde{\Sigma_2}$.

 Let $\tilde{Z}^1=Z^1\oplus \{ 0\}\subset \widetilde{\Sigma_2}$
%. So $\tilde{Z}_r=\oplus \tilde{Z}_r^{x,y}$ where $\tilde{Z}_r^{x,y}:=Z_r\oplus \{ 0\}$ 
with generators $\tilde{D}_1^{\{ x,y\}}(a,b):= (D_1^{\{ x,y\}}(a,b),0)$. Similarly 
 Let $\tilde{Z}^2=\{ 0\} \oplus Z^2\subset \widetilde{\Sigma_2}$
%. So $\tilde{Z}_\ell=\oplus \tilde{Z}_\ell^{x,y}$ where $\tilde{Z}_\ell^{x,y}:=\{ 0\}\oplus Z_\ell^{x,y}$
 with generators $\tilde{D}_2^{\{ x,y\}}(a,b):= (0, D_2^{\{ x,y\}}(a,b))$. Let $\tilde{Z}^c$ denote the subgroup of $\widetilde{\Sigma_2}$ generated by 
the elements $\tilde{C}^{\{ x,y\}}(a):=(R_1^{\{ x,y\}}(a),-R_2^{\{ x,y\}}(a))$ for any $\{ x,y\}\in Y_1$, $a\in Y^{\{ x,y\}}_0$.
 Note that $\tilde{C}(\{ x,y,z\})= \tilde{C}^{\{ x,y\}}(z)+\tilde{C}^{\{ x,z\}}(y)+\tilde{C}^{\{ y,z\}}(x)$. Since $K\subset \tilde{Z}^c$ and $\pi(\tilde{C}^{\{ x,y\}}(a))={C}^{\{ x,y\}}(a)$, $\pi$ induces 
an exact sequence $0\to K\to \tilde{Z}^c\to Z^c\to 0$. Furthermore, if we let $\tilde{Z}_2(\Sigma_\bullet):=\tilde{Z}^1+\tilde{Z}^2+\tilde{Z}^c\subset \widetilde{\Sigma_2}$ then we have an exact sequence
\[
\xymatrix{
0\ar[r]&K\ar[r]& \tilde{Z}_2(\Sigma_\bullet)\ar^-{\pi}[r]&Z_2(\Sigma_\bullet)\ar[r]&0\\
}
\]

Thus it is enough to prove that $\tilde{Z}_2(\Sigma_\bullet)$ is generated by the elements  $\tilde{D}_1^{\{ x,y\}}(a,b)$, $\tilde{D}_2^{\{ x,y\}}(a,b)$, $\tilde{C}^{\{ x,y\}}(a)$ satisfying (the analogue of) relations (1), (2), (3).

 First observe that, for $i=1,2$, $\tilde{Z}^i=\oplus_{\{ x,y\}\in Y_1}\tilde{Z}^i(x,y)$ where $\tilde{Z}^i(x,y) \subset \widetilde{\Sigma_2}(x,y):= \Sigma_2^1(x,y)\oplus \Sigma_2^2(x,y)$ is the group generated by the set 
 $\{\tilde{D}_i^{\{ x,y\}}(a,b)\ |\ a,b\in Y^{\{ x,y\}}_0\}$. 
%Analogously, $\tilde{Z}_r=\oplus_{\{ x,y\}\in Y_1}\tilde{Z}_r^{x,y}$ where $\tilde{Z}^{x,y}_r\subset \tilde{K}^{x,y}_r$ is the group generated by the set 
% $\{\tilde{R}_{x,y}(a,b)\ |\ a,b\in \}$.
 For $\{ x,y\}\in Y_1$, let $\tilde{Z}^c(x,y)\subset \widetilde{\Sigma_2}(x,y)$ denote the group generated by
$\{ \tilde{C}^{\{ x,y\}}(a)\ |  a\in Y^{\{ x,y\}}_0\}$.

Clearly $\tilde{Z}_2(\Sigma)=\oplus \tilde{Z}_2(\Sigma)(x,y)$ where $\tilde{Z}_2(\Sigma)(x,y):= \tilde{Z}^1(x,y)+\tilde{Z}^2(x,y)+\tilde{Z}^c(x,y)$ and it is enough the prove the given presentation for each fixed
 $\{x,y\}\in Y_1$.
So we fix $\{ x,y\}\in Y_1$ (and we will drop most of the references to $\{x,y\}$ from our notation below). Let $S:= Y^{\{ x,y\}}_0$.

Now, for $i=1,2$,  $\widetilde{\Sigma}_2^i(x,y)\cong \Z[S]$ via $\tilde{R}_i(a)\leftrightarrow a$.  
%For any $a\in S$ we have $d_2(L(a))=d_2(R(a))=C_{x,y}\in L_1$. 
Under these isomorphisms, $\tilde{Z}^i(x,y)$  corresponds to the subgroup $\Z[S]^0:=\ker{\epsilon:\Z[S]\to \Z}$, where $\epsilon(a)=1$ for all $a$: i.e., it is generated by 
the elements  $\tilde{D}_i(a,b)=\tilde{R}_i(b)-\tilde{R}_i(a)\leftrightarrow (b)-(a)=d_2(a,b)\in \Z[S]^0$.

However,
by Lemma \ref{lem:exact}, the sequence
\[
\Z[S\times S\times S]\to \Z[S\times S]\to \Z[S]^0\to 0
\]
is exact. From this it follows that $\tilde{Z}^i(x,y)$ is generated by the $\tilde{D}_i(a,b)$ subject to relations (1) when $i=1$ or (2) when $i=2$.
%, similarly, that $\tilde{Z}^{x,y}_r$ is generated by the $\tilde{R}(a,b)$ subject to relations (2). 

Finally, fix $a\in Y_0$ and observe that $\Z\cdot \tilde{C}(a)\cap (\tilde{Z}^1(x,y)\oplus \tilde{Z}^2(x,y))=\{ 0\}$
 (since $\Z\cdot\tilde{R}_i(a)\cap \tilde{Z}^i(x,y)=\{ 0\}$ ).
Since $\tilde{C}(b)=\tilde{C}(a)+\tilde{D}_1(a,b)-\tilde{D}_2(a,b)$ for all $b\in S$ it follows that $\tilde{Z}_2(\Sigma)(x,y)=\tilde{Z}^1(x,y)\oplus \tilde{Z}^2(x,y)\oplus \Z\cdot \tilde{C}(a)$ and the required presentation follows immediately.
\end{proof}

\begin{cor}\label{cor:z2pres}
$Z_2(\Sigma_\bullet)=Z^c\oplus Z^1=Z^c\oplus Z^2$  where $Z^c$ has the presentation:

Generators $C^{\{ x,y\}}(a)$ , $\{ x,y,a\}\in Y_2$\\
Relations:  For all $\{ x,y,z\}\in Y_2$, $C^{\{ x,y\}}(z)+C^{\{ x,z\}}(y)+C^{\{ y,z\}}(x)=0$, and

$Z^i$ has the presentation

Generators: $D_i^{\{ x,y\}}(a,b)$   with $\{ x,y\}\in Y_1$, $a,b\in Y^{\{ x,y\}}_0$\\
 Relations: $D_i^{\{ x,y\}}(a,b)+D_i^{\{ x,y\}}(b,c)=D_i^{\{ x,y\}}(a,c)$\\
%(resp. $R_{x,y}(a,b)+R_{x,y}(b,c)=R_{x,y}(a,c)$ ).
\end{cor}
\begin{proof}
Observe that in relations (1)-(4) of Proposition \ref{prop:pres},  (3) can be removed and instead taken as a definition of the terms  $D_2^{\{ x,y\}}(a,b)$. With this definition (2) follows from (1), so that (2) can be removed.

 Likewise, if we remove (3) and instead use it to define  the terms  $D_1^{\{ x,y\}}(a)$, then relation (1) becomes redundant.
\end{proof}

%%%%%%%%%%%%%%%%%%%%%%%%%%%%%%%%%%%%%%%%%%%%%%%%%%%%%%%%%%%%%%%%%%%%
%%%%%%%%%%%%%%%%%%%%%%%%%%%%%%%%%%%%%%%%%%%%%%%%%%%%%%%%%%%%%%%%%%%

\section{Calculating $\rpb{A}$ for rings $A$ such that $\wn{A}=\emptyset$}\label{sec:wna=0}
%%%%%%%%%%%%%%%%%%%%%%%%%%%%%%%%%%%%%%%%%%%%%%%%%%%%%%%%%%%%%%

We begin by observing that $\wn{\Z/2}=\emptyset$ and thus $\wn{A}=\emptyset$ for any ring $A$ which has an ideal of index $2$.

For convenience, we introduce the following notations for a ring $A$:
\[
\crpb{A}:= Z_2(\Sigma_\bullet)_{\pspl{2}{A}},\quad \cpb{A}:= Z_2(\Sigma_\bullet)_{\pgl{2}{A}}.
\]

By 
%Proposition \ref{prop:z2} 
Corollary \ref{cor:2dim} and Lemma \ref{lem:ya} (2), taking $\spl{2}{A}$-coinvariants, we have 
\begin{prop}\label{prop:exacta} Let $A$ be a commutative ring such that $\wn{A}=\emptyset$.  Then there is an  exact sequences of $\sgr{A}$-modules
\[
\crpb{A}
\to \rpb{A}\to H_2(|Y(A)|)_{\spl{2}{A}}\to 0
\]
and an exact sequence of abelian groups
\[
\cpb{A} \to \pb{A}\to H_2(|Y(A)|)_{\gl{2}{A}}\to 0
\]
Furthermore, if $H_2(|Y(A)|)=0$  then $\rpb{A}=\crpb{A}$ and $\pb{A}=\cpb{A}$.
\end{prop}

%\begin{rem} Recall that $L_3(A)=0\iff Y_3(A)=\emptyset\iff \wn{A}=\emptyset$\\
%  (where $\wn{A}=\{ u\in A^\times\ |\ 1-u\in A^\times\}$).
%\end{rem}

Recall that $Z_2(\Sigma_\bullet)=Z^c\oplus Z^1$. For the generators, we naturally have $C^{\{ x,y\}}(a)\cdot g=C^{\{ xg,yg\}}(ag)$ for all $\{x,y,a\}\in X_2$, $g\in \pgl{2}{A}$ and 
$D_i^{\{ x,y\}}(a,b)\cdot g=D_i^{\{ xg,yg\}}(ag,bg)$ when $\{ x,y\}\in Y_1$, $a,b\in Y^{\{ x,y\}}_0$. Thus this is a decomposition of $\pgl{2}{A}$-modules and thus 
\[
Z_2(\Sigma_\bullet)_G=(Z^c\oplus Z^1)_G=(Z^c)_G\oplus (Z^1)_G
\]
for $G=\pgl{2}{A}$ or $\pspl{2}{A}$.

We will denote the module $(Z^c)_{\pspl{2}{A}}\subset \crpb{A}$ by $\ccconstmod{A}$ and we will denote the submodules $(Z^i)_{\pspl{2}{A}}$ by $\cks{i}{A}$ for $i=1,2$. Thus
$\crpb{A}=\ccconstmod{A}\oplus \cks{1}{A}=\ccconstmod{A}\oplus \cks{2}{A}$.

%%%%%%%%%%%%%%%%%%%%%%%%%%%%%%%%%%%%%%%%%%%%%%%%%%%%%%%%%%%%%%%%%%%%%%%%%%%%%%%%
\subsection{A presentation of $\crpb{A}$}
%%%%%%%%%%%%%%%%%%%%%%%%%%%%%%%%%%%%%%%%%%%%%%%%%%%%%%%%%%%%%%%%%%%%%%%%%%%%%%%%%%%
Let $\cconst{A}:=C^{\{ \infty,0\}}(1)$ in $Z^c$. In this section we prove:

\begin{prop}\label{prop:crpbpres}  Let $A$ be a commutative ring. Then the $\sgr{A}$-module $\crpb{A}$ has the following presentation:

\begin{enumerate}
\item[Generators:] The elements $\suss{i}{u}$, $u\in A^\times$, $i\in \{ 1,2\}$ and the constant $\cconst{A}$
\item[Relations:]
\begin{enumerate}
\item $\an{-1}\cconst{A}=\cconst{A}$ and $\an{-1}\suss{i}{u}=\suss{i}{u^{-1}}$,
\item $3\cconst{A}=0$,
\item For $i\in \{ 1,2\}$ and $u,v\in A^\times$, $\suss{i}{uv}=\an{u}\suss{i}{v}+\suss{i}{u}$.
\item For all $u\in A^\times$, $\pf{u}\cconst{A}=\suss{1}{u}-\suss{2}{u}$.
\end{enumerate}
\end{enumerate}
\end{prop}

Since $\crpb{A}=\ccconstmod{A}\oplus \cks{i}{A}$ it is enough to give presentations of the modules $\ccconstmod{A}$ and $\cks{1}{A}$, which we do below (Corollary \ref{cor:zcpres} and 
Corollary \ref{cor:zipres}). Proposition \ref{prop:crpbpres} follows from these two corollaries and Corollary \ref{cor:key}.

In fact, in order to prove Proposition \ref{prop:crpbpres}, we will begin with a presentation of $Z_2(\Sigma_\bullet)$ as a $\pgl{2}{A}$-module and then  apply the following 
principle:

\begin{lem}\label{lem:presprin} Let $G$ be a group and let $M$ be a $\Z[G]$ module with presentation of the following form: $G$ is generated (as a $\Z[G]$-module) by the set $\mathcal{E}$, subject to the family of relations $\sum_{e\in\mathcal{E}}\alpha^{(r)}_e e=0$, where $r$ varies over an index set $R$ and $\alpha^{(r)}_e\in \Z[G]$. 

Let $N$ be a normal subgroup of $G$. Then the $\Z[G/N]$-module $M_N$ is generated by the set $\mathcal{E}$ subject to the relations $\sum_{e\in\mathcal{E}}\epsilon_{G/N}(\alpha^{(r)}_e)e=0$, $r\in R$, where $\epsilon_{G/N}:\Z[G]\to \Z[G/N]$ is the canonical surjection.
\end{lem}

We apply Lemma \ref{lem:presprin} at several places below in the case $G=\pgl{2}{A}$ and $N=\pgl{2}{A}$ or $\pspl{2}{A}$.

%Let $G_{\{ \infty,0\}}$ denote the subgroup 
%\[
%\an{ \left[
%\begin{array}{cc}
%u&0\\
%0&1\\
%\end{array}
%\right],\left[
%\begin{array}{cc}
%0&u\\
%1&0\\
%\end{array}
%\right]\ |\ u\in A^\times}\subset \pgl{2}{A}.
%\]

\begin{lem}\label{lem:zcpres} The $\pgl{2}{A}$-module $Z^c$ has the following presentation: It is generated by the element $\cconst{A}$ subject to the relations
\begin{enumerate}
\item  
 $\cconst{A}=\cconst{A}\cdot \omega$ where  $\omega:=
\left[
\begin{array}{cc}
0&1\\
1&0\\
\end{array}
\right]\in \pgl{2}{A}$.

\item 
\[\cconst{A}\cdot \left( I+M_{0,1,\infty}+M_{1,\infty,0}\right)=0.
\]
\end{enumerate}

\end{lem}

\begin{proof}   By Corollary \ref{cor:z2pres}, $Z^c$ is generated as an abelian group by the elements $C^{\{ x,y\}}(a)$, $\{ x,y,a\}\in Y_2$, subject to the relations 
\[
C^{\{ x,y\}}(z)+C^{\{ y,z\}}(x)+C^{\{ x,z\}}(y)=0
\]
for $\{ x,y,z\}\in Y_2$.

The $\pgl{2}{A}$-set of the generators $C^{\{ x,y\}}(a)\in \Z_2(\Sigma_\bullet)$ is in bijective correspondence with the indexing set 
$\mathcal{A}:=\{ (\{ x,y\},a)\ |\  \{ x,y\}\in Y_1, a\in Y^{\{ x,y\}}_0\}$. The group $\pgl{2}{A}$ acts transitively on $\mathcal{A}$  since $(\{ x,y\},a)=(\{ \infty,0\},1)\cdot M_{x,y,a}$.
Thus $Z^c$ is generated by the element $\cconst{A}=C^{\{ \infty,0\}}(1)$. 

Relation (1) follows from the fact that the stabilizer in $\pgl{2}{A}$ of $(\{ \infty,0\},1)\in \mathcal{A}$ is the subgroup of order $2$ generated by $\omega$. 

Relation (2) follows from the identity
\begin{eqnarray*}
C^{\{ x,y\}}(z)+C^{\{ y,z\}}(x)+C^{\{ x,z\}}(y)&=&\left( C^{\{ \infty,0\}}(1)+C^{\{ 0,1\}}(\infty)+C^{\{ \infty,1\}}(0)\right)\cdot M_{x,y,z}\\
&=&C^{\{ \infty,0\}}(1)\left( I+M_{0,1,\infty}+M_{1,\infty,0}\right)\cdot M_{x,y,z}.\\
\end{eqnarray*}
for all $\{ x,y,z\}\in Y_2(A)$.
\end{proof}

Note that 
\[
M_{0,1,\infty}= \left[
\begin{array}{cc}
0&-1\\
1&1\\
\end{array}
\right]\in\pspl{2}{A}\mbox{ and } M_{1,\infty,0}=
\left[
\begin{array}{cc}
1&1\\
-1&0\\
\end{array}
\right]=M_{0,1,\infty}^2\in\pspl{2}{A}.
\]

Recall that $\sgr{A}$ is the group ring $\Z[\sqcl{A}]$ where $\sqcl{A}:=A^\times/(A^\times)^2$ is the group of square classes in $A$. 
Let $\psqcl{A}$ denote the multiplicative group $A^\times/\pm(A^\times)^2$.

\begin{cor}\label{cor:zcpres}  The $\sgr{A}$-module $\ccconstmod{A}$ has the presentation: It is generated by $\cconst{A}$ subject to the relations
\begin{enumerate}
\item $\an{-1}\cconst{A}=\cconst{A}$, and
\item $3\cconst{A}=0$.
\end{enumerate}
Equivalently, the $\sgr{A}$-module $\ccconstmod{A}$ is a free rank one module over the group ring\\
 $\F{3}[\psqcl{A}]$ with generator $\cconst{A}$. 
\end{cor}

\begin{proof} Relation (1) follows from the observation that
 \[
G_{\{ \infty,0\}}\cap\pspl{2}{A}=
\an{ \left[
\begin{array}{cc}
u&0\\
0&u^{-1}\\
\end{array}
\right],\left[
\begin{array}{cc}
0&u\\
u^{-1}&0\\
\end{array}
\right]\ |\ u\in A^\times}
\]
and $\an{\mathrm{det}\left[
\begin{array}{cc}
u&0\\
0&u^{-1}\\
\end{array}
\right]}=\an{1}$ while $\an{\mathrm{det}\left[
\begin{array}{cc}
0&u\\
u^{-1}&0\\
\end{array}
\right]}=\an{-1}$ in $A^\times/(A^\times)^2$.
\end{proof}

Of course, taking $A^\times/(A^\times)^2$-coinvariants, we deduce:
\begin{cor}\label{cor:zcpresp} The $\Z$-module  $(Z^c)_{\pgl{2}{A}}\subset \cpb{A}$ is cyclic of order $3$ with generator $\cconst{A}$. 
\end{cor}

For $u\in A^\times$, let $\suss{i}{u}:=D_i^{\{ \infty,0\}}(1,u)\in Z^i\subset Z_2(\Sigma_\bullet)$.

\begin{lem}\label{lem:zipres} For $i=1,2$, the $\pgl{2}{A}$-module $Z^i$ has the following presentation: It is generated by the elements $\suss{i}{u},u\in A^\times$ subject to the relations
\begin{enumerate}
\item For all $u\in A^\times$, $\suss{i}{u}\cdot \left[
\begin{array}{cc}
0&1\\
1&0\\
\end{array}
\right]
=\suss{i}{u^{-1}}$ and,
\item for all $u,v\in A^\times$, $\suss{i}{uv}=\suss{i}{u}+\suss{i}{v}\cdot
\left[
\begin{array}{cc}
u&0\\
0&1\\
\end{array}
\right]$.
\end{enumerate}
\end{lem}

\begin{proof} By Corollary \ref{cor:z2pres}, $Z^i$ is generated as an abelian group by the elements $D_i^{\{ x,y\}}(a,b)$, $\{ x,y\}\in Y_1$, $a,b\in Y^{\{ x,y\}}_0$,
 subject to the relations $D_i^{\{ x,y\}}(a,b)+D_i^{\{ x,y\}}(b,c)=D_i^{\{ x,y\}}(a,c)$   for $\{x,y\}\in Y_1, a,b,c\in Y^{\{ x,y\}}_0$.

The generators are indexed by the right $\pgl{2}{A}$-set $\mathcal{B}:= \{ (\{x,y\}, a,b)\ |\ \{ x,y\}\in Y_1, a,b\in  Y^{\{ x,y\}}_0\}$. 

We note first that for 
$(\{ x,y\}, a,b)\in\mathcal{B}$ we have 
\[
(\{ x,y\}, a,b)=(\{\infty,0\},1,b\cdot M_{x,y,a}^{-1})\cdot M_{x,y,a}=((\{\infty,0\},1,u(x,y,a,b))\cdot M_{x,y,a}
\]
and hence
\[
D_i^{\{ x,y\}}(a,b)=\suss{i}{u(x,y,a,b)}\cdot M_{x,y,a}.
\]
This justifies the assertion about $\pgl{2}{A}$-generators. 

Relation (1) the follows from the observation that the stabilizer in $\pgl{2}{A}$ of $(\{ \infty,0\},1)$ is the cyclic group of order $2$ generated by ${\omega}:=
\left[
\begin{array}{cc}
0&1\\
1&0\\
\end{array}
\right]$ and that $u\cdot \omega=u^{-1}$ for all $u\in A^\times$. 

Relation (2) follows from the calculation: Given $\{x,y\}\in Y_1$ and $ a,b,c\in Y^{\{ x,y\}}_0$,\\
$D_i^{\{ x,y\}}(a,b)+D_i^{\{ x,y\}}(b,c)-D_i^{\{ x,y\}}(a,c)$
\begin{eqnarray*}
=\suss{i}{u(x,y,a,b)}\cdot M_{x,y,a}+\suss{i}{u(x,y,b,c)}\cdot M_{x,y,b}-\suss{i}{u(x,y,a,c)}\cdot M_{x,y,a}\\
=\left(
 \suss{i}{u(x,y,a,b)}+\suss{i}{u(x,y,b,c)}\cdot M_{x,y,b}M_{x,y,a}^{-1}-\suss{i}{u(x,y,a,c)}
\right)\cdot M_{x,y,a}\\
\end{eqnarray*}
and recall that $u(x,y,a,c)=u(x,y,a,b)\cdot u(x,y,b,c)$ by Corollary \ref{cor:cr} and note that 
\[
 M_{x,y,b}M_{x,y,a}^{-1}= \left[
\begin{array}{cc}
d(y,b)&0\\
0&d(b,x)\\
\end{array}
\right]\cdot T_{x,y}\cdot T_{x,y}^{-1}\cdot
\left[
\begin{array}{cc}
d(a,x)&0\\
0&d(y,a))\\
\end{array}
\right]=
\left[
\begin{array}{cc}
u(x,y,a,b)&0\\
0&1\\
\end{array}
\right]
\]
in $\pgl{2}{A}$. Furthermore, $M_{\infty,0,1}=I$ and hence $u(\infty,0,1,u)=u$ for any $u\in A^\times$.
\end{proof}

\begin{cor}\label{cor:zipres} For $i=1,2$, the $\sgr{A}$-module $\cks{i}{A}$ is generated by the elements $\suss{i}{u}$, $u\in A^\times$ subject only to the relations
\begin{enumerate}
\item $\an{-1}\suss{i}{u}=\suss{i}{u^{-1}}$ for all $u$ and,
\item $\suss{i}{uv}=\suss{i}{u}+\an{u}\suss{i}{v}$ for all $u,v\in A^\times$.
\end{enumerate}
 
The $\Z$-module $(Z^i)_{\pgl{2}{A}}\subset \cpb{A}$ is generated by the elements $\suss{i}{u}$, $u\in A^\times$ subject only to the relations
\begin{enumerate}
\item $\suss{i}{u}=\suss{i}{u^{-1}}$ for all $u$ and,
\item $\suss{i}{uv}=\suss{i}{u}+\suss{i}{v}$ for all $u,v\in A^\times$.
\end{enumerate}
\end{cor}

\begin{cor}\label{cor:zipres2}  For any ring $A$ and $i\in \{ 1,2\}$, there is a  natural isomorphism of groups $(Z^i)_{\pgl{2}{A}}\cong A^\times/(A^\times)^2$, given by
$\suss{i}{u}\mapsto \bar{u}$ for $u\in A^\times$.

Thus $\cpb{A}\cong A^\times/(A^\times)^2\oplus \Z/3$.
\end{cor}

\begin{proof} For any $u\in U$, $\suss{i}{u^2}=\suss{i}{u}+\suss{i}{u^{-1}}=\suss{i}{1}=0$ in $(Z^i)_{\pgl{2}{A}}$ by Corollary \ref{cor:zipres}. The second assertion follows from 
Corollary \ref{cor:zcpresp}.
\end{proof}

\begin{lem} \label{lem:psi} Fix $i\in \{ 1,2\}$ and, for $u\in A^\times$,  let $\psi(u):=\suss{i}{u}\in \cks{i}{A}$.  Then
\begin{enumerate}
\item $\pf{u}\psi(v)=\pf{v}\psi(u)$ for all $u,v\in U$.
\item $\psi(uv^2)=\psi(u)+\psi(v^2)$ for all $u.v$.
\item $2\psi(-1)=0$.
\item $\psi(u^2)=\pf{-1}\psi(u)=\pf{u}\psi(-1)$ for all $u$.
\item $2\psi(u^2)=0$ and $\psi(u^4)=0$ for all $u$.
\item  If $-1\in (A^\times)^2$, $\psi(u^2)=0$ for all $u$.
\item $\pf{u}\pf{v}\psi(-1)=0$ for all $u,v$. 
\item If $u^2=1$, then $\pf{u}\psi(-1)=0$.
\end{enumerate}
\end{lem}

\begin{proof}
\begin{enumerate}
\item The cocycle identity gives $\psi(uv)=\an{u}\psi(v)+\psi(u)=\psi(vu)=\an{v}\psi(u)+\psi(v)$.
\item This follows from the cocycle identity since $\an{v^2}=1$ in $\mathcal{R}$.
\item $0=\psi(1)=\psi(-1\cdot -1)=\an{-1}\psi(-1)+\psi(-1)=\psi((-1)^{-1})+\psi(-1)=2\psi(-1)$.
\item For any $u\in U$, we have $\psi(u)=\psi(u^{-1}\cdot u^2)=\psi(u^{-1})+\psi(u^2)=\an{-1}\psi(u)+\psi(u^2)$. Thus 
$\psi(u^2)=-\pf{-1}\psi(u)=-\pf{u}\psi(-1)=\pf{u}\psi(-1)$ by (1),(3).
\item By (3),(4), $2\psi(u^2)=2\pf{u}\psi(-1)=0$. Furthermore, $\psi(u^4)=2\psi(u^2)=0$ using (2).
\item If $-1\in (A^\times)^2$, then $\pf{-1}=0$ in $\mathcal{R}$ and hence $\psi(u^2)=\pf{-1}\psi(u)=0$ by (4).
\item $\pf{u}\pf{v}\psi(-1)=\pf{u}\psi(v^2)=\pf{v^2}\psi(u)=0$.
\item $\pf{u}\psi(-1)=\psi(u^2)=\psi(1)=0$.
\end{enumerate}
\end{proof}

We recall also that, for any $u\in A^\times$, the elements $\Suss{1}{u},\Suss{2}{u}$ and $\cconst{A}$ are related in $Z_2(\Sigma_\bullet)\subset Z_2(L_\bullet)$:

For $u\in A^\times$, we let
\[
M_u:=M_{\infty,0,u}= \left[
\begin{array}{cc}
u&0\\
0&1\\
\end{array}
\right]= \left[
\begin{array}{cc}
1&0\\
0&u^{-1}\\
\end{array}
\right]\mbox{ in } \pgl{2}{A}.
\]
\begin{lem}\label{lem:key} For any $u\in A^\times$ we have
\[
\cconst{A}\cdot \left( M_u-I\right)=\suss{1}{u}-\suss{2}{u}
\]
in $Z_2(\Sigma_\bullet)\subset Z_2(L_\bullet)$.
\end{lem}

\begin{proof} In $ Z_2(\Sigma_\bullet)$ we have 
\begin{eqnarray*}
\suss{1}{u}-\suss{2}{u}&=&D_1^{ \{ \infty,0\}}(1,u)-D_2^{ \{ \infty,0\}}(1,u)\\
&=& C^{ \{ \infty,0\}}(u)-C^{ \{ \infty,0\}}(1)\\
&=&C^{ \{ \infty,0\}}(1)\cdot M_u-C^{ \{ \infty,0\}}(1)\\
&=& \cconst{A}\cdot M_u-\cconst{A}.
\end{eqnarray*}
\end{proof}

Taking $\pspl{2}{A}$-coinvariants, we deduce
\begin{cor}\label{cor:key}
 In $ \crpb{A}$ we have the identity
\[
\pf{u}\cconst{A}:=\an{u}\cconst{A}-\cconst{A}=\suss{1}{u}-\suss{2}{u}.
\]
\end{cor}

\begin{rem} It follows that $\suss{1}{-1}=\suss{2}{-1}$ in $\crpb{A}$ for any ring $A$.
\end{rem}
%%%%%%%%%%%%%%%%%%%%%%%%%%%%%%%%%%%%%%%%%%%%%%%%%%%%%%%%%%%%%%%%%%%%%%%%%%%
%%%%%%%%%%%%%%%%%%%%%%%%%%%%%%%%%%%%%%%%%%%%%%%%%%%%%%%%%%%%%%%%%%%%%%%%
\section{The structure of $\crpb{A}$ and the homology groups $H_i(G,Z_2(\Sigma_\bullet))$}\label{sec:crpba}
%\subsection{Explicit structure of the $\sgr{A}$-module $\crpb{A}$}
%%%%%%%%%%%%%%%%%%%%%%%%%%%%%%%%%%%%%%%%%%%%%%%%%%%%%%%%%%%%%%%
Recall that the complex of $\pgl{2}{A}$-modules $L_\bullet^\tau$  gives rise, for any group $G$ mapping to $\pgl{2}{A}$,  to a spectral sequence $E(G, L_\bullet^\tau)$ converging to $H_n(G,L_\bullet^\tau)$. This latter group is isomorphic to $H_n(G,\Z)$ for many rings, and, in particular, for all fields and local rings. The spectral sequence has the form $E^1_{p,q}=H_p(G, L_q^\tau)\Longrightarrow H_{p+q}(G, L_\bullet^\tau)$. Now, $L_3^\tau=Z_2(L_\bullet)$ and thus the groups $H_p(G,Z_2(L_\bullet))$ occur on row $3$ of the $E^1$-term of the spectral sequence. In particular, $\rpb{A}=E^1_{0,3}(\pspl{2}{A}, L_\bullet^\tau)=H_0(\pspl{2}{A}, Z_2(L_\bullet))$. 

By Corollary \ref{cor:2dim}, if $A$ is a ring for which $\wn{A}=\emptyset$ ( including any ring for which $\F{2}$ occurs as a residue field ) then there is a short exact sequence of $\pgl{2}{A}$-modules 
\[
0\to Z_2(\Sigma_\bullet)\to Z_2(L_\bullet)=L_3^\tau\to H_2(|Y(A)|)\to 0.
\]
Thus, for any group $G$ mapping to $\pgl{2}{A}$ there is a long exact sequence of homology groups
\[
\cdots \to H_p(G,Z_2(\Sigma_\bullet))\to H_p(G,Z_2(L_\bullet))\to H_p(G,|Y(A)|)\to H_{n-1}(G,Z_2(\Sigma_\bullet))\to\cdots
\]
Thus it is useful to calculate the homology groups $H_p(G,Z_2(\Sigma_\bullet))$. Of course, in particular, $\crpb{A}=H_0(\pspl{2}{A},Z_2(\Sigma_\bullet))$.

Since $Z_2(\Sigma_\bullet)$ decomposes over $\pgl{2}{A}$ as the direct sum $Z^c\oplus Z^i$ for $i\in \{ 1,2\}$, we have $H_p(G,Z_2(\Sigma_\bullet))=H_p(G, Z^c)\oplus H_p(G,Z^i)$ 
for any group $G$ mapping to $\pgl{2}{A}$. We now consider the structure of the modules $Z^c$ and $Z^i$ separately. 
\subsection{Conjugation}\label{sub:conj}
 Let $H$ be a subgroup of a group $G$. For $g\in G$, we let $^gh:=ghg^{-1}$ (resp. $h^g:=g^{-1}hg$) for any $h\in H$ and $^gH:=gHg^{-1}$ (resp. $H^g:=g^{-1}Hg$).
We denote by $\lconj{g}:H\to {}^gH$ the isomorphism -- left conjugation -- $h\mapsto {}^gh$ and by $\rconj{g}:H\to H^g$ the right conjugation $h\mapsto h^g$. (Of course,
$\rconj{g}=\lconj{g^{-1}}$).

We recall that if $H$ is a subgroup of $G$ and if $M$ is a (right) $\Z[G]$-module, then the pair of maps $\lconj{g}:H\to {}^gH$ and $\alpha_g:M\to M$, $m\mapsto mg^{-1}$ induce a map on homology $\ho{\bullet}{H}{M}\to \ho{\bullet}{{}^gH}{M}$, which we will denote $L_g$.  Let $P_\bullet$ be a projective resolution 
of $\Z$ over $\Z[G]$. Then the map $L_g$ is induced by the map $f:M\otimes_{\Z[H]}P_\bullet\to M\otimes_{\Z[{}^gH]}P_\bullet$, $f(m\otimes p):=mg^{-1}\otimes gp$. 
(In the case where $M=\Z$ with trivial action, we will also use the notation $g\cdot$ instead of $L_g$.)

If $H$ is normal in $G$, this gives rise to an action of $G$ by conjugation on $\ho{\bullet}{H}{M}$, $g\mapsto L_g$. This action clearly restricts to the trivial action of $H$ and hence induces a (left) $\Z[G/H]$-module 
structure on $\ho{\bullet}{H}{M}$.

The particular case of most relevance below is the following: Let $A$ be a ring and let $M$ be a right $\pgl{2}{A}$-module (and hence a right $\gl{2}{A}$-module). Then the action of $\gl{2}{A}$ by conjugation makes 
$\ho{\bullet}{\spl{2}{A}}{M}$ a  module over $\Z[\gl{2}{A}/\spl{2}{A}]\cong\Z[A^\times]$. However,  the central group of scalar matrices $Z$ acts trivially and hence this action descends to an action of the group ring $\sgr{A}:=\Z[A^\times/(A^\times)^2]\cong \Z[\pgl{2}{A}/\pspl{2}{A}]$.  Thus the action of $u\in A^\times/(A^\times)^2$ is induced by conjugation by  any element of $\gl{2}{A}$ with determinant $u\mod{(A^\times)^2}$.

%We let $g\cdot$ denote the functorial map $\ho{\bullet}{H}{\Z}\to\ho{\bullet}{{}^gH}{\Z}$ induced by the conjugation by $g$  isomorphism 
%$H\to {}^gH$.
\subsection{Permutation modules}\label{sub:oerm} We will describe these modules in terms of permutation modules. We review first a few basic facts that we will use:

Let $X_1$ and $X_2$ be transitive right $\Z[G]$-modules. For $i=1,2$, let $x_i\in X_i$ have stabiizer subgroup $H_i$ in $G$.

Now let $\phi:\Z[X_1]\to \Z[X_2]$ be a homomorphism of $\Z[G]$-modules.  Since $G$ act transitively on $X_1$, $\phi$ is determined by its value at $x_1$. Suppose that 
$\phi(x_1)=\sum_{g\in H_2\backslash G}n_g\cdot (x_2g)$ with $n_g\in \Z$.  Let 
\[
\phi_\bullet:\ho{\bullet}{H_1}{\Z}\cong \ho{\bullet}{G}{\Z[X_1]}\to\ho{\bullet}{G}{\Z[X_2]}\cong\ho{\bullet}{H_2}{\Z}
\]
be the map induced by $\phi$, and the Shapiro isomorphisms, on the homology of the groups $H_i$. Let $\supp{\phi}:=\{ g\in H_2\backslash G\ |\ n_g\not=0\}$.
\begin{lem}\label{lem:phi}\cite[Lemma 3]{hut:mat}
\begin{enumerate}
\item Let $E:=H_2\backslash G/H_1$. Then $n_g$ depends only on the class of $g$ in $E$.
\item For all $g\in \supp{\phi}$, the index $[H_1:H_1\cap H_2^g]$ is finite.
\item For $g\in G$, Let $K_g:={}^gH_1\cap H_2$. Then 
\[
\phi_\bullet(z)=\sum_{g\in E}n_g\Cor{K_g}{H_2}\Res{^gH_1}{K_g}(g\cdot z), \quad z\in \ho{\bullet}{H_1}{\Z}.
\]
\end{enumerate}
\end{lem}

If $G$ is a group and if $H$ is a subgroup the right $G$-set $H\backslash G$ gives rise to the right permutation $\Z[G]$-module $\Z[H\backslash G]$. There is a natural isomorphism of right 
$\Z[G]$-modules $\mathrm{Ind}^G_H \Z=\Z\otimes_{\Z[H]}\Z[G]\cong \Z[H\backslash G]$ given by $1\otimes g\leftrightarrow Hg$.  By Shapiro's Lemma, there are thus isomorphisms 
$H_p(G,\Z[H\backslash G])\cong H_p(G,\mathrm{Ind}^G_H\Z)\cong H_p(H,\Z)$.

If $H$ and $K$ are subgroups of $G$  then the right $G$-set $H\backslash G$ is also a right $K$-set. For any $g\in G$, the orbit under $K$ of the coset $Hg$ is $H\backslash Hg{K}$ and the stabilizer of $Hg\in H\backslash G$ is $H^g\cap K$. We easily deduce:

\begin{lem}\label{lem:permmod} Let $H$ and $K$ be  subgroups of the group $G$ . Let $\epsilon:H\backslash G/K\to G$ be a section 
of the canonical map $G\to H\backslash G/K$ with image $E\subset G$. For each $g\in G$ let $\rho(g)\in K$ satisfy $Hg=H\epsilon(\bar{g}){\rho(g)}$. Then there is a natural isomorphism 
of right $\Z[K]$-modules 
\begin{eqnarray*}
\mathrm{Ind}_H^G\Z& \cong &\bigoplus_{e\in E}\mathrm{Ind}^K_{H^e\cap K}\Z\\
1\otimes g&\rightarrow& (1\otimes \rho(g))_{\epsilon(\bar{g})}\\
1\otimes e{k}&\leftarrow& (1\otimes k)_e\\
\end{eqnarray*}
\end{lem}

\begin{rem}\label{rem:erho} Observe that the direct sum decomposition occurring here is non-canonical and depends on the choice of section $\epsilon$:
If $HeK=HfK$ then it does not generally follow that $H^e\cap K=H^{f}\cap K$. (Consider the trivial case $K=G$.)

However, note, on the other hand, that the isomorphisms given are  independent of the choice of $\rho(g)$ (once $\epsilon$ is given): If $\rho$, $\rho'\in K$ satisfy
$He\rho=He\tilde{\rho}$ then $H^e\rho=H^e\tilde{\rho}$ and hence $(H^e\cap K)\rho=(H^e\cap K)\tilde{\rho}$. Thus $1\otimes \rho=1\otimes \tilde{\rho}$ in $\Z\otimes_{\Z[H^e\cap K]}\Z[K]$.
\end{rem}
\begin{cor}\label{cor:shapiro} Let $G$ be a group with subgroups $H$ and $K$. Let $\epsilon$, $E$, $\rho$ be given as in Lemma \ref{lem:permmod}. Let $P_\bullet$ be a $\Z[K]$-projective 
resolution of $\Z$. Then there is an isomorphism of complexes 
\[
(\Ind{H}{G}{\Z})\otimes_{\Z[K]}P_\bullet\cong \bigoplus_{e\in E}\Z\otimes_{\Z[H^e\cap K]}P_\bullet
\]
given by $(1\otimes g)\otimes p\mapsto (1\otimes\rho(g)p)_{\epsilon(\bar{g})}\in \Z\otimes_{\Z[H^{\epsilon{\bar{g}}}\cap K]}P_\bullet$ and thus  inducing an isomorphism
\[
\ho{\bullet}{K}{\Ind{H}{G}{\Z}}\cong\bigoplus_{e\in E}\ho{\bullet}{H^e\cap K}{\Z}.
\]
\end{cor}

%%%%%%%%%%%%%%%%%%%%%%%%%%%%%%%%%%%%%%%%%%%%%%%%%%%%%%%%%%%%%%%%
\subsection{The homology of $Z^c$}
%%%%%%%%%%%%%%%%%%%%%%%%%%%%%%%%%%%%%%%%%%%%%%%%%%%%%%%%%%%%%%%%%%%%%

Recall that for any ring $A$ there is  a  nonabelian subgroup, $S$ say, of order $6$ in $\gl{2}{A}$ generated by 
$\omega:=
\left[
\begin{array}{cc}
0&1\\
1&0\\
\end{array}
\right]$ (of order $2$) and 
$M:=M_{0,1,\infty}=
\left[
\begin{array}{cc}
0&-1\\
1&1\\
\end{array}
\right]$ of order $3$.  Let $C_2$ denote the subgroup of order $2$ generated by $\omega$ and let $C_3$ denote the subgroup of order $3$  generated by $M$.  Thus $S=C_3\rtimes C_2$. We also denote the (isomorphic) image of $S$ in $\pgl{2}{A}$ by $S$.

From Lemma \ref{lem:zcpres}, we can immediately deduce:

\begin{lem}\label{lem:zcperm}  Let $A$ be a ring. Then there is a natural short exact sequence of right $\pgl{2}{A}$-modules
\[
\xymatrix{
0\ar[r]&\Ind{S}{\pgl{2}{A}}\Z\ar^-{\beta}[r]&\Ind{C_2}{\pgl{2}{A}}\Z\ar^-{\pi}[r]&Z^c\ar[r]&0\\
%0\ar[r]&\Z\otimes_{\Z[S]}\Z[\pgl{2}{A}]\ar^-{\beta}[r]&\Z\otimes_{\Z[C_2]}\Z[\pgl{2}{A}]\ar^-{\pi}[r]&Z^c\ar[r]&0\\
}
\]
where $\beta(1\otimes g):=1\otimes  (I+M+M^2)\cdot g$ and $\pi(1\otimes g)=\cconst{A}\cdot g$.
\end{lem}

Since $\pgl{2}{A}=\gl{2}{A}/Z=Z\backslash \gl{2}{A}$ where $Z$ is the group of scalar matrices,  we thus deduce :
\begin{cor}\label{cor:zcperm}
 Let $A$ be a ring. Then there is a natural short exact sequence of right $\gl{2}{A}$-modules
\[
\xymatrix{
0\ar[r]&\Ind{S\times Z}{\gl{2}{A}}\Z\ar^-{\beta}[r]&\Ind{C_2\times Z}{\gl{2}{A}}\Z\ar^-{\pi}[r]&Z^c\ar[r]&0\\
%0\ar[r]&\Z\otimes_{\Z[S\times Z]}\Z[\gl{2}{A}]\ar^-{\beta}[r]&\Z\otimes_{\Z[C_2\times Z]}\Z[\gl{2}{A}]\ar^-{\pi}[r]&Z^c\ar[r]&0\\
}
\]
where $\beta(1\otimes g):=1\otimes  (I+M+M^2)\cdot g$ and $\pi(1\otimes g)=\cconst{A}\cdot g$.
\end{cor}

Now we wish to restrict to $\spl{2}{A}$ and use Corollary \ref{cor:shapiro}. There are two cases to consider: (a) $-1\not\in (A^\times)^2$ and (b) $-1\in (A^\times)^2$:

For simplicity, we define the following subgroups of $\spl{2}{A}$:
\[
\tilde{C}_2:=
\left\{
\begin{array}{ll}
1,& -1\not\in (A^\times)^2\\
\an{i\omega}, & -1=i^2 \mbox{ in }A\\
\end{array}
\right.
\]
and
\[
\tilde{S}:=C_3\rtimes \tilde{C}_2=
\left\{
\begin{array}{ll}
C_3,& -1\not\in (A^\times)^2\\
C_3\rtimes\an{i\omega}, & -1=i^2 \mbox{ in }A.\\
\end{array}
\right.
\]
Thus, for any ring $A$, we have 
\[
(S\times Z)\cap \spl{2}{A}=\tilde{S}\times\mu_2 \mbox{ and }(C_2\times Z)\cap\spl{2}{A}=\tilde{C}_2\times \mu_2
\]
  (where $\mu_2=\mu_2(A)\subset Z\cong A^\times$ 
is the subgroup $\{ z\in A\ |\ z^2=1\}$).

Now let $\oodet:\gl{2}{A}\to A^\times/\pm(A^\times)^2$ be the map induced by the determinant. (Of course, $\oodet=\odet$ if and only if $-1\in (A^\times)^2$.) Since $\spl{2}{A}$ is normal in $\gl{2}{A}$ and since $C_3\subset \spl{2}{A}$ we have natural bijections of double coset spaces
\[
(S\times Z)\backslash \gl{2}{A}/\spl{2}{A}\leftrightarrow (C_2\times Z)\backslash \gl{2}{A}/\spl{2}{A}
\]
and both are in turn in bijective correspondence with the groups
\[
\xymatrix{
\gl{2}{A}/(C_2\times Z)\cdot\spl{2}{A}\cong \pgl{2}{A}/C_2\cdot\pspl{2}{A}\ar_-{\cong}^-{\oodet}[r]&A^\times/\pm(A^\times)^2.
}
\]
Let $\epsilon:A^\times/\pm(A^\times)^2\to \gl{2}{A}$ be a section of $\oodet$ with image $U$ and such that $\epsilon(1)=1$. Thus, by definition, there is a bijection
\begin{eqnarray*}
U&\leftrightarrow&A^\times/\pm(A^\times)^2\\
u&\rightarrow&\tilde{u}:=\oodet{u}\\
\epsilon(v)&\leftarrow&v.\\
\end{eqnarray*}
For convenience, if $g\in \gl{2}{A}$ we will simply let $\epsilon(g)$ denote $\epsilon(\oodet(g))$ below. 

For any $g\in \gl{2}{A}$ choose 
$\rho(g)\in\spl{2}{A}$ such that $(C_2\times Z)\epsilon(\oodet(g))\rho(g)=(C_2\times Z)g$ (and hence also $(S\times Z)\epsilon(\oodet(g))\rho(g)=(S\times Z)g$). 

By 
Corollary \ref{cor:shapiro} we have isomorphisms of right $\Z[\spl{2}{A}]$-modules
\begin{eqnarray*}
\Ind{S\times Z}{\gl{2}{A}}\Z&\cong & \bigoplus_{u\in U}\Ind{\tilde{S}^u\times\mu_2}{\spl{2}{A}}\Z\\
\Ind{C_2\times Z}{\gl{2}{A}}\Z&\cong & \bigoplus_{u\in U}\Ind{\tilde{C}_2^u\times\mu_2}{\spl{2}{A}}\Z\\
\end{eqnarray*}
via the maps
\begin{eqnarray*}
1\otimes g&\rightarrow & (1\otimes\rho(g))_{\epsilon(g)}\\
1\otimes ug&\leftarrow& (1\otimes g)_u.\\
\end{eqnarray*}

In view of Corollary \ref{cor:zcperm}, we immediately deduce:
\begin{cor}\label{cor:U} There is a short exact sequence of right $\Z[\spl{2}{A}]$-modules 
\[
\xymatrix{
0\ar[r]&\bigoplus_{u\in U}\Ind{\tilde{S}^u\times\mu_2}{\spl{2}{A}}\Z\ar^-{\beta}[r]&\bigoplus_{u\in U}\Ind{\tilde{C}_2^u\times\mu_2}{\spl{2}{A}}\Z\ar^-{\pi}[r]&Z^c\ar[r]&0\\
}
\]
where $\beta=\oplus_{u\in U}\beta^u$ with 
\[
\beta^u((1\otimes g)_u):= (1\otimes (1+M^u+(M^u)^2)g)_u \mbox{ and  } \pi((1\otimes g)_u)=\cconst{A}\cdot ug.
\]
\end{cor}

\begin{proof}
For the calculation of $\beta$, we note that $(1\otimes g)_u\in \Ind{\tilde{S}^u\times\mu_2}{\spl{2}{A}}\Z$ corresponds to $1\otimes ug\in \Ind{S\times Z}{\gl{2}{A}}\Z$ 
which maps under $\beta$ to 
\[
1\otimes (1+M+M^2)ug=1\otimes u\cdot (1+M^u+(M^u)^2)g=\sum_{i=0}^2 1\otimes u\cdot (M^u)^ig\in \Ind{C_2\times Z}{\gl{2}{A}}\Z
\]
which in turn corresponds to \[
\sum_{i=0}^2 (1\otimes \rho((M^u)^ig))_u=\sum_{i=0}^2 (1\otimes (M^u)^ig)_u\in \Ind{\tilde{C}_2^u\times\mu_2}{\spl{2}{A}}\Z.
\]
This last equality follows  from 
\begin{eqnarray*}
(C_2\times Z)u\cdot (M^u)^ig&=&(C_2\times Z)u\rho((M^u)^ig)\\
\imp (C_2\times Z)^u\cdot (M^u)^ig&=&(C_2\times Z)^u\rho((M^u)^ig)\\
\imp \left((C_2\times Z)^u\cap \spl{2}{A}\right)\cdot (M^u)^ig&=&\left((C_2\times Z)^u\cap\spl{2}{A}\right)\rho((M^u)^ig)\\
\imp (\tilde{C}_2^u\times\mu_2)\cdot (M^u)^ig&=&(\tilde{C}_2^u\times\mu_2)\rho((M^u)^ig)\\
\end{eqnarray*}
(since $(M^u)^ig$ and $\rho((M^u)^ig)$ both  lie in $\spl{2}{A}$).
\end{proof}

Taking $\spl{2}{A}$-homology we thus obtain a long exact homology sequence:

\begin{prop} \label{prop:leszc}There is a long exact sequence of the form 
\[
\xymatrix{
\cdots\ar[r]&\bigoplus_{u\in U}\ho{k}{\tilde{S}^u\times\mu_2}{\Z}\ar^-{\beta_k}[r]&\bigoplus_{u\in U}\ho{k}{\tilde{C}_2^u\times\mu_2}{\Z}\ar[r]&\ho{k}{\spl{2}{A}}{Z^c}\\
\ar^-{\delta_k}[r]&\bigoplus_{u\in U}\ho{k-1}{\tilde{S}^u\times\mu_2}{\Z}\ar^-{\beta_{k-1}}[r]&\bigoplus_{u\in U}\ho{k-1}{\tilde{C}_2^u\times\mu_2}{\Z}\ar[r]&\cdots
}
\]
where $\beta_k=\oplus_{u\in U}\beta_k^u$ and 
\[
\beta_k^u=
\left\{
\begin{array}{ll}
\Res{\tilde{C}_2^u\times\mu_2}{\tilde{S}^u\times\mu_2},& \mbox{ if }k>0,\\
\mbox{ mult. by } 3:\Z\to \Z,& \mbox{ if }k=0.\\
\end{array}
\right.
\]
\end{prop}

\begin{proof} We need only explain the identification of $\beta_k^u:\ho{k}{\tilde{S}^u\times\mu_2}{\Z}\to \ho{k}{\tilde{C}_2^u\times\mu_2}{\Z}$:

 By Corollary \ref{cor:U} and Lemma \ref{lem:phi}
\begin{eqnarray*}
\beta_k^u&=& \sum_{i=0}^2\Res{\tilde{C}_2^u\times\mu_2}{\tilde{S}^u\times\mu_2}(M^u)^i\cdot\\
&=& 3\Res{\tilde{C}_2^u\times\mu_2}{\tilde{S}^u\times\mu_2}\mbox{ since } M^u\in \tilde{S}\times\mu_2\\
&=&
\left\{
\begin{array}{ll}
\Res{\tilde{C}_2^u\times\mu_2}{\tilde{S}^u\times\mu_2},& \mbox{if } k>0 (\mbox{ since } 2\cdot  \ho{k}{\tilde{C}_2^u\times\mu_2}{\Z}=0)\\
\mbox{ mult. by } 3:\Z\to \Z,& \mbox{ if }k=0.
\end{array}
\right.\\
\end{eqnarray*}
\end{proof}

\begin{lem}\label{lem:res} Let $A$ be a ring such that $\mu_2=\mu_2(A)$ is finite.  Let $\tilde{S}$, $\tilde{C}_2$ be as above. For all $k>0$ the map\\
 $\Res{\tilde{C}_2\times\mu_2}{\tilde{S}\times\mu_2}:\ho{k}{\tilde{S}\times\mu_2}{\Z}\to\ho{k}{\tilde{C}_2\times\mu_2}{\Z}$ is surjective with kernel
\[
\ker{\Res{\tilde{C}_2\times\mu_2}{\tilde{S}\times\mu_2}}\cong 
\left\{
\begin{array}{ll}
\ho{k}{C_3}{\Z}^{C_2},& -1=i^2 \mbox{ in } A^\times,\\
\ho{k}{C_3}{\Z},& -1\not\in (A^\times)^2\\
\end{array}
\right.
\]
where $C_2$ acts non-trivially on $C_3$. 
\end{lem}

\begin{proof} 
The elementary abelian $2$-group $\tilde{C}_2\times \mu_2$ is the $2$-Sylow subgroup of  $\tilde{S}\times \mu_2$. Since the corresponding inclusion of groups is split, it follows that 
$\Cor{\tilde{C}_2\times\mu_2}{\tilde{S}\times \mu_2}:\ho{k}{\tilde{C}_2\times\mu_2}{\Z}\to\ho{k}{\tilde{S}\times \mu_2}{\Z}$ is injective and, when $k>0$,  its image is the $2$-primary part of 
$\ho{k}{\tilde{S}\times \mu_2}{\Z}$. Now, for $k>0$,  the image of $\Res{\tilde{C}_2\times\mu_2}{\tilde{S}\times \mu_2}:\ho{k}{\tilde{S}\times \mu_2}{\Z}\to\ho{k}{\tilde{C}_2\times\mu_2}{\Z}$
is isomorphic to the $2$-Sylow subgroup of $\ho{k}{\tilde{S}\times \mu_2}{\Z}$ (see, eg., \cite[Theorem III.10.3]{brown:coh}).
It follows immediately that $\Res{\tilde{C}_2\times\mu_2}{\tilde{S}\times\mu_2}$ is surjective and its kernel is the $3$-primary part of $\ho{k}{\tilde{S}\times\mu_2}{\Z}$.

By \cite[Theorem III.10.3]{brown:coh} again, the $3$-primary part of $\ho{k}{\tilde{S}\times\mu_2}{\Z}$ is the subgroup of invariant homology classes in $\ho{k}{C_3}{\Z}$. If $-1\not\in (A^\times)^2$,
$C_3$ is central and hence this is all of $\ho{k}{C_3}{\Z}$. If $-1\in (A^\times)^2$, then $i\omega$ acts as by conjugation on $C_3$ as inversion, and the $3$-primary part is precisely the 
group of homology classes invariant under this action. 
\end{proof}
\begin{cor}\label{cor:kerres} For all $k>0$,
\[
\ker{\Res{\tilde{C}_2\times\mu_2}{\tilde{S}\times\mu_2}}\cong 
\left\{
\begin{array}{ll}
\Z/3,&  $k$\mbox{ odd and either } -1\not\in (A^\times)^2 \mbox{ or } k\equiv 3\pmod{4}\\
0,& \mbox{ otherwise.}\\
\end{array}
\right.
\]
\end{cor}
\begin{proof} It is straightforward  that the nontrivial automorphism of $C_3$ induces the map $-1$ on $\ho{k}{C_3}{\Z}\cong \Z/3$ when $k\equiv 1\pmod{4}$ and the identity map when 
$k\equiv 3\pmod{4}$. 
\end{proof}

\begin{cor}\label{cor:hksl2zc}  Let $A$ be a ring such that $\mu_2=\mu_2(A)$ is finite. Then
\[
\ho{k}{\spl{2}{A}}{Z^c}\cong
\left\{
\begin{array}{ll}
\bigoplus_{u\in U}\Z/3,&  $k$\mbox{ even and either } -1\not\in (A^\times)^2 \mbox{ or } k\equiv 0\pmod{4}\\
0,&\mbox{ otherwise}
\end{array}
\right.
\]
\end{cor}

\begin{proof} The long exact homology sequence of Proposition \ref{prop:leszc}, together with Lemma \ref{lem:res} gives isomorphisms
\[
\ho{k}{\spl{2}{A}}{\Z}\cong \oplus_{u\in U}\Ker{\Res{\tilde{C}_2^u\times\mu_2}{\tilde{S}^u\times\mu_2}:\ho{k-1}{\tilde{S}^u\times\mu_2}{\Z}\to\ho{k-1}{\tilde{C}_2^u\times\mu_2}{\Z}}
\]
for $k\geq 1$ and a short exact sequence
\[
\xymatrix{
0\ar[r]&\bigoplus_{u\in U}\Z\ar^-{3\cdot}[r]&\bigoplus_{u\in U}\Z\ar[r]&\ho{0}{\spl{2}{A}}{Z^c}\ar[r]&0.\\
}
\]
The result then follows from Corollary \ref{cor:kerres}.
\end{proof}

Recall that if $M$ is any $\pgl{2}{A}$-module then $\ho{k}{\spl{2}{A}}{M}$ is naturally a $\sgr{A}$-module. It follows that the long exact sequence of Proposition \ref{prop:leszc} is a sequence of $\sgr{A}$-modules. 
We will now make this $\sgr{A}$-action explicit. 

We begin by observing that for any $v\in A^\times/(A^\times)^2$ we have $v\cdot\odet(\epsilon({v}))\in \pm(A^\times)^2/(A^\times)^2$.
%  (where $\bar{v}$ denotes the image of $v$ in $A^\times/\pm(A^\times)^2$). 

Thus, for $v\in A^\times/(A^\times)^2$ we can define $\chi(v)=\chi_{\epsilon}(v)\in \Z/2$ as follows:\\
If $-1\not\in (A^\times)^2$, we require that $(-1)^{\chi(v)}=v\cdot\odet(\epsilon({v}))$. Otherwise $\chi(v):=0$.  

More generally, for $g\in \gl{2}{A}$, we can define $\chi(g):=\chi(\odet(g))$. Observe that if $-1\not\in (A^\times)^2$, since $1=\epsilon(1)=\epsilon(-1)$, we have $\chi(g)=-1$ 
whenever $\odet(g)=-1$. In particular, $\chi(\omega)=\chi(-1)=-1$. 

On the other hand, for $u\in U$, we have $u=\epsilon(\oodet(u))=\epsilon({\odet(u)})$ and hence $\chi(u)=\chi(\odet(u))=1$.

In the proof of following Lemma we will use the fact that $\odet(\omega^{\chi(g)})=(-1)^{\chi(g)}=\odet(g\epsilon(g))$ in $A^\times/(A^\times)^2$ and hence
$\odet(\epsilon(g)^{-1}\omega^{\chi(g)})= \odet(g)$.

\begin{lem}\label{lem:raact} Let $M=\Ind{H}{\gl{2}{A}}\Z$ where $H=S\times Z$ or $C_2\times Z$. Let $H_0=H\cap \spl{2}{A}$. So 
\[
\ho{k}{\spl{2}{A}}{M}\cong\bigoplus_{u\in U}\ho{k}{H_0^u}{\Z}
\]
and, for any $g\in\gl{2}{A}$, the resulting action of $\odet(g)\in A^\times/(A^\times)^2$ on the right-hand side is as follows: 

For any $u\in U$, $\odet(g)$ maps $\ho{k}{H_0^u}{\Z}$ to $\ho{k}{H_0^{\epsilon(ug)}}{\Z}$ via the map induced by right conjugation by $u^{-1}\omega^{\chi(ug)}\epsilon(ug)$. 
\end{lem}
\begin{proof}  Let  $v\in A^\times/(A^\times)^2$.  The action of $g$ on $\ho{\bullet}{\spl{2}{A}}{M}$ is equal to the homomorphism $L_g$, which is induced by the chain map sending  $m\otimes p$ to $mg^{-1}\otimes gp$  where $m\otimes p\in M\otimes_{\Z[\spl{2}{A}]}P_\bullet$ and $P_\bullet$ is a projective resolution of $\Z$ over $\Z[\gl{2}{A}]$.

Now let $u\in U$ and let the homology class $z\in\ho{\bullet}{H_0^u}{\Z}$ be represented by the cycle $1\otimes p\in \Z\otimes_{\Z[H_0^u]}P_\bullet$. This corresponds to the term
$(1\otimes 1)_u\otimes p\in \Ind{H_0^u}{\spl{2}{A}}\Z\otimes_{\Z[\spl{2}{A}]}P_\bullet$, which in turn corresponds to $(1\otimes u)\otimes p\in \Ind{H}{\gl{2}{A}}\Z\otimes_{\Z[\spl{2}{A}]}P_\bullet=M\otimes_{\Z[\spl{2}{A}]}P_\bullet$ which represents the corresponding homology class in $\ho{\bullet}{\spl{2}{A}}{M}$. The image of this class under $L_g$ is thus represented by 
the cycle $(1\otimes ug^{-1})\otimes gp\in M\otimes_{\Z[\spl{2}{A}]}P_\bullet$. 

Observe that $\epsilon(ug^{-1})=\epsilon(ug)$ since $\oodet(g^{-1})=\oodet(g)$. Thus the cycle $(1\otimes ug^{-1})\otimes gp$ in turn corresponds to 
\[
(1\otimes \rho(ug^{-1}))_{\epsilon(ug)}\otimes gp \mbox{ in }\Ind{H_0^{\epsilon(ug)}}{\spl{2}{A}}\Z\otimes_{\Z[\spl{2}{A}]}P_\bullet
\]
 which corresponds to 
$1\otimes \rho(ug^{-1})gp$ in $\Z\otimes_{\Z[H_0^{\epsilon(ug)}]}P_\bullet$. Thus the overall effect is to send $1\otimes p\in \Z\otimes_{\Z[H_0^u]}P_\bullet$ to 
$1\otimes \rho(ug^{-1})gp\in \Z\otimes_{\Z[H_0^{\epsilon(ug)}]}P_\bullet$. On the corresponding homology groups, this is just left conjugation by $\rho(ug^{-1})g$.

Finally, consider $\rho:=\epsilon(ug)^{-1}\omega^{\chi(ug)}ug^{-1}\in \gl{2}{A}$. Observe that $\odet(\rho)=1$ and that 
\[
H\epsilon(ug^{-1})\rho=H\epsilon(ug)\rho=H\omega^{\chi(ug)}ug^{-1}=Hug^{-1}.
\]
It follows that we can let $\rho(ug^{-1}):=\rho \cdot zI$ where $z\in A^\times$ satisfies $z^{-2}=\det(\rho)$. Hence $\rho(ug^{-1})g=\epsilon(ug)^{-1}\omega^{\chi(ug)}u\cdot zI$. Thus left conjugation by $\rho(ug^{-1})g$ is equal to right conjugation by $u^{-1}\omega^{\chi(ug)}\epsilon(ug)$ (since $zI$ is central) as required.
\end{proof}

For a ring $A$, let $\qsgr{A}$ denote the quotient group ring 
$\Z[A^\times/\pm(A^\times)^2)]$ of $\sgr{A}$. Note that $\sgr{A}/\sgr{A}\cdot\pf{-1}\cong \qsgr{A}$ and that $\qsgr{A}=\sgr{A}$ if $-1\in(A^\times)^2$.

If $-1\not\in (A^\times)^2$, define $\msgr{A}:=\sgr{A}/\sgr{A}\cdot(\an{-1}+1))$. Observe that $\msgr{A}\otimes\Z/2=\qsgr{A}\otimes \Z/2$ and, more generally, if $M$ is any $\Z/2$-module then $\msgr{A}\otimes  M=\qsgr{A}\otimes M$ (when $-1\not\in (A^\times)^2$).

It is clear from the definitions that, for any ring $A$ with $-1\not\in (A^\times)^2$,  the set $\{ \an{\odet(u)}\}_{u\in U}$ is a $\Z$-basis of $\qsgr{A}$ or of $\msgr{A}$.

For $k>0$ odd, let $u_k$ denote the image of a generator of $\ho{k}{C_3}{\Z}$ in $\ho{k}{\tilde{S}\times\mu_2}{\Z}$. 

\begin{prop}\label{prop:hokzc} Let $A$ be a ring for which $\mu_2(A)$ is finite. We have isomorphisms of $\sgr{A}$-modules:
\begin{enumerate}
\item
$\ho{k}{\spl{2}{A}}{\Ind{{C}_2\times Z}{\gl{2}{A}}{\Z}}\cong \qsgr{A}\otimes\ho{k}{\tilde{C}_2\times\mu_2}{\Z}$ for all $k\geq 0$.
\item
\[
\ho{k}{\spl{2}{A}}{\Ind{S\times Z}{\gl{2}{A}}{\Z}}\cong 
\left\{
\begin{array}{ll}
\msgr{A}\otimes \ho{k}{\tilde{S}\times\mu_2}{\Z}, &  -1\not\in (A^\times)^2\mbox{ and } k\equiv 1\pmod{4},\\
\qsgr{A}\otimes \ho{k}{\tilde{S}\times\mu_2}{\Z}, &\mbox{ otherwise.}\\
\end{array}
\right.
\]
\item 
\[
\ho{k}{\spl{2}{A}}{Z^c}\cong 
\left\{
\begin{array}{ll}
\qsgr{A}\otimes \Z/3, &  k\equiv 0\pmod{4},\\
\msgr{A}\otimes \Z/3, & k\equiv 2\pmod{4} \mbox{ and } -1\not\in (A^\times)^2,\\
0,&\mbox{ otherwise.}\\
\end{array}
\right.
\]
with generator $\cconst{A}$ when $k=0$ and $\delta_k^{-1}(u_{k-1})$ when $k>0$.
\end{enumerate}
\end{prop}

\begin{proof}
\begin{enumerate}
\item By Lemma \ref{lem:raact} the isomorphism is given by:
\begin{eqnarray*}
\bigoplus_{u\in U}\ho{k}{\tilde{C}_2^u\times\mu_2}{\Z}&\cong&\qsgr{A}\otimes_{\Z}\ho{k}{\tilde{C}_2\times\mu_2}{\Z}\\
(z_u)_{u\in U}&\leftrightarrow&\sum_{u\in U}\an{\odet(u)}\otimes u^{-1}\cdot z_u.\\
\end{eqnarray*}
\item Again the isomorphism is 
\begin{eqnarray*}
\bigoplus_{u\in U}\ho{k}{\tilde{S}^u\times\mu_2}{\Z}&\cong&\sgr{A}^{\pm}\otimes_{\Z}\ho{k}{\tilde{S}\times\mu_2}{\Z}\\
(z_u)_{u\in U}&\leftrightarrow&\sum_{u\in U}\an{\odet(u)}\otimes u^{-1}\cdot z_u.\\
\end{eqnarray*}
However, conjugation by $\omega$ acts as inversion on $C_3$ and hence acts on the $3$-Sylow subgroup of $\ho{k}{\tilde{S}\times\mu_2}{\Z}$ as $-1$ when $k\equiv 1\pmod{4}$. (See the proofs of Lemma \ref{lem:res} and Corollary \ref{cor:kerres}. )
\item  Combine part (2) with  Corollary \ref{cor:kerres}.
% conjugation by $\omega$ acts trivially on $\Ker{\Res{\tilde{C}_2\times\mu_2}{\tilde{S}\times\mu_2}}$ and hence on $\ho{k}{\spl{2}{A}}{Z^c}$. Thus these groups are $\qsgr{A}$-%modules. Lemma \ref{lem:raact} then shows that $v\in A^\times/\pm(A^\times)$ acts on 
%$\bigoplus_{u\in U}\ho{k}{C_3^u}{\Z}\cong \bigoplus_{u\in U}\Z/3\cong \bigoplus_{u\in A^\times/\pm(A^\times)^2}\Z/3$ as multiplication of the index by $v$.
\end{enumerate}
\end{proof}

In particular, we have:
\begin{cor}\label{cor:ccconstmod}  Let $A$ be a ring for which $\mu_2(A)$ is finite. Then 
%\begin{enumerate}
%\item 
$\ccconstmod{A}=\ho{0}{\spl{2}{A}}{Z^c}=\left(\qsgr{A}\otimes\Z/3\right)\cdot\cconst{A}$.
%\item $\ho{1}{\spl{2}{A}}{Z^c}=0$.
%\end{enumerate}
\end{cor}
%%%%%%%%%%%%%%%%%%%%%%%%%%%%%%%%%%%%%%%%%%%%%%%%%%%%%%%%%%%%%%%%%%%%%
\subsection{The homology of $Z^i$}
%%%%%%%%%%%%%%%%%%%%%%%%%%%%%%%%%%%%%%%%%%%%%%%%%%%%%%%%%%%%%%%%%%%%%%%%%%%%%%%

We begin by observing that the element $\omega\in \gl{2}{A}$ acts by conjugation on the diagonal subgroup $T_A\cong A^\times\times A^\times$ by interchanging factors. We let 
$D\subset \gl{2}{A}$ denote the resulting subgroup $T_A\rtimes C$. Let $\bar{D}$ denote the image of $D$ in $\pgl{2}{A}$: this is generated by $\omega$ and the subgroup
$\{ M_u:=\mathrm{diag}(u,1)\ |\ u\in A^\times\}\cong A^\times$ on which $\omega$ acts (by conjugation) as inversion.

\begin{prop}\label{prop:zises} For any ring $A$ there is a short exact sequence of right $\Z[\pgl{2}{A}]$-modules
\[
\xymatrix{
0\ar[r] &Z^i\ar^-{\gamma}[r]& \Ind{C_2}{\pgl{2}{A}}{\Z}\ar^-{\eta}[r]&\Ind{\bar{D}}{\pgl{2}{A}}{\Z}\ar[r]&0,
}
\]
where $\gamma(\suss{i}{u}):=1\otimes(M_u-I)$ for $u\in A^\times$ and $\eta(1\otimes g):=1\otimes g$.
\end{prop}

\begin{proof} In view of the presentation of $Z^i$ (as a $\pgl{2}{A}$-module) given in Lemma \ref{lem:zipres}, this follows from the observations that 
$(M_u-I)M_v= M_{uv}-M_v$ in $\Z[\pgl{2}{A}]$ and that 
\begin{eqnarray*}
(1\otimes (M_u-I))\cdot \omega&=&1\otimes \omega\cdot (M_{u^{-1}}-I)\\ 
&=&1\otimes (M_{u^{-1}}-I)\\
\end{eqnarray*}
in $\Z\otimes_{\Z[C_2]}\Z[\pgl{2}{A}]=\Ind{C_2}{\pgl{2}{A}}{\Z}$.
\end{proof}

%Given a group $G$ mapping to $\pgl{2}{A}$, let $\eta^G_n$ denote the map  
%\[
%\ho{n}{G}{\Ind{C_2}{\pgl{2}{A}}{\Z}}\to \ho{n}{G}{\Ind{\bar{D}}{\pgl{2}{A}}{\Z}}
%\]
% induced by 
%$\eta$. Thus the  long exact homology sequence associated to the short exact sequence of modules in Proposition \ref{prop:zises} gives a short exact sequence
%\[
%0\to \coker{\eta_{p+1}}\to H_p(G,Z^i)\to \ker{\eta_p}\to 0
%\]
%for all $p\geq 0$.
We immediately deduce:
\begin{cor}\label{cor:zises}  For any ring $A$ there is a short exact sequence of right $\Z[\gl{2}{A}]$-modules
\[
\xymatrix{
0\ar[r] &Z^i\ar^-{\gamma}[r]& \Ind{C_2\times Z}{\gl{2}{A}}{\Z}\ar^-{\eta}[r]&\Ind{{D}}{\gl{2}{A}}{\Z}\ar[r]&0,
}
\]
where $\gamma(\suss{i}{u}):=1\otimes(M_u-I)$ for $u\in A^\times$ and $\eta(1\otimes g):=1\otimes g$.
\end{cor}

Let $d:A^\times\to \spl{2}{A}$ be the embedding  $d(u):=\mathrm{diag}(u,u^{-1})$. Then $\tilde{D}:=D\cap\spl{2}{A}=d(A^\times)\rtimes \an{\tilde{\omega}}$
(where $\tilde{\omega}:=\omega\cdot M_{-1}\in\spl{2}{A}$). Since $D\backslash \gl{2}{A}/\spl{2}{A}$ is a singleton, we deduce that 
\[
\Ind{D}{\gl{2}{A}}{\Z}=\Ind{\tilde{D}}{\spl{2}{A}}{\Z}
\]
as a right $\Z[\spl{2}{A}]$-module.

As above, let $\epsilon:A^\times/\pm(A^\times)^2\to \gl{2}{A}$ be a section of $\oodet$ with image $U$ and such that $\epsilon(1)=1$ so that 
$\Ind{C_2\times Z}{\gl{2}{A}}{\Z}\cong \oplus_{u\in U}\Ind{\tilde{C}_2^u\times\mu_2}{\spl{2}{A}}{\Z}$ as a right $\Z[\spl{2}{A}]$-module. We thus have:

\begin{cor}\label{cor:zispl}  For any ring $A$ there is a short exact sequence of right $\Z[\spl{2}{A}]$-modules
\[
\xymatrix{
0\ar[r] &Z^i\ar^-{\gamma}[r]& \oplus_{u\in U}\Ind{\tilde{C}_2^u\times \mu_2}{\spl{2}{A}}{\Z}\ar^-{\eta}[r]&\Ind{\tilde{D}}{\spl{2}{A}}{\Z}\ar[r]&0.
}
\]
%where $\gamma(\suss{i}{u}):=1\otimes(M_u-I)$ for $u\in A^\times$ and $\eta(1\otimes g):=1\otimes g$.
\end{cor}

 Let $\eta_p$ denote the map  
\[
\oplus_{u\in U}\ho{p}{\tilde{C}_2^u\times \mu_2}{\Z}\cong \oplus_{u\in U}\ho{p}{\spl{2}{A}}{\Ind{\tilde{C}_2^u\times \mu_2}{\spl{2}{A}}{\Z}}\to \ho{p}{\spl{2}{A}}{\Ind{\tilde{D}}{\spl{2}{A}}{\Z}}
\cong\ho{p}{\tilde{D}}{\Z}
\]
 induced by 
$\eta$. Thus the  long exact homology sequence associated to the short exact sequence of modules in Corollary \ref{cor:zispl} gives a short exact sequence (of $\sgr{A}$-modules):
\begin{cor}\label{cor:zieta} For any ring $A$ there are natural short exact sequences of $\sgr{A}$-modules
\[
0\to \coker{\eta_{p+1}}\to \ho{p}{\spl{2}{A}}{Z^i}\to \ker{\eta_p}\to 0
\]
for all $p\geq 0$.
\end{cor}

We begin with $\ker{\eta_0}$:
For any ring $A$, let $\qaug{A}$ denote the augmentation ideal of the group ring $\qsgr{A}$.

\begin{lem}\label{lem:eta0} For any ring $A$, there is a $\sgr{A}$-module isomorphism $\ker{\eta_0}\cong\qaug{A}$ such that the (surjective) homomorphism 
\[
\xymatrix{
\cks{i}{A}:=\ho{0}{\spl{2}{A}}{Z^i}\ar^-{\gamma_0}[r]&\ker{\eta_0}\cong\qaug{A}\\
}
\]
sends $\suss{i}{v}$ to $\pf{v}\in \qaug{A}$.
\end{lem}

\begin{proof} The map $\oplus_{u\in U}\ho{0}{\tilde{C}_2\times\mu_2}{\Z}=\oplus_{u\in U}\Z\to \qsgr{A}$, $(1)_u\mapsto \an{u}$ is an isomorphism of $\sgr{A}$-modules by 
Lemma \ref{lem:raact}. The map resulting homomorphism
 \[
\qsgr{A}\cong \oplus_{u\in U}\ho{0}{\tilde{C}_2\times\mu_2}{\Z}\to\ho{0}{\tilde{D}}{\Z}=\Z
\]
 is clearly 
the augmentation map. 

Now, let $v\in A^\times$. Then $\mathrm{det}(M_v)=v$ and hence if $u':=\epsilon(M_v)\in U$ then $\an{v}=\an{u'}$ in $\qsgr{A}$ since $u'\equiv v \pmod{\pm (A^\times)^2}$.

The map $\eta$ sends $\suss{i}{v}\in Z^i$ to
 $1\otimes M_v-1\otimes I\in Z\otimes_{\Z[c_2\times Z]}\Z[\gl{2}{A}]$, which in turn corresponds to
 $(1\otimes \rho(M_v))_{u'}-(1\otimes I)_1\in \oplus_{u\in U}\Z\otimes_{\Z[\tilde{C}_2\times\mu_2]}\Z[\spl{2}{A}]$. Taking $\spl{2}{A}$-coinvariants, it follows that $\suss{i}{v}\in\cks{i}{A}$ maps to 
$\an{u'}-\an{1}=\pf{u'}=\pf{v}\in \qsgr{A}$ as required.
\end{proof}

Next we describe the (connecting) $\sgr{A}$-homomorphism $\delta_1:\ho{1}{\tilde{D}}{\Z}\to\ho{0}{\spl{2}{A}}{Z^i}:=\cks{i}{A}$:

\begin{lem}\label{lem:delta1} $\ho{1}{\tilde{D}}{\Z}=\gab{\tilde{D}}$ is a cyclic $\sgr{A}$-module with generator $\tilde{\omega}$. More precisely,
\[
\gab{\tilde{D}}\cong \frac{d(A^\times)\times \an{\tilde{\omega}}}{d((A^\times)^2)\cdot\an{d(-1)\tilde{\omega}^2}}
\]
and the $\sgr{A}$-action is determined by $\an{u}\cdot\tilde{\omega}=d(u)\cdot\tilde{\omega}$ and $\an{u}d(v)=d(v)$ for $u,v\in A^\times$.
%such that the $\sgr{A}$-homomorphism 
%$\rsgr{A}\to \gab{\tilde{D}}$ sending $1$ to $\tilde{\omega}$ is an isomorphism. 

The $\sgr{A}$-homomorphism $\delta_1$ is determined by $\delta_1(\tilde{\omega})=\suss{i}{-1}$. 
\end{lem}

\begin{proof}  Clearly $\gab{\tilde{D}}$ is as described since $d(u)$, $d(v)$ commute for any $u,v$ and 
\[
[d(u),\tilde{\omega}]=d(u)\cdot d(u^{-1})^{\tilde{\omega}}=d(u)\cdot d(u)=d(u^2)
\]
for any $u\in A^\times$. Furthermore, $\tilde{\omega}^2=-I=d(-1)=d(-1)^{-1}$.

For $u\in A^\times$, the action of $\an{u}$ on $\ho{1}{\tilde{D}}{\Z}=\gab{\tilde{D}}$ is induced from conjugation by $M_u:=\mathrm{diag}(u,1)$.
Such conjugation is trivial on $d(A^\times)$ while 
\[
\an{u}\cdot\tilde{\omega}=\tilde{\omega}^{M_u}=
\left[
\begin{array}{cc}
0&u^{-1}\\
-u&0\\
\end{array}
\right]
=\tilde{\omega}\cdot d(u).
\]
The first statement follows easily.

Now (the class of ) $\tilde{\omega}$ in $\gab{\tilde{D}}=\ho{1}{\tilde{D}}{\Z}$ is represented by the cycle 
\[
1\otimes (1,\tilde{\omega})\in \Z\otimes_{\Z[\tilde{D}]}C_1= \left(\Z\otimes_{\Z[\tilde{D}]}\Z[\spl{2}{A}]\right)\otimes_{\Z[\spl{2}{A}]}C_1=\Ind{\tilde{D}}{\spl{2}{A}}{\Z}\otimes_
{\Z[\spl{2}{A}]}C_1
\]
where here $C_i:=C_i(\spl{2}{A})$ is the standard homogeneous resolution of $\Z$ over $\spl{2}{A}$. This cycle in turn is represented by 
$(1\otimes 1)\otimes (1,\tilde{\omega})\in \Ind{D}{\gl{2}{A}}{\Z}\otimes C_1\cong \Ind{\tilde{D}}{\spl{2}{A}}{\Z}\otimes C_1$. It lifts to 
$(1\otimes 1)\otimes (1,\tilde{\omega})\in \Ind{C_2\times Z}{\gl{2}{A}}{\Z}\otimes C_1$ which maps under the boundary homomorphism to 
$(1\otimes 1)\otimes \left[ (\tilde{\omega})-(1)\right]\in \Ind{C_2\times Z}{\gl{2}{A}}{\Z}\otimes C_0$, and hence to $1\otimes \tilde{\omega}-1\otimes I$ 
under the isomorphism $\Ind{C_2\times Z}{\gl{2}{A}}{\Z}\otimes C_0\cong \Ind{C_2\times Z}{\gl{2}{A}}{\Z}$. Since $\tilde{\omega}=\omega\cdot M_{-1}$ and $\omega\in C_2\times Z$ this latter term is $1\otimes M_{-1}-1\otimes I=1\otimes (M_{-1}-I)$ which is the image of $\suss{i}{-1}$ under the map $Z^i\to\Ind{C_2\times Z}{\gl{2}{A}}{\Z}$. This establishes the second statement in the lemma.
\end{proof}

\begin{rem} Let $u\in A^\times$. Since $d(u)\in\gab{\tilde{D}}$ is equal to $\pf{u}\tilde{\omega}$, it follows that $\delta_1(d(u))=\pf{u}\suss{i}{-1} (=\suss{i}{u^2})$ in $\cks{i}{A}$.
\end{rem}

%Let $\rsgr{A}$ denote the ring $\F{2}[A^\times/(A^\times)^2]/\aug{A}^2$ which is a quotient ring of $\sgr{A}$, and hence a cyclic $\sgr{A}$-module.

Let $\tsgr{A}$ denote the quotient ring $\F{2}[A^\times/\mu_2(A)\cdot(A^\times)^2]/I^2$ of $\rsgr{A}$. 

\begin{rem} Let $G$ be a multiplicative elementary abelian $2$-group and let $L(G)$ be the isomorphic additive group with elements $\{ L(g)\ | g\in G\}$, obtained by replacing the operation of multiplication in $G$ with addition. If $I$ is the augmentation ideal if $\F{2}[G]$, then $I/I^2\cong L(G)$ via $\pf{g}\mapsto g$ and as a ring $\F{2}[G]/I^2\cong \F{2}\oplus L(G)$, where the second factor is a square zero ideal, via the map $\an{g}\mapsto 1+L(g)$.

In particular,  $\tsgr{A}\cong \F{2}\oplus L(A^\times/\mu_2(A)\cdot(A^\times)^2)$.
% $\rsgr{A}\cong\F{2}\oplus L( A^\times/(A^\times)^2)$ where the second term is a square zero ideal.
\end{rem}

%Thus
% $\tsgr{A}\cong \F{2}\oplus L(A^\times/\mu_2(A)\cdot(A^\times)^2)$.
\begin{lem}\label{lem:eta1} As an $\sgr{A}$-module
\[
\coker{\eta_1:\qsgr{A}\otimes\ho{1}{\tilde{C}_2\times\mu_2}{\Z}\to \ho{1}{\tilde{D}}{\Z}}\cong
\left\{
\begin{array}{ll}
\tsgr{A}\cdot\tilde{\omega},& -1\not\in (A^\times)^2\\
0,& -1\in (A^\times)^2.
\end{array}
\right.
\]
\end{lem}

\begin{proof} By Lemma \ref{lem:permmod}, the map composite map of $\Z[\spl{2}{A}]$-modules
\[
\xymatrix{
\Ind{\tilde{C}_2\times\mu_2}{\spl{2}{A}}{\Z}\ar[r]&\Ind{C_2\times Z}{\gl{2}{A}}{\Z}\ar^-{\eta}[r]&\Ind{D}{\gl{2}{A}}{\Z}\ar^-{\cong}[r]&\Ind{\tilde{D}}{\spl{2}{A}}{\Z}\\
}
\]
sends $1\otimes g$ to $1\otimes g$, for $g\in \spl{2}{A}$. Thus the induced map\\
 $\eta_k:\ho{k}{\tilde{C}_2\times\mu_2}{\Z}\to\ho{k}{\tilde{D}}{\Z}$ is corestriction. The corestriction map $\tilde{C}_2\times\mu_2=\ho{1}{\tilde{C}_2\times \mu_2}{\Z}\to
\ho{1}{\tilde{D}}{\Z}=\gab{\tilde{D}}$ sends $d(\zeta)$ to $d(\zeta)=\pf{\zeta}\tilde{\omega}\in \gab{\tilde{D}}$ for any $\zeta\in \mu_2$. If $-1\not\in (A^\times)^2$,  the image of $\eta_1$ consists of these elements -- since the $\sgr{A}$-action on the elements $d(u)\in\gab{\tilde{D}}$ is trivial --  and the cokernel is as described in the statement of the lemma.

On the other hand if $i^2=-1$ in $A^\times$ then $1\not= i\omega\in \tilde{C}_2$. Note that $i\omega=\tilde{\omega}^{M_{-i}}$ maps to $\an{-i}\tilde{\omega}\in\gab{\tilde{D}}$.
Since this latter generates the cyclic $\sgr{A}$-module $\gab{\tilde{D}}$, the result follows.
\end{proof}

Combining Lemma \ref{lem:eta1} with  Corollary \ref{cor:zieta} and Lemmas \ref{lem:eta0} and \ref{lem:delta1}, we  deduce:
\begin{prop}\label{prop:cksia}
\begin{enumerate}
\item  If $-1\in A^\times$,  then the map
\[
\cks{i}{A}:=H_0(\pspl{2}{A},Z^i)\to \aug{A}\subset\sgr{A},\quad \suss{i}{u}\mapsto \pf{u}
\]
is an isomomorphism of $\sgr{A}$-modules.
\item If $-1\not\in (A^\times)^2$ then  there is a natural short exact sequence of $\sgr{A}$-modules 
\[
\xymatrix{
0\ar[r]& \tsgr{A}\ar^-{f_1}[r]&\cks{i}{A}\ar[r]&\tilde{\aug{A}}\ar[r]& 0\\
}
\]
where $f_1(\alpha):=\alpha\cdot\suss{i}{-1}$ for $\alpha\in \tsgr{A}$ and the second map is $\suss{i}{v}\mapsto \pf{v}$.
\end{enumerate}
\end{prop}

\begin{rem} The short exact sequence of (2) is split as a sequence of $\Z$-modules since $\tilde{\aug{A}}$ is a free abelian group (with basis $\{ \pf{u}\ |\ u\in U\setminus\{ 1\} \}$ where 
$U$ is a set of representatives of $A^\times/\pm(A^\times)^2$).  However, it is not in general split as a sequence of $\sgr{A}$-modules.
\end{rem}

%%%%%%%%%%%%%%%%%%%%%%%%%%%%%%%%%%%%%%%%%%%%%%%%%%%%%%%%%%%%%%%%%%%%%
\subsection{The explicit structure of $\crpb{A}$}
%%%%%%%%%%%%%%%%%%%%%%%%%%%%%%%%%%%%%%%%%%%%%%%%%%%%%%%%%%%%%%%%%%%%%%%%%%%%

%Recall that
% \[
%\crpb{A}=(Z^c)_{\pspl{2}{A}}\oplus (Z^i)_{\pspl{2}{A}}=\ccconstmod{A}\oplus \cks{i}{A}  \mbox{ for }i=1,2.
%\]

For convenience, we collect together the conclusions of Corollary \ref{cor:ccconstmod}  and Proposition \ref{prop:cksia}. Recall that $\qsgr{A}:=\Z[A^\times/\pm(A^\times)^2]$ and $\tsgr{A}=\F{2}[H_A]/I^2\cong \F{2}\oplus L(H_A)$ where $H_A:=A^\times/\mu_2(A)\cdot (A^\times)^2$.

\begin{prop}\label{prop:crpbstruct}. Let $A$ be a ring. Then, for $i\in \{ 1,2\}$,  $\crpb{A}=\ccconstmod{A}\oplus\cks{i}{A}$ as an $\sgr{A}$-module, where
\begin{enumerate}
\item $\ccconstmod{A}$ is a free $\qsgr{A}\otimes \Z/3$-module with generator $\cconst{A}$ and,
\item There is a natural short exact sequence of $\sgr{A}$-modules
\[
0\to \sgr{A}\cdot\suss{i}{-1}\to\cks{i}{A}\to \tilde{\aug{A}}\to 0
\]
where 
\[
\sgr{A}\cdot\suss{i}{-1}=
\left\{
\begin{array}{ll}
\tsgr{A}\cdot\suss{i}{-1}\cong\tsgr{A}\cong\F{2}\oplus L(H_A),& -1\not\in (A^\times)^2\\
0,&-1\in (A^\times)^2.
\end{array}
\right.
\]

\end{enumerate}
\end{prop}
%By Corollary \ref{cor:zcpres}, $\ccconstmod{A}\cong \F{3}[A^\times/\pm (A^\times)^2]=\F{3}[\psqcl{A}]$ as an $\sgr{A}$-module.

%%%%%%%%%%%%%%%%%%%%%%%%%%%%%%%%%%%%%%%%%%%%%%%%%%%%%%%%%
\subsection{The field $\F{2}=\Z/2$ again}%  Example \F{2}
%%%%%%%%%%%%%%%%%%%%%%%%%%%%%%

 Take $A=\Z/2=\F{2}$.   Then $\wn{A}=\emptyset$ and $|Y(A)|$ is a contractible $2$-simplex, so that $H_2(|Y(A)|)=0$  and hence $\crpb{\F{2}}=\rpb{\F{2}}$ by Proposition \ref{prop:exacta}. 

We have $\cks{i}{\F{2}}=0$ and 
$\qsgr{\F{2}}=\Z/3$. Thus 
\[
\rpb{\F{2}}=\pb{\F{2}}=(\Z/3)\cdot \cconst{\Z/2}
\]
 is cyclic of order $3$ with generator $\cconst{\F{2}}$.

It easily follows that $\rbl{\F{2}}=\rpbker{\F{2}}=\rpb{\F{2}}$ since $\cconst{A}\in\rbl{A}$ always.

%%%%%%%%%%%%%%%%%%%%%%%%%%%%%%%%%%%%%%%%%%%
\subsection{The ring $\Z$ again}
%%%%%%%%%%%%%%%%%%%%%%%%%%%%%%%%%%
%\begin{exa}\label{exa:z}
 Since $\wn{\Z}=\emptyset$ we have $L_3(\Z)=0$. Furthermore $|Y(\Z)|$ is contractible and thus $\LL_3=Z_2(L_\bullet)= Z_2(\Sigma_\bullet)$ by Corollary \ref{cor:2dim}. 
The low-dimensional terms of the spectral sequence $E(\spl{2}{\Z},\LL)$ -- and hence the low-dimensional homology of $\spl{2}{\Z}$ -- are now easy to calculate (without appealing to any external calculations):
The top non-vanishing row of the $E^1$-page is $E^1_{p,3}=\ho{p}{\spl{2}{\Z}}{Z_2(L_\bullet)}=\ho{p}{\spl{2}{\Z}}{Z_2(\Sigma_\bullet)}$.

Thus $\rpb{\Z}=E^1_{0,3}=\ho{0}{\spl{2}{\Z}}{Z_2(\Sigma_\bullet)}=\crpb{\Z}$. We have $\qsgr{A}=\Z/3$, while $-1$ is not a square and $H_\Z=\{ 1\}$. 
We have $\cks{i}{\Z}=\ks{i}{\Z}=\Z/2\cdot \suss{1}{-1}$. Thus 
\[
\rpb{\Z}=\pb{\Z}=\Z/3\cdot\cconst{\Z}\oplus(\Z/2)\cdot\suss{1}{-1}
\]
 is cyclic of order $6$, with generator $\bconst{\Z}:=\suss{1}{-1}-\cconst{\Z}$.

Likewise $E^1_{1,3}=\ho{1}{\spl{2}{\Z}}{Z_2(\Sigma_\bullet)}=\ho{1}{\spl{2}{\Z}}{Z^c}\oplus\ho{1}{\spl{2}{\Z}}{Z^1}$. Now\\
 $\ho{1}{\spl{2}{\Z}}{Z^c}=0$ by  Proposition \ref{prop:hokzc}. Furthermore the group $\tilde{D}\subset\spl{2}{\Z}$ is cyclic of order $4$ with generator $\tilde{\omega}$. Thus the long exact homology sequence (arising from Corollary \ref{cor:zises})
\[
\cdots \to\ho{2}{\tilde{D}}{\Z}\to \ho{1}{\spl{2}{\Z}}{Z^1}\to\ho{1}{\mu_2}{\Z}\to\ho{1}{\tilde{D}}{\Z}\to\cdots 
\]
forces $\ho{1}{\spl{2}{\Z}}{Z^1}=0$. So $E^1_{1,3}=0$.

Note that  $d^1=\lambda_1:E^1_{0,3}\rpb{\Z}\to\sgr{\Z}=E^1_{0,2}$  is the zero map since $\lambda_1(\cconst{A})=\lambda_1(\suss{i}{-1})=0$ always. Thus $\rpb{\Z}=\rpbker{\Z}$ and the $E^2$-page of the spectral sequence loooks like:
\[
\xymatrix{
0&0&0&\cdots\\
\rpbker{\Z}\ar_-{d^2}[rdd]&0&E^2_{2,3}&\cdots\\
\aug{\Z}\ar_-{d^2}[rdd]&\aug{\Z}\otimes\ho{1}{\mu_2}{\Z}\quad\ar_-{d^2}[rdd]&0&\cdots\\
0&0&0&\cdots\\
\Z&B&\ho{2}{B}{\Z}\cong \Z/2&\ho{3}{\mu_2}{\Z}=\ho{3}{B}{\Z}\\
}
\]
Here $B$ -- the stabilizer of $\infty$ --  is the group of lower triangular matrices and is equal to $\mu_2\times \an{\rho}$ where $\rho:=
\left(
\begin{array}{cc}
1&0\\
1&1\\
\end{array}
\right)$ has infinite order. By the K\"unneth formula,\\
 $\ho{2}{B}{\Z}\cong\mu_2\otimes\an{\rho}$ is cyclic of order $2$ with generator $-I\otimes\rho$, and the inclusion 
$\mu_2\to B$ induces an isomorphism $\ho{3}{\mu_2}{\Z}\cong\ho{3}{B}{\Z}$.

$\aug{\Z}$ is infinite cyclic with generator $\pf{-1}$ and $d^2(\pf{-1})=-I\cdot \rho^{-6}\in B$ by \cite[Lemma A.6]{hut:bloch1}. It follows that $E^3_{0,2}=0$ and $E^3_{1,1}=E^\infty_{1,1}$ is cyclic of order $12$ with generator $\rho$. $\aug{\Z}\otimes\ho{1}{\mu_2}{\Z}$ is cyclic or order $2$ with generator $\pf{-1}\otimes -I$ and 
$d^2(\pf{-1}\otimes -I)=0$ by \cite[Corollary A.9]{hut:bloch1}. 

So the $E^3$-page of the spectral sequence has the form:
\[
\xymatrix{
0&0&0&\cdots\\
\rpbker{\Z}\ar^-{d^3}[rrddd]&0&E^3_{2,3}&\cdots\\
0&\aug{\Z}\otimes{\mu_2}&0&\cdots\\
0&0&0&\cdots\\
\Z&\Z/12&\mu_2\otimes\an{\rho}&\ho{3}{\mu_2}{\Z}\\
}
\]
Now $d^3(\cconst{\Z})=0$ and $d^3(\suss{i}{-1})=-I\otimes \rho$ (by \cite[Lemma A.12]{hut:bloch1}). It follows that $E^4_{0,3}=E^\infty_{0,3}=\rbl{\Z}$ is cyclic of order $3$ 
with generator $\cconst{\Z}$ and that $E^4_{2,0}=E^\infty_{2,0}=0$. Thus we can read off the low-dimensional homology groups: 

$\ho{1}{\spl{2}{\Z}}{\Z}\cong E^\infty_{1,1}\cong \Z/12$.

$\ho{2}{\spl{2}{\Z}}{\Z}=0$. 

Looking at the $E^\infty$ terms of total degree $3$, we see that $\ho{3}{\spl{2}{\Z}}{\Z}$ maps onto $\rbl{\Z}=\Z/3\cdot\cconst{\Z}$ with kernel of order $4$.
The inclusion $\mu_2\to \spl{2}{\Z}$ induces an injective map $E^\infty_{3,0}=\ho{3}{\mu_2}{\Z}\to\ho{3}{\spl{2}{\Z}}{\Z}$. Since $\tilde{\omega}^2=-I\in\mu_2$ it follows that the natural map $\Z/4\cong\ho{3}{\an{\tilde{\omega}}}{\Z}\to\ho{3}{\spl{2}{\Z}}{\Z}$ is injective and hence its image is the kernel of the map to $\rbl{\Z}$. We conclude:
\begin{prop}\label{prop:h3sl2z} There is a natural short exact sequence 
\[
0\to\ho{3}{\an{\tilde{\omega}}}{\Z}\to \ho{3}{\spl{2}{\Z}}{\Z}\to\rbl{\Z}\to 0.
\]
\end{prop}
\begin{rem} The sequence is naturally split since there is a distinguished element of order $3$ in $\ho{3}{\spl{2}{\Z}}{\Z}$, also denoted $\cconst{\Z}$, mapping to $\cconst{\Z}\in\rbl{\Z}$
(\cite[Proposition 7.3]{hut:bloch1}).
\end{rem}
%\end{exa}
%\begin{rem} In \cite[Lemma 8.13]{hut:bloch1} it is shown however that $\suss{1}{-1}\not\in \rbl{\Z}$ (or $\bl{\Z}$). This fact also follows independently from Corollary \ref{cor:suss1z4} below.
%It follows again that $\rbl{\Z}=\bl{\Z}=\Z/3\cdot\cconst{\Z}$ is cyclic of order $3$.
%\end{rem}
%%%%%%%%%%%%%%%%%%%%%%%%%%%%%%%%%%%%%%%%%%%%%%%%%%%%%%
\subsection{The map $\crpb{A}\to \rpb{A}$}
%%%%%%%%%%%%%%%%%%%%%%%%%%%%%%%%%%%%%%%%%%%%%%%%%%%%%%%%%%%%%
For any ring $A$, the natural inclusion homomorphism $Z_2(\Sigma_\bullet)\to Z_2(L_\bullet)$ induces a natural homomorphism
 \[
\crpb{A}=(Z_2(\Sigma_\bullet))_{\pspl{2}{A}}\to 
(Z_2(L_\bullet))_{\pspl{2}{A}}=\rpb{A}.
\]
Since $\crpb{A}=\ccconstmod{A}\oplus \cks{i}{A}$ it easily follows that the image of this homomorphism is the submodule $\cconstmod{A}\oplus \ks{i}{A}$ of $\rpb{A}$ generated by the element $\cconst{A}$ and the elements $\suss{i}{u}$, $u\in A^\times$ (where $i=1$ or $2$; recall that $\pf{u}\cconst{A}=\suss{1}{u}-\suss{2}{u}$.).

Recall that the inclusion map $Z_2(L_\bullet)\to L_2(A)$ induces a natural homomorphism $\lambda_1:\rpb{A}=Z_2(L_\bullet)_{\pspl{2}{A}}\to L_2(A)_{\pspl{2}{A}}=\sgr{A}$. Thus there is a compatible homomorphism $\lambda_1:\crpb{A}\to \sgr{A}$ induced by the inclusion $Z_2(\Sigma_\bullet)\to L_2(A)$.  It is easily verified that $\lambda_1(\cconst{A})=0$ and $\lambda_1(\suss{i}{u})=(1+\an{-1})\pf{u}$. Thus if we let $\crpbker{A}:=\ker{\lambda_1:\crpb{A}\to \sgr{A}}$, we have a commutative diagram with exact rows
\[
\xymatrix{
0\ar[r]&\crpbker{A}\ar[r]\ar[d]&\crpb{A}\ar^-{\lambda_1}[r]\ar[d]&(1+\an{-1})\aug{A}\ar[r]\ar[d]&0\\
0\ar[r]&\rpbker{A}\ar[r]&\rpb{A}\ar^-{\lambda_1}[r]&\sgr{A}&\\
}
\]
where the righmost vertical arrow is an inclusion map. Thus $\ker{\crpb{A}\to\rpb{A}}=\ker{\crpbker{A}\to\rpbker{A}}$.  Furthermore, note that the map of $\sgr{A}$-modules
$\qsgr{A}\to\sgr{A}$, $x\mapsto (1+\an{-1})x$ is well-defined and induces an isomorphism $\tilde{\aug{A}}\cong (1+\an{-1})\aug{A}$. We deduce:
\begin{prop}\label{prop:crpbker}
 Let $A$ be a ring such that $\mu_2(A)$ is finite. Then
 \[
\crpbker{A}=\ccconstmod{A}\oplus \sgr{A}\cdot\suss{i}{-1}\cong
\left\{
\begin{array}{ll}
\left(\qsgr{A}\otimes\Z/3\right)\cdot\cconst{A}\oplus \tsgr{A}\cdot\suss{i}{-1},& -1\not\in (A^\times)^2\\
\left(\qsgr{A}\otimes\Z/3\right)\cdot\cconst{A}, &-1\in (A^\times)^2.\\
\end{array}
\right.
\]
\end{prop}

In general, the map $\crpbker{A}\to \rpbker{A}$ is not injective. Here are some examples:

\begin{exa}  Let $A$ be any ring for which $\cconst{A}=0$; eg., let $A$ be a field in which the quadratic $x^2-x+1$ has a root. Then $\ccconstmod{A}\not=0$ but the 
map $\ccconstmod{A}\to \rpb{A}$ is the zero map. Thus $0\not= \ccconstmod{A}\subset \ker{\crpb{A}}\to \rpb{A}$  for such a ring.

\end{exa}

\begin{exa} Consider the ring $A=\Q$. Then \[
\ccconstmod{\Q}\cong \Z/3[\Q^\times/\pm(\Q^\times)^2]\cong\Z/3\left[\bigoplus p^{\Z/2}\right]
\]
 where the sum is over all primes. However, by \cite[Lemma 6.5]{hut:sl2Q}, if we let $E:=\Q(\sqrt{-3})$ then
\[
\cconstmod{\Q}\cong\Z/3[\Q^\times/\pm N_{E/\Q}(E^\times)]\cong  \Z/3\left[\bigoplus_{p=2,p\equiv -1\mod 3}p^{\Z/2}\right],
\]
and hence the homomorphism $\ccconstmod{\Q}\to\cconstmod{\Q}\subset \rpb{\Q}$ has an infinite kernel.
\end{exa}
\begin{exa}  Let $A$ be a ring satisfying 
\begin{enumerate}
\item $-1$ is not a square in $A^\times$, and
\item there exists $u\in A^\times$ such that $w:=1+u^2\in A^\times$ but $w^2\not\in (A^\times)^4$.
\end{enumerate}
(For example, $A=\Z[1/2]$ and $u=-1$, or $A=\Q$ and $u\not=0$ with $1+u^2\not\in (\Q^\times)^2$.)

Then $\pf{w}\suss{i}{-1}\not=0$ in $\sgr{A}\cdot\suss{i}{-1}\subset\crpb{A}$ by Proposition \ref{prop:crpbstruct}, but $\pf{w}\suss{i}{-1}=0$ in $\rpb{A}$ by \cite[Corollary 3.22]{hut:bloch1}.
\end{exa}
\begin{exa} Below in section \ref{sec:z8} it is shown that $\pf{3}\cconst{\Z/8}\in \crpb{\Z/8}$ has order $3$ while its image in $\rpb{\Z/8}$ is zero.
\end{exa}
%%%%%%%%%%%%%%%%%%%%%%%%%%%%%%%%%%%%%%%%%%%%%%%%%%%%%%%%%%%%%%%%%%%%%%%%
\section{Local rings with residue field $\F{2}$}\label{sec:f2}
%%%%%%%%%%%%%%%%%%%%%%%%%%%%%%%%%%%%%%%%%%%%%%%%%%%%%%%%%
Let $A$ be a local ring with maximal ideal $M\not= 0$ and residue field $A/M=\F{2}$. Let $Y=Y(A)$.

%%%%%%%%%%%%%%%%%%%%%%%%%%%%%%%%%%%%%%%%%%%%%%%%%%%%%%%%%%%%%%%%%%%
\subsection{Generators of $\rpb{A}$}
%%%%%%%%%%%%%%%%%%%%%%%%%%%%%%%%%%%%%%%%%%%%%%%%%%%%%%%%%%%%%%%%%%%%%%%%%%
By Corollary \ref{cor:2dim} we have a short exact sequence
\[
0\to Z_2(\Sigma_\bullet)\to Z_2(L_\bullet)\to H_2(|Y|)\to 0.
\]
Furthermore, Proposition \ref{prop:localhom} gives an explicit basis of the free $\Z$-module $H_2(|Y|)$ represented by elements $\gamma$ which lie in $Z_2(L_\bullet)$. In particular, we can write down an explicit $\Z$-splitting of this short exact sequence. 

By Proposition \ref{prop:localhom}, elements of $Z_2(L_\bullet)$ representing generators of $H_2(|Y|)$ can be obtained as follows:

We let $T_0:=M=\{ (m,1)\in Y_0\ |\  m\in M\}\subset Y_0$, $T_1:=A^\times=\{ (u,1)\in Y_0\  |\  u\in A^\times\}$ and 
$T_\infty:= 1/M:= \{ (1,m)\in Y_0\ |\ m\in M\}$. Thus $Y_0=T_\infty\sqcup T_0\sqcup T_1$.

 Let $\oct:= (t_{i,j})$, $i=1,2$, $j=1,2,3$ be a $2\times 3$ matrix with entries 
in $Y_0$ satisfying the following conditions:
\begin{enumerate}
\item There exists a bijection $\phi=\phi_\oct:\{ 1,2,3\}\to\{ \infty, 0,1\}$ such that $t_{i,j}\in T_{\phi(j)}$ for all $i,j$, \and
\item $t_{1,j}\not= t_{2,j}$  for $j=1,2,3$.
\end{enumerate}
Thus the entries of $\oct$ are the vertices of an octahedron in the space $Y$ whose surface represents a top-dimensional cycle. The pairs $\{ t_{1,j},t_{2,j}\}$ are antipodal points in the octahedron (and are not edges in $Y$). The corresponding homology class in $H_2(|Y|)$ is represented by the cycle
\[
\cyc{\oct}=\cyc{t_{i,j}}:= \sum_{i=(i_1,i_2,i_3)\in \{ 1,2\}^3}(-1)^{\delta(i)}(t_{i_1,1},t_{i_2,2},t_{i_3,3})\in Z_2(L_\bullet)
\]
where $\delta(i)=\sum_{r=1}^3 i_r\pmod{2}$. 

By Proposition \ref{prop:localhom}, if we fix the bijection $\phi:\{1,2,3\}\to \{\infty,0,1\}$ and row $1$ of the octahedron $\oct$, then as the second row varies among all possible 
values,  the family $\{ \cyc{\oct}\}$ in $Z_2(L_\bullet)$  maps to a basis of $H_2(|Y|)$.

We will call an octahedron $\oct$ \emph{standard} if $\phi_\oct=\Phi$ where $\Phi(1)=\infty$, $\Phi(2)=0$ and $\Phi(3)=1$. In general, we let $\sigma_\oct\in S_3$ be the permutation
$\Phi^{-1}\circ \phi_\oct$. Thus $\oct$ is standard if and only if $\sigma_\oct$ is the identity.

We will say that $\oct$ is \emph{normalized} if $t_{1,j}=\phi(j)$ for all $j$. We will say that $\oct$ is \emph{basic} if it is both standard and normalized. 
Thus $\oct$ is basic if and only if the top row 
is $[ \infty,0,1]$. In particular, the family $\{ \cyc{\oct}\}_{\oct \mbox{ basic}}$ maps to a basis of $H_2(|Y|)$ which we call the \emph{standard basis}.

Thus, we can let $H(A)$ denote the free additive subgroup of $Z_2(L_\bullet)$ with basis the family $\{ \cyc{\oct}\}_{\oct \mbox{ basic}}$. By Corollary \ref{cor:2dim} 
we deduce:

\begin{prop}\label{prop:localf2} Let $A$ be a local ring with residue field $\F{2}$. Then $Z_2(L_\bullet)=Z_2(\Sigma_\bullet)\oplus H(A)$ where $H(A)$ is free with basis 
$\{ \cyc{\oct}\}_{\oct \mbox{ basic}}$. 

In particular, $Z_2(L_\bullet)$ is generated as an additive group by the elements $D_1^{x,y}(a,b)$, $C^{x,y}(a)$ and $\{ \cyc{\oct}\}_{\oct \mbox{ basic}}$.
\end{prop}

\begin{cor}\label{cor:localf2} Let $A$ be a local ring with residue field $\F{2}$. Then $\rpb{A}$ is generated as an abelian group by the elements
$\an{a}\suss{1}{u}$, $\an{a}\cconst{A}$ and $\cyc{\oct}$ for $a,u\in A^\times$, and $\oct$ a basic octahedron
\end{cor}

For convenience, we will extend the above notations in the following ways below:
\begin{enumerate}
\item  We will allow degenerate octahedra $\oct$, where $t_{1,j}=t_{2,j}$ for at least one $j$. In such a case, the corresponding cycle $\cyc{\oct}$ is equal to $0$.
\item Since the first row of a basic octahedron is determined, we will denote basic octahedra by their second row $(1/m, n,u)$ where $m,n\in M$ and $u\in A^\times$, and we will denote the corresponding cycle by $\scyc{1/m\ n\ u}$.
\end{enumerate}

%%%%%%%%%%%%%%%%%%%%%%%%%%%%%%%%%%%%%%%%%%%%%%%%%%%%%%%%%%%%%%%%%%%%%%%%
\subsection{A complete family of relations}
%%%%%%%%%%%%%%%%%%%%%%%%%%%%%%%%%%%%%%%%%%%%%%%%%%%%%%%%%%%%%%%%%%%%%%%%%%%%

In order to calculate the action of $\pgl{2}{A}$ on (cycles associated to) octahedra below, however, we first need to calculate some algebraic relations among the $\cyc{\oct}$ for general 
octahedra $\oct$.  The following lemma is straightforward to verify from the defining formula for the cycles $\cyc{\oct}$:

\begin{lem}\label{lem:oct1} Fix an octahedron ${\oct}=( t_{i,j})$ and let $\phi=\phi_\oct$. Given $j\in \{ 1,2,3\}$ and $s_{i,j}\in T_{\phi(j)}$ for $i=1,2$ we let
$
\cyc{\cdots
\begin{array}{c}
s_{1,j}\\
s_{2,j}\\
\end{array}
\cdots
}
$
denote the cycle associated to the octahadron obtained from $\oct$ by replacing column $j$ with the given column. 

\begin{enumerate}
\item
\[
\cyc{\cdots
\begin{array}{c}
t_{2,j}\\
t_{1,j}\\
\end{array}
\cdots
}=
-\cyc{\cdots
\begin{array}{c}
t_{1,j}\\
t_{2,j}\\
\end{array}
\cdots
}=-\cyc{\oct} \mbox{ in } Z_2(L_\bullet).
\]

\item For any $s_{1,j}\not=t_{1,j}\in T_{\phi(j)}$ we have 
\begin{eqnarray*}
\cyc{\cdots
\begin{array}{c}
s_{1,j}\\
t_{2,j}\\
\end{array}
\cdots
}
&=&
\cyc{\cdots
\begin{array}{c}
t_{1,j}\\
t_{2,j}\\
\end{array}
\cdots}
-
\cyc{\cdots
\begin{array}{c}
t_{1,j}\\
s_{1,j}\\
\end{array}
\cdots
}\\
&=&
\cyc{\oct}
-
\cyc{\cdots
\begin{array}{c}
t_{1,j}\\
s_{1,j}\\
\end{array}
\cdots
}\mbox{ in } Z_2(L_\bullet).\\
\end{eqnarray*}
\end{enumerate}
\end{lem}

\begin{cor} 
 If $\oct$ is a standard octahedron, then $\cyc{\oct}\in Z_2(L_\bullet)$ is an integer linear combination of basic octahedra; i.e., $\cyc{\oct}\in H(A)$.
\end{cor}

Now, for $\sigma\in S_3$ (acting on the left on $\{ 1,2,3\}$) we let $\sigma$ act on octahedra by permuting the columns: $\sigma\cyc{ t_{i,j}}:= \cyc{ t_{i,\sigma^{-1}(j)}}$. In this way $S_3$ 
acts on the elements $\cyc{\oct}\in Z_2(L_\bullet)$. 

Now for any $\sigma\in S_3$ and any octahedron $\oct$, let 
\[
R_\sigma(\cyc{\oct}):= \cyc{\oct}-\mathrm{sgn}(\sigma)\sigma\cdot \cyc{\oct}\in Z_2(\Sigma_\bullet)\subset Z_2(L_\bullet).
\]

\begin{lem}\label{lem:oct2} 
\begin{enumerate}

%\item For any $\sigma\in S_3$ and for any octahedron $\oct$, $R_\sigma(\cyc{\oct})\in Z_2(\Sigma_\bullet)$.
\item For any octahedron $\oct=(t_{i,j})$, $R_{(12)}(\cyc{\oct})=D_1(\oct)$ where 
\[
D_1(\oct):= \sum_{(i_1,i_2)\in \{ 1,2\}^2}(-1)^{i_1+i_2}D_1^{\{ t_{i_1,1},t_{1_2,2}\}}(t_{1,3},t_{2,3})\in Z_2(\Sigma_\bullet).
\]
\item For any octahedron $\oct=(t_{i,j})$, $R_{(23)}(\cyc{\oct})=D_2(\oct)$ where 
\[
D_2(\oct):= \sum_{(i_2,i_3)\in \{ 1,2\}^2}(-1)^{i_2+i_3}D_2^{\{ t_{i_2,2},t_{1_3,3}\}}(t_{1,1},t_{2,1})\in Z_2(\Sigma_\bullet).
\]
\end{enumerate}
\end{lem}

\begin{proof}\ \\
\begin{enumerate}
%\item For any $\sigma\in S_3$, $R_\sigma(\cyc{\oct})\in \Sigma_2\cap Z_2(L_\bullet)=Z_2(\Sigma_\bullet)$.
\item We have
\begin{eqnarray*}
\cyc{\oct}&=&  \sum_{(i_1,i_2)\in \{ 1,2\}^2}(-1)^{i_1+i_2}\left[ (t_{1_1,1},t_{i_2,2},t_{2,3})-(t_{1_1,1},t_{i_2,2},t_{1,3})\right]\\
\end{eqnarray*}
and hence
\begin{eqnarray*}
R_{(12)}(\cyc{\oct})=\cyc{\oct}+(12)\cdot\cyc{\oct}&=&\sum_{(i_1,i_2)\in \{ 1,2\}^2}(-1)^{i_1+i_2}\left[R_1^{\{ t_{i_1,1},t_{1_2,2}\}}(t_{2,3})-R_1^{\{ t_{i_1,1},t_{1_2,2}\}}(t_{1,3})\right]\\
&=&\sum_{(i_1,i_2)\in \{ 1,2\}^2}(-1)^{i_1+i_2} D_1^{\{ t_{i_1,1},t_{1_2,2}\}}(t_{1,3},t_{2,3})=D_1(\oct).\\
\end{eqnarray*}
\item Entirely analogous to (1).

\end{enumerate}
\end{proof}
\begin{rem}\label{rem:13}
 From the definition, we see that for all $\sigma,\tau\in S_3$ and for all octahedra $\oct$, we have
\[
R_{\sigma\tau}(\cyc{\oct})=R_{\tau}(\cyc{\oct})+\mathrm{sgn}(\tau)R_{\sigma}(\tau\cdot\cyc{\oct}).
\]
Since $(12)$ and $(23)$ generate $S_3$ and since $R_{(12)}(\cyc{\oct})$ and $R_{(23)}(\cyc{\oct})$ are known for all octahedra $\oct$, we can thus calculate $R_{\sigma}(\cyc{\oct})$ for any 
$\sigma\in S_3$.

For example, since $(13)=(12)(23)(12)=(12)(132)$ in $S_3$ we obtain, using the formula of (3) twice:
\begin{eqnarray*}
R_{(13)}(\cyc{\oct})&=& R_{(23)(12)}(\cyc{\oct})+R_{(12)}( (132)\cdot\cyc{\oct})\\
&=& R_{(12)}(\cyc{\oct})-R_{(23)}((12)\cdot\cyc{\oct})+R_{(12)}((132)\cdot\cyc{\oct})\\
&=& D_1(\oct)-D_2((12)\cdot\oct)+D_1((132)\cdot\oct)\in Z_2(\Sigma_\bullet).
\end{eqnarray*}
\end{rem}
Since $Z_2(\Sigma_\bullet)$ is the kernel of the map $Z_2(L_\bullet)\to H_2(|Y|)$ we deduce:
\begin{cor} In $H_2(|Y|)$ we have $\sigma\cdot\cyc{\oct}=\mathrm{sgn}(\sigma)\cyc{\oct}$ for all permutations $\sigma\in S_3$ and for all octahadra $\oct$. 
\end{cor}

We are now in a position to calculate the action of $\pgl{2}{A}$ on the (cycles associated to the) basic octahedra in $Z_2(L_\bullet)$:

Note that since $A$ is a local ring with residue field $\F{2}$, $\pspl{2}{A}$ is generated by the elementary matrices $E_{12}(w)$ and $E_{21}(w)$ where $w\in A^\times$ runs over any set of units which additively generate $A$. Thus $\pgl{2}{A}$ is generated by these elementary matices and the matrices $M_w:=\mathrm{diag}(w,1)$, where $w$ run over any set of generators of the group $A^\times$.
\begin{lem}\label{lem:action} Let $\scyc{ 1/m\  n\  u}$ be a basic octahedron in $Z_2(L_\bullet)$.

\begin{enumerate}
\item For any $w\in A^\times$ we have $\scyc{ 1/m\  n\ u}\cdot M_w=\scyc{w/m\ wn\ wu}-\scyc{ w/m\ wn\ w}$.
\item For any $w\in A^\times$, 
\begin{eqnarray*}
\scyc{ 1/m\  n\ u}\cdot E_{12}(w)&=& \scyc{\frac{1}{1+w}\ \frac{n}{1+nw}\ \frac{1}{w+m}}-\scyc{\frac{1}{1+w}\ \frac{n}{1+nw}\ \frac{1}{w}} \\
&&+\scyc{\frac{u}{1+uw}\ \frac{n}{1+nw}\ \frac{1}{w}}-\scyc{\frac{u}{1+uw}\ \frac{n}{1+nw}\ \frac{1}{w+m}}+R_{(13)}(\cyc{\oct})\\
\end{eqnarray*}
where 
\[
\oct=\left(
\begin{array}{lll}
\frac{1}{1+w}& 0&\frac{1}{w}\\
\frac{u}{1+uw}&\frac{n}{1+nw}&\frac{1}{w+m}\\
\end{array}
\right).
\]
\item For any $w\in A^\times$,
\begin{eqnarray*}
 \scyc{ 1/m\  n\ u}\cdot E_{21}(w)&=&    \scyc{\frac{1+mw}{m}\ u+w\ w}-\scyc{\frac{1+mw}{m}\ u+w\ n+w} \\
&&+\scyc{\frac{1+mw}{m}\ 1+w\ n+w}-\scyc{\frac{1+mw}{m}\ 1+w\ w}+R_{(23)}(\cyc{\oct'})\\    \\
\end{eqnarray*}
where 
\[
\oct'=\left(
\begin{array}{lll}
\infty& 1+w&{w}\\
\frac{1+mw}{m}&u+w&{n+w}\\
\end{array}
\right).
\]
\end{enumerate}
\end{lem}

\begin{proof}
\begin{enumerate}
\item We have 
\begin{eqnarray*}
\left(
\begin{array}{lll}
\infty&0&1\\
\frac{1}{m}&n&u\\
\end{array}
\right)\cdot M_w&=&
\left(
\begin{array}{lll}
\infty&0&w\\
\frac{w}{m}&wn&wu\\
\end{array}
\right)\\
&=& \left(\frac{w}{m}\ wn\ wu\right)-\left(\frac{w}{m}\ wn\ w\right)\\
\end{eqnarray*}
using Lemma \ref{lem:oct1} (2) in the second line.
\item We have
\[
\left(
\begin{array}{lll}
\infty&0&1\\
\frac{1}{m}&n&u\\
\end{array}
\right)
\cdot E_{12}(w)=
\left(
\begin{array}{lll}
\frac{1}{w}&0&\frac{1}{1+w}\\
\frac{1}{w+m}&\frac{n}{1+nw}&\frac{u}{1+uw}\\
\end{array}
\right)
=(1 3)\cdot \oct.
\]
By definition of $R_{(13)}(\cyc{\oct})$ we have $(13)\cdot \cyc{\oct} = -\cyc{\oct}+ R_{(13)}(\cyc{\oct})$.

 Finally, we have 
\begin{eqnarray*}
\cyc{\oct}&=& \cyc{
\begin{array}{lll}
\frac{1}{1+w}& 0&\frac{1}{w}\\
\frac{u}{1+uw}&\frac{n}{1+nw}&\frac{1}{w+m}\\
\end{array}
}\\
&=&  \cyc{
\begin{array}{lll}
\frac{1}{1+w}& 0&1\\
\frac{u}{1+uw}&\frac{n}{1+nw}&\frac{1}{w+m}\\
\end{array}
}-
 \cyc{
\begin{array}{lll}
\frac{1}{1+w}& 0&1\\
\frac{u}{1+uw}&\frac{n}{1+nw}&\frac{1}{w}\\
\end{array}
}\\
&=& \left(\scyc{\frac{u}{1+uw}\ \frac{n}{1+nw}\ \frac{1}{w+m}}-\scyc{\frac{1}{1+w}\ \frac{n}{1+nw}\ \frac{1}{w+m}}\right)\\
&&- \left( \scyc{\frac{u}{1+uw}\ \frac{n}{1+nw}\ \frac{1}{w}}-\scyc{\frac{1}{1+w}\ \frac{n}{1+nw}\ \frac{1}{w}} \right)
\end{eqnarray*}
using Lemma \ref{lem:oct1} (2) twice.
\item This identity is arrived at in an entirely analogous way as  (2), since
\[
\left(
\begin{array}{lll}
\infty&0&1\\
\frac{1}{m}&n&u\\
\end{array}
\right)
\cdot E_{21}(w)=
\left(
\begin{array}{lll}
\infty&w&{1+w}\\
\frac{1+mw}{m}&n+w&u+w\\
\end{array}
\right)
=(2 3)\cdot \oct'.
\]
\end{enumerate}
\end{proof}

Recall that when  $A$ is a local ring with residue field $\F{2}$,
\[
Z_2(L_\bullet)=Z_2(\Sigma_\bullet)\oplus\left( \bigoplus_{\oct \mbox{ basic }}\Z\cdot\scyc{\oct}\right)=Z^c\oplus Z^i\oplus \left( \bigoplus_{\oct \mbox{ basic }}\Z\cdot\scyc{\oct}\right)
\mbox{($i=1$ or $2$)}
\]
as an abelian group. Thus combining the $\pgl{2}{A}$-presentation of $Z^c$ and $Z^i$  (Lemmas \ref{lem:zcpres} and \ref{lem:zipres}) with Lemma \ref{lem:key} and Lemma \ref{lem:action},
we obtain the following $\Z[\pgl{2}{A}]$-presentation of $Z_2(L_\bullet)$:

\begin{prop}\label{prop:lz2lpres}
 Let $A$ be a local ring with maximal ideal $M$ and residue field $\F{2}$. Then $Z_2(L_\bullet)$ has the following presentation as a $\Z[\pgl{2}{A}]$-module:

\begin{enumerate}
\item[Generators:]

$\cconst{A}$, $\Suss{i}{u}$, $i\in \{ 1,2\}$, $u\in A^\times$, $\scyc{ 1/m\ n\  u}$, $m,n\in M$, $u\in A^\times$.

\item[Relations:]  Let $W\subset A^\times$ be any subset which additively generates $A$.
\begin{enumerate}
\item For $i\in \{ 1,2\}$, $\Suss{i}{uv}=\Suss{i}{v}\cdot M_u+\Suss{i}{u}$ for all $u,v\in A^\times$,
\item $\cconst{A}\cdot\omega=\cconst{A}$ and $\Suss{i}{u}\cdot\omega=\Suss{i}{u^{-1}}$ for $i\in \{ 1,2\}$ $u\in A^\times$,
\item $\cconst{A}\cdot (M_u-I)=\Suss{1}{u}-\Suss{2}{u}$ for all $u\in A^\times$,
\item $\cconst{A}\cdot(I+M_{0,1,\infty}+M_{1,\infty,0})=0$,
\item For any $w\in W$,
\begin{eqnarray*}
\scyc{ 1/m\  n\ u}\cdot E_{12}(w)&=& \scyc{\frac{1}{1+w}\ \frac{n}{1+nw}\ \frac{1}{w+m}}-\scyc{\frac{1}{1+w}\ \frac{n}{1+nw}\ \frac{1}{w}} \\
&&+\scyc{\frac{u}{1+uw}\ \frac{n}{1+nw}\ \frac{1}{w}}-\scyc{\frac{u}{1+uw}\ \frac{n}{1+nw}\ \frac{1}{w+m}}\\
&&+ D_1(\oct)-D_2((12)\cdot\oct)+D_1((132)\cdot\oct)\in Z_2(\Sigma_\bullet)
%\suss{1}{\oct}-\suss{2}{(12)\oct}+\suss{1}{(13)\oct}
\end{eqnarray*}
where 
\[
\oct=\left(
\begin{array}{lll}
\frac{1}{1+w}& 0&\frac{1}{w}\\
\frac{u}{1+uw}&\frac{n}{1+nw}&\frac{1}{w+m}\\
\end{array}
\right).
\]
\item For any $w\in W$,
\begin{eqnarray*}
 \scyc{ 1/m\  n\ u}\cdot E_{21}(w)&=&    \scyc{\frac{1+mw}{m}\ u+w\ w}-\scyc{\frac{1+mw}{m}\ u+w\ n+w} \\
&&+\scyc{\frac{1+mw}{m}\ 1+w\ n+w}-\scyc{\frac{1+mw}{m}\ 1+w\ w}\\    
&&+D_2(\oct')\\
\end{eqnarray*}
where 
\[
\oct'=\left(
\begin{array}{lll}
\infty& 1+w&{w}\\
\frac{1+mw}{m}&u+w&{n+w}\\
\end{array}
\right).
\]
\item For any $w\in A^\times$,
\[
\scyc{ 1/m\  n\ u}\cdot M_w=\scyc{w/m\ wn\ wu}-\scyc{ w/m\ wn\ w}.
\]
\end{enumerate}
\end{enumerate}
\end{prop}

For $i=1,2$ and for any octahedron $\oct$, we let $\suss{i}{\oct}$ denote the image of $D_i(\oct)$ under the composite homomorphism $Z_2(\Sigma_\bullet)\to \crpb{A}\to\rpb{A}$.  Recall that 
\[
D_i^{\{ x,y\}}(a,b)\mapsto \an{d(x,y,a)}\suss{i}{\frac{d(x,a)d(y,b)}{d(x,b)d(y,a)}}
\]
under this composite, where $d(x,y,a):=d(x,y)d(y,a)d(a,x)\mod{(A^\times)^2}$. Thus, when $\oct =(t_{i,j})$, we have 
\begin{eqnarray*}
\suss{1}{\oct}&=&\sum_{(i_1,i_2)\in \{ 1,2\}^2}(-1)^{i_1+i_2}\an{d(t_{i_1,1},t_{i_2,2},t_{1,3})}\suss{1}{\frac{d(t_{i_1,1},t_{1,3})d(t_{i_2,2},t_{2,3})}{d(t_{i_1,1},t_{2,3})d(t_{i_2,2},t_{1,3})}}\\
\suss{2}{\oct}&=&\sum_{(i_2,i_3)\in \{ 1,2\}^2}(-1)^{i_2+i_3}\an{d(t_{i_2,2},t_{i_3,3},t_{1,1})}\suss{2}{\frac{d(t_{i_2,2},t_{1,1})d(t_{i_3,3},t_{2,1})}{d(t_{i_2,2},t_{2,1})d(t_{i_3,3},t_{1,1})}}\\
\end{eqnarray*}

Now applying Lemma \ref{lem:presprin} with $G=\pgl{2}{A}$ and $N=\pspl{2}{A}$, we immediately deduce:

\begin{thm}\label{thm:rpaf2} Let $A$ be a local ring with maximal ideal $M$ and residue field $\F{2}$. Then $\rpb{A}$ has the following presentation as a $\sgr{A}$-module:
\begin{enumerate}
\item[Generators:]

$\cconst{A}$, $\suss{i}{u}$, $i\in \{ 1,2\}$, $u\in A^\times$, $\scyc{ 1/m\ n\  u}$, $m,n\in M$, $u\in A^\times$.

\item[Relations:]  Let $W\subset A^\times$ be any subset which additively generates $A$.
\begin{enumerate}
\item For $i\in \{ 1,2\}$, $\suss{i}{uv}=\an{u}\suss{i}{v}+\suss{i}{u}$ for all $u,v\in A^\times$,
\item $\an{-1}\cconst{A}=\cconst{A}$ and $\an{-1}\suss{i}{u}=\suss{i}{u^{-1}}$ for $i\in \{ 1,2\}$ $u\in A^\times$,
\item $\pf{u}\cconst{A}=\suss{1}{u}-\suss{2}{u}$ for all $u\in A^\times$,
\item $3\cconst{A}=0$,
\item For any $w\in W$,
\begin{eqnarray*}
\scyc{ 1/m\  n\ u}&=& \scyc{\frac{1}{1+w}\ \frac{n}{1+nw}\ \frac{1}{w+m}}-\scyc{\frac{1}{1+w}\ \frac{n}{1+nw}\ \frac{1}{w}} \\
&&+\scyc{\frac{u}{1+uw}\ \frac{n}{1+nw}\ \frac{1}{w}}-\scyc{\frac{u}{1+uw}\ \frac{n}{1+nw}\ \frac{1}{w+m}}\\
&&+\suss{1}{\oct}-\suss{2}{(12)\oct}+\suss{1}{(13)\oct}
\end{eqnarray*}
where 
\[
\oct=\left(
\begin{array}{lll}
\frac{1}{1+w}& 0&\frac{1}{w}\\
\frac{u}{1+uw}&\frac{n}{1+nw}&\frac{1}{w+m}\\
\end{array}
\right).
\]
\item For any $w\in W$,
\begin{eqnarray*}
 \scyc{ 1/m\  n\ u}&=&    \scyc{\frac{1+mw}{m}\ u+w\ w}-\scyc{\frac{1+mw}{m}\ u+w\ n+w} \\
&&+\scyc{\frac{1+mw}{m}\ 1+w\ n+w}-\scyc{\frac{1+mw}{m}\ 1+w\ w}\\    
&&+\suss{2}{\oct'}\\
\end{eqnarray*}
where 
\[
\oct'=\left(
\begin{array}{lll}
\infty& 1+w&{w}\\
\frac{1+mw}{m}&u+w&{n+w}\\
\end{array}
\right).
\]
\item For any $w\in A^\times$,
\[
\an{w}\cdot\scyc{ 1/m\  n\ u}=\scyc{w/m\ wn\ wu}-\scyc{ w/m\ wn\ w}.
\]
\end{enumerate}
\end{enumerate}
\end{thm}

For the purpose of applications below, we calculate the elements $\suss{i}{\oct}\in \rpb{A}$ explicitly:

\begin{prop}\label{prop:susscalc} Let $A$ be a local ring with maximal ideal $M$ and residue field $\F{2}$.
\begin{enumerate}
\item For any octahedron $\oct=(t_{i,j})$ and for $i=1,2$, 
\[
\suss{i}{\oct}=\an{d(t_{1,1},t_{1,2},t_{1,3})}\suss{i}{\scyc{ 1/m\ n\ u}}
\]
 where
$1/m:=t_{2,1}\cdot N^{-1}\in 1/M$, $n:=t_{2,2}\cdot N^{-1}\in M$ and $u:=t_{3,1}\cdot N^{-1}\in A^\times$ with $N:=M_{t_{1,1},t_{1,2},t_{1,3}}\in \pgl{2}{A}$.
\item For any $m,n\in M$, $u\in A^\times$, 
\[
\suss{1}{\scyc{ 1/m\ n\ u}} = \alpha_1\suss{1}{1-m}+\alpha_2\suss{1}{1-mn}+\alpha_3\suss{1}{1-mu}
\]
 with
$\alpha_1=\an{1-m}\pf{\frac{1-mn}{1-n}}$, $\alpha_2=\an{u-n}-\an{1-n}$, $\alpha_3=\an{\frac{1-mn}{u-n}}-\an{u}$.
\item For any $m,n\in M$, $u\in A^\times$, 
\[
\suss{2}{\scyc{ 1/m\ n\ u}} = \beta_1\suss{2}{1-m}+\beta_2\suss{2}{1-mn}+\beta_3\suss{2}{1-mu}
\]
 with
$\beta_1=-\pf{1-n}$, $\beta_2=\an{\frac{1-n}{1-m}}-\an{\frac{u-n}{1-mu}}$, $\beta_3=\an{{u-n}}-\an{u}$.
\end{enumerate}
\end{prop}

\begin{proof}
\begin{enumerate}
\item This follows from the fact that 
\[
\oct=(t_{i,j})=
\left(
\begin{array}{lll}
\infty&0&1\\
\frac{1}{m}&n&u
\end{array}
\right)\cdot N
\]
and that $d(a\cdot N, b\cdot N)=d(a,b)\mathrm{det}(N)$ for any unimodular rows $a,b$.
\item We have
\begin{eqnarray*}
\suss{1}{\scyc{\frac{1}{m}\  n\  u}}&=& D_1^{\{ \infty, 0\}}(1,u)-D_1^{\{ \frac{1}{m},0\}}(1,u)-D_1^{\{\infty,n\}}(1,u)+D_1^{\{ \frac{1}{m},n\}}(1,u)\\
&=&\an{d(\infty,0,1)}\suss{1}{\frac{d(\infty,1)d(0,u)}{d(\infty,u)d(0,1)}}-\an{d\left(\frac{1}{m},0,1\right)}\suss{1}{\frac{d(\frac{1}{m},1)d(0,u)}{d(\frac{1}{m},u)d(0,1)}}\\
&&-\an{d(\infty,n,1)}\suss{1}{\frac{d(\infty,1)d(n,u)}{d(\infty,u)d(n,1)}}-\an{d\left(\frac{1}{m},n,1\right)}\suss{1}{\frac{d(\frac{1}{m},1)d(n,u)}{d(\frac{1}{m},u)d(n,1)}}\\
&=&\suss{1}{u}-\an{1-m}\suss{1}{u\cdot\frac{1-m}{1-mu}}-\an{1-n}\suss{1}{\frac{u-n}{1-n}}\\
&&+\an{\frac{1-mn}{(1-m)(1-n)}}\suss{1}{\frac{1-m}{1-mu}\cdot\frac{u-n}{1-n}}\\
&=& \suss{1}{u}-\an{1-m}\left(\an{\frac{1-m}{1-mu}}\suss{1}{u}+\suss{1}{\frac{1-m}{1-mu}}\right)-\an{1-n}\suss{1}{\frac{u-n}{1-n}}\\
&&+\an{\frac{1-mn}{(1-m)(1-n)}}\left(\an{\frac{-1m}{1-mu}}\suss{1}{\frac{u-n}{1-n}}+\suss{1}{\frac{1-m}{1-mu}}\right)\\
&=&-\pf{1-mu}\suss{1}{u}+\an{1-m}\pf{\frac{1-mn}{1-n}}\suss{1}{\frac{1-m}{1-mu}}\\
&&+\an{1-n}\pf{\frac{1-mn}{1-mu}}\suss{1}{\frac{u-n}{1-n}}\\
&=&-\pf{u}\suss{1}{1-mu}+\an{1-m}\pf{\frac{1-mn}{1-n}}\suss{1}{\frac{1-m}{1-mu}}\\
&&+\an{1-n}\pf{\frac{u-n}{1-n}}\suss{1}{\frac{1-mn}{1-mu}}\\
&& \mbox{ (using $\pf{a}\suss{1}{b}=\pf{b}\suss{1}{a}$)}\\
&=&-\pf{u}\suss{1}{1-mu}+\an{1-m}\pf{\frac{1-mn}{1-n}}\left( \an{1-m}\suss{1}{\frac{1}{1-mu}}+\suss{1}{1-m}\right)\\
&&+\an{1-n}\pf{\frac{u-n}{1-n}}\left(\an{1-mn}\suss{1}{\frac{1}{1-mu}}+\suss{1}{1-mn}\right)\\
&=&\left(-\pf{u}+\pf{\frac{1-mn}{1-n}}+\an{\frac{1-mn}{1-n}}\pf{\frac{u-n}{1-n}}\right)\suss{1}{1-mu}\\
&&+\an{1-m}\pf{\frac{1-mn}{1-n}}\suss{1}{1-m}+\an{1-n}\pf{\frac{u-n}{1-n}}\suss{1}{1-mn}\\
&&\mbox{ (using $\pf{a}\suss{1}{b}=\pf{a}\suss{1}{b^{-1}}$)}\\
&=& \left(\an{\frac{1-mn}{u-n}}-\an{u}\right)\suss{1}{1-mu}+\an{1-m}\pf{\frac{1-mn}{1-n}}\suss{1}{1-m}\\
&&+(\an{u-n}-\an{1-n})\suss{1}{1-mn}\\
\end{eqnarray*}
as required.
\item We have
\begin{eqnarray*}
\suss{2}{\scyc{\frac{1}{m}\  n\  u}}&=& D_2^{\{0,1\}}\left(\infty,\frac{1}{m}\right)-D_2^{\{ n,1\}}\left(\infty,\frac{1}{m}\right) -D_2^{\{0,u\}}\left(\infty,\frac{1}{m}\right)
+D_2^{\{n,u\}}\left(\infty,\frac{1}{m}\right) \\
&=&\suss{2}{1-m}-\an{1-n}\suss{2}{\frac{1-m}{1-mn}}-\an{u}\suss{2}{1-mu}+\an{u-n}\suss{2}{\frac{1-mu}{1-mn}}\\
&=&\suss{2}{1-m}-\an{1-n}\left[\an{m-1}\suss{2}{1-mn}+\suss{2}{1-m}\right]\\
&&-\an{u}\suss{2}{1-mu}+\an{u-n}\left[\an{mu-1}\suss{2}{1-mn}+\suss{2}{1-mu}\right]\\
&=&-\an{1-n}\suss{2}{1-m}+\left(\an{\frac{u-n}{mu-1}}-\an{\frac{1-n}{m-1}}\right)\suss{2}{1-mn}\\
&&+\left(\an{u-n}-\an{u}\right)\suss{2}{1-mu}.\\
\end{eqnarray*}

\end{enumerate}
\end{proof}

%%%%%%%%%%%%%%%%%%%%%%%%%%%%%%%%%%%%%%%%%%%%%%%%%%%%%%%%%%%%%%%%%
\section{Examples}\label{sec:exa}
%%%%%%%%%%%%%%%%%%%%%%%%%%%%%%%%%%%%%%%%%%%%%%%%%%%%%%%%%%%%%%%%%%%%%%%%%%%%%%%
\subsection{The ring $\Z/4$}
%%%%%%%%%%%%%%%%%%%%%%%%%%%%%%%%%%%%%%%%%%%%%%%%%%%%%%%%%%%%%%%%%%%%%%.

%As above for $a\in A$  we will write $a:=(a,1)$ in $Y_0$. Furthermore we will denote $(1,a)\in Y_0$ by $1/a$ (even when $a\not\in A^\times$).  We let $\infty:= 1/0$, as above.

%Thus, when $A$ is a local ring $Y_0(A)=A\cup 1/A$ and observe that $A\cap 1/A=A^\times$. When the residue field $k(A)=\Z/2$, we also have $A=1-A^\times\sqcup A^\times$ and thus 
%$Y_0(A)=A\sqcup 1/(1-A^\times)$.

%Thus $Y_0(\Z/4)=\Z/4\sqcup \{ \infty, 1/2\}=\{ 0,1,2,-1,1/2,\infty\}$.  The neighbours of $\infty$ are the elements of $A$ annd the neighbours of $0$ are $1/A$. By the transitive action of $%\pgl{2}{\Z/4}$,
%each of the six vertices has four neighbours. 

Let $Y$ be the simplicial complex $Y(\Z/4)$. By Proposition \ref{prop:localya}, $Y$ is the clique complex of the complete tripartite graph $K:=K_{\{ 0,2\}, \{ 1,-1\}, \{ \infty,1/2\}}$. 
Thus $Y$ is the surface of an octahedron (and hence, as expected, topologically a $2$-sphere) with the given sets as pairs of antipodes. It has a unique basic octahedron represented by the $2\times 3$ matrix
\[
\oct=(t_{i,j}):=
\left(
\begin{array}{lll}
\infty&0&1\\
\frac{1}{2}&2&-1\\
\end{array}
\right).
\]

The group of automorphisms $\Gamma$ of $Y$ is equal to the group of automorphisms of the underlying graph  (acting on the right). It is thus the wreath product $S_3\wr S_2:=
S_3\rtimes S_2^3\subset S_6$. We will write elements of $S_2^3$ as maps $f:\{ 1,2,3\}\to S_2$. The multiplication is thus: $(\sigma_1,f_1)\cdot (\sigma_2,f_2)=(\sigma_1\sigma_2, (f_1\ast \sigma_2^{-1})f_2)$ where $(f\ast\sigma)(j):=f(j\cdot\sigma^{-1})$ for $j\in \{ 1,2,3\}$. The (right)  action is then as follows:
\[
t_{i,j}\cdot (\sigma,f):=t_{i\cdot f(j\cdot\sigma),j\cdot\sigma}.
\]
Thus the subgroups $S_2$ interchange pairs of antipodes and the $S_3$ permutes the three pairs of antipodes. Observe that the automorphism $(\sigma,f)$ preserves orientation if and only if the number of exchanges of antipodes has the same parity as the number of transpositions in $\sigma$.
% Let $\Omega=\{ 0,1,\infty\}=\projl{\F{2}}$ 
%Thus $\Gamma$ is isomorphic to the semidirect product $\mu_2^\Omega\ltimes S(\Omega)$ where  and $S(\Omega)$ denotes the full group of permutations of the set $\Omega$.  The product is given by $(f_1,\sigma_1)\cdot (f_2,\sigma_2)=(f_1\cdot(f_2\ast \sigma_1^{-1}),\sigma_1\sigma_2)$ where $(f\ast\sigma)(i):=f(i\cdot\sigma^{-1})$ for $f\in C^\Omega$, $\sigma\in S(\Omega)$, $i\in \Omega$.

 The (right) action of $\pgl{2}{\Z/4}$ on $K$ induces a homomorphism $\iota: \pgl{2}{\Z/4}\to\Gamma$. Let $\pi:\pgl{2}{\Z/4}\to\pgl{2}{\Z/2}$ be the natural surjections and let $U:=\ker{\pi}$. 

For convenience, we identify the groups $S_2$ and $\mu_2$, so that $\Gamma=S_3\rtimes \mu_2^3$. Then there is an induced commutative diagram of group extensions
\[
\xymatrix{
1\ar[r]&U\ar[r]\ar^-{\iota}[d]&\pgl{2}{\Z/4}\ar^-{\pi}[r]\ar^-{\iota}[d]&\pgl{2}{\Z/2}\ar[r]\ar^-{\cong}[d]&1\\
1\ar[r]&\mu_2^3\ar[r]&\Gamma\ar[r]&S_3\ar[r]&1.\\
}
\]
For $j\in \{ 1,2,3\}$, let $g_j\in \mu_2^3$ be the function
\[
g_j(k):=
\left\{
\begin{array}{ll}
-1,& k=j,\\
1,&k\not=j
\end{array}
\right.
\]
for all $i,j\in \Omega$. One easily verifies that $\iota(\rho_j)=g_j$ for $i\in \Omega$, where 
\[
\rho_1:=
\left[
\begin{array}{cc}
-1&2\\
0&1\\
\end{array}
\right],\quad
\rho_2:=
\left[
\begin{array}{cc}
-1&0\\
2&1\\
\end{array}
\right]
,\quad
\rho_3:=
\left[
\begin{array}{cc}
-1&0\\
0&1\\
\end{array}
\right]
\in U.
\]
It follows that $\iota:\pgl{2}{\Z/4}\to\Gamma$ is an isomorphism. Furthermore, one verifies that the diagram 
\[
\xymatrix{
\pgl{2}{\Z/4}\ar^-{\mathrm{det}}[rd]\ar^-{\iota}_-{\cong}[rr]&&\Gamma\ar_-{\delta}[ld]\\
&\mu_2&\\
}
\]
commutes where $\delta((\sigma,f)):=\mathrm{sgn}(\sigma)\prod_{j=1}^3f(j)\in \mu_2$. 

Thus $\pspl{2}{\Z/4}\cong \Gamma^0:= \{ (\sigma,f)\in \Gamma\ |\ \prod_{j=1}^3f(j)=\mathrm{sgn}(\sigma)\}$.
%$(C^\Omega)^0:= \{ f\in C^\Omega\ |\ \prod f(i)=1\}$. 
It follows that $\pspl{2}{\Z/4}$, and hence also $\spl{2}{\Z/4}$,  acts on the octahedron $Y$ via orientation-preserving automorphisms (and is also isomorphic to $S_4$).

%%%%%%%%%%%%%%%%%%%%%%%%%%%%%%%%%%%%%%%%%%%%%%%%%%%%%%%%%%%%%%%%%%%%%%%
\subsubsection{Calculation of $\rpb{\Z/4}$}  
%%%%%%%%%%%%%%%%%%%%%%%%%%%%%%%%%%%%%%%%%%%%%%%%%%%%%%%%%%%%%%%%%%%%

As noted above, there is a unique basic octahedron with corresponding cycle
 $\scyc{\oct}=\scyc{ 1/2\ 2\ -1}$.

\begin{lem}\label{lem:octz4} For any octahedron $\oct'$ and for $i\in \{ 1,2\}$, $\suss{i}{\oct'}=0$ in $\crpb{\Z/4}$. 
\end{lem}
\begin{proof} The only nontrivial unit is $-1$ and thus $\suss{i}{u}=0$ or $\suss{i}{-1}$ if $u$ is a unit. The augmentation ideal $\aug{\Z/4}$ of $\sgr{\Z/4}$ is generated by $\pf{-1}$.
Since $\pf{-1}\suss{i}{-1}=0$ always and since the coefficients $\alpha_i,\beta_i$ in Proposition \ref{prop:susscalc} (2) and (3) always lie in this augmentation ideal, it follows that 
$\suss{i}{\oct'}=0$ when $\oct$ is the basic octahedron. By Proposition \ref{prop:susscalc} (1), it follows that $\suss{i}{\oct'}=0$ for any octahedron.
\end{proof}

\begin{prop} Let $\oct$ be the basic octahedron for $\Z/4$. As a $\sgr{\Z/4}$-module, $\rpb{\Z/4}$ decomposes as a direct sum $\crpb{\Z/4}\oplus \Z\cdot\scyc{\oct}$ where
$\an{-1}\scyc{\oct}=-\scyc{\oct}$ and $\crpb{\Z/4}$ is cyclic of order $6$ with generator $\cconst{\Z/4}+\suss{1}{-1}$ with trivial action of $\sgr{\Z/4}$.
\end{prop}

\begin{proof} We apply Theorem \ref{thm:rpaf2}. Recalling that, by definition, $\scyc{1/m\ n\ u}= 0$ whenever $m=0$ or $n=0$ or $u=1$, and given Lemma \ref{lem:octz4}, we see that 
relations (e) and (f) both reduce to the trivial relation $\scyc{\oct}=\scyc{\oct}$. Relation (g) becomes $\an{-1}\scyc{\oct}=-\scyc{\oct}$ (for the basic octahedron  $\oct$). 

The remaining generators and relations are just the presentation of $\crpb{A}$ in Proposition \ref{prop:crpbpres}. It follows that $\crpb{\Z/4}$ factors as a direct summand of $\rpb{\Z/4}$.
Proposition \ref{prop:crpbstruct} then implies that it is equal to $(\Z/3\cdot\cconst{\Z/4})\oplus (\Z/2\cdot \suss{1}{-1})$. 
\end{proof}

%%%%%%%%%%%%%%%%%%%%%%%%%%%%%%%%%%%%%%%%%%%%%%%%%%%%%%%%%%%%%%%%%%%%
\subsubsection{Calculation of $\rpbker{\Z/4}$ and $\rbl{\Z/4}$}  
%%%%%%%%%%%%%%%%%%%%%%%%%%%%%%%%%%%%%%%%%%%%%%%%%%%%%%%%%%%%%%%%%%%%
Recall that $\lambda_1:\rpb{A}\to \sgr{A}$ is the map $\rpb{A}:=Z_2(L_\bullet)_{\pspl{2}{A}}\to (L_2)_{\pspl{2}{A}}\cong\sgr{A}$ 
induced by the inclusion $Z_2(L_\bullet)\to L_2$. The isomorphism $ (L_2)_{\pspl{2}{A}}\cong\sgr{A}$ sends the class of $(a,b,c)$ to $D(a,b,c):=\an{\mathrm{det}(M_{a,b,c})}=
\an{d(a,b)d(b,c)d(c,a)}$.

We have 
\[
\rpb{\Z/4}=\Z/2\cdot\suss{1}{-1}\oplus \Z/3\cdot \cconst{\Z/4}\oplus \Z\cdot \Omega_1\cong \Z/2\oplus\Z/3\oplus \Z.
\]
As always, $0=\lambda_1(\suss{1}{-1})=\lambda_1(\cconst{\Z/4})$. On the other hand, if we let $\oct$ denote the basic octahedron for $\Z/4$, then
\begin{eqnarray*}
\lambda_1(\scyc{\oct})&=& D(0,\infty,1)+D(\infty,2,1)+D(2,1/2,1)+D(1/2,0,1)\\
&&-D(0,\infty,-1)-D(\infty,2,-1)-D(2,1/2,-1)-D(1/2,0,-1)\\
&=& \an{-1\cdot 1\cdot 1}+\an{1\cdot 1\cdot -1}+\an{-1\cdot -1\cdot -1}+\an{1\cdot -1\cdot 1}\\
&&-\an{-1\cdot 1\cdot -1}-\an{1\cdot -1\cdot -1}-\an{-1\cdot -1\cdot 1}-\an{1\cdot 1\cdot 1}\\
&=& 4(\an{-1}-\an{1})=4\pf{-1}\in \sgr{\Z/4}.
\end{eqnarray*}
Since this element has infinite order, we deduce that 
\[
\rpbker{\Z/4}=\crpb{\Z/4}=\Z/2\cdot\suss{1}{-1}\oplus \Z/3\cdot \cconst{\Z/4}\cong \Z/6.
\]

Finally, recall that $\rbl{\Z/4}=\ker{d^3}$ where 
 $d^3:E^3_{0,3}=\rpbker{\Z/4}\to E^3_{2,0}$ is the appropriate differential in the spectral sequence $E(\spl{2}{\Z/4},L^\tau(\Z/4))$.

For any ring $A$, $\cconst{A}\in\rbl{A}$. On the other hand, by \cite[Lemma 5.8, Lemma A.12]{hut:bloch1}, $d^3(\suss{1}{-1})$ has order $2$ in $E^3_{2,0}(\spl{2}{\Z/4},L^\tau(\Z/4))$. Thus:
%%%%%%%%%%%%%%%%%%%%%%%%%%%%%%%%%%%%%%%%%%%%%%%%%%%%%%%%%%%%%%%%%%%%%%%%%%%
\begin{cor} \label{cor:rblz4} The group  $\rbl{\Z/4}$ is cyclic of order $3$ generated by $\cconst{\Z/4}$. It does not contain the element $\suss{1}{-1}$.
\end{cor}
%%%%%%%%%%%%%%%%%%%%%%%%%%%%%%%%%%%%%%%%%%%%%%%%%%%%%%%%%%%%%%%%%%%%%%%%
\subsubsection{Low-dimensional homology of $\spl{2}{\Z/4}$}
Let $\mathcal{R}:=\sgr{\Z/4}=\Z[C]$ where $C$ is cyclic of order $2$ with generator $\an{-1}$. The augmentation ideal $\mathcal{I}$ is a free abelian group of rank $1$ with generator $\pf{-1}$.
We consider the spectral sequence (of $\mathcal{R}$-modules) $E(\Z/4, \LL)$ converging to 
$\ho{\bullet}{\spl{2}{\Z/4}}{\Z}$.

 The group $B$ of lower triangular matrices is just the abelian group $\mu_2\times \an{\rho}\cong\mu_2\times (\Z/4)$. Thus, $E^1_{1,0}=\ho{1}{B}{\Z}\cong B$ and, by the K\"unneth formula, $E^1_{2,0}=\ho{2}{B}{\Z}=\mu_2\otimes \an{\rho}\cong \Z/2$ and 
\[
E^1_{3,0}\cong\ho{3}{B}{\Z}\cong \ho{3}{\mu_2}{\Z}\oplus\ho{3}{\an{\rho}}{\Z}\oplus\mathrm{Tor}(\mu_2,\an{\rho})\cong \Z/2\oplus\Z/4\oplus\Z/2.
\]
By the calculations above, the differential $d^1=\lambda_1:E^1_{0,3}=\rpb{\Z/4}\to\mathcal{R}=E^1_{0,2}$ has kernel $\rpbker{\Z/4}\cong\Z/6$ and image $\mathcal{I}/4\mathcal{I}\cong \Z/4$. 

We calculate the term $E^1_{1,3}=\ho{1}{\spl{2}{\Z/4}}{Z_2(L_\bullet)}$:\\
 the short exact sequence of 
$\Z[\spl{2}{\Z/4}]$-modules 
\[
0\to Z_2(\Sigma_\bullet)\to Z_2(L_\bullet)\to H_2(|Y(\Z/4)|\cong \Z\to 0
\]
gives a long exact homology sequence
\[
\ho{1}{\spl{2}{\Z/4}}{Z_2(\Sigma_\bullet)}\to\ho{1}{\spl{2}{\Z/4}}{Z_2(L_\bullet)}\to\ho{1}{\spl{2}{\Z/4}}{\Z}\to\crpb{\Z/4}\to\rpb{\Z/4}
\]
Since the right-hand map is injective and $ \ho{1}{\spl{2}{\Z/4}}{Z_2(\Sigma_\bullet)}=0$ (by a calculation identical to that for $\spl{2}{\Z}$ above), we deduce that 
$E^1_{1,3}\cong\ho{1}{\spl{2}{\Z/4}}{\Z}$. 

Recall that the differentials $d^1:E^1_{p,2}\cong  \mathcal{R}\otimes\ho{p}{\mu_2}{\Z}\to\ho{p}{\mu_2}{\Z}\cong E^1_{p,1}$ are the maps $\epsilon\otimes\mathrm{Id}$ which are surjective with kernel $\mathcal{I}\otimes\ho{p}{\mu_2}{\Z}$.

Thus $E^2_{1,2}$ is the cokernel of the map $d^1:E^1_{1,3}\cong\Z/4\to \mathcal{I}\otimes\ho{1}{\mu_2}{\Z}\cong\Z/2$ and $E^2_{1,3}$ is the kernel of this map. (We will see shortly that this map is the zero map.)

The $E^2$-page has the form:
\[
\xymatrix{
0&0&0&\cdots\\
\rpbker{\Z/4}&E^2_{1,3}&E^2_{2,3}&\cdots\\
\mathcal{I}/4\cong\Z/4\ar_-{d^2}[rdd]&E^2_{1,2}\ar_-{d^2}[rdd]&0&\cdots\\
0&0&0&\cdots\\
\Z&B=\mu_2\times\an{\rho}&\ho{2}{B}{\Z}\cong \Z/2&\ho{3}{B}{\Z}\cong (\Z/2)^{ 2}\oplus\Z/4\\
}
\]
The differential $d^2:E^2_{0,2}=\mathcal{I}/4\to \mu_2\times\an{\rho}=E^2_{1,0}$ sends $\pf{-1}$ to $-I\cdot \rho^{-6}=-I\cdot\rho^2$ (see \cite[Lemma A.6]{hut:bloch1}). 
Thus $E^3_{0,1}=E^\infty_{0,1}\cong \Z/4$ and $E^\infty_{0,2}\cong \Z/2$ (with generator $2\pf{-1}$).

It follows that the natural map $B=\ho{1}{B}{\Z}\to \ho{1}{\spl{2}{\Z/4}}{\Z}$ induces an isomorphism $\Z/4\cong B/\an{-I\cdot\rho^2}\cong\ho{1}{\spl{2}{\Z/4}}{\Z}$. Note that it follows that the induced map $\ho{1}{\mu_2}{\Z}\cong\mu_2\to\gab{\spl{2}{\Z/4}}\cong \ho{1}{\spl{2}{\Z/4}}{\Z}$ is injective. Since $\tilde{\omega}^2=-I\in\mu_2$ it follows that there is a commutative triangle
\[
\xymatrix{
\an{\tilde{\omega}}\ar[r]\ar^-{\cong}[rd]&\spl{2}{\Z/4}\ar@{>>}[d]\\
&\gab{\spl{2}{\Z/4}}\\
}
\]
from which it follows that the horizontal map is a split injection and hence that this map induces injections $\ho{k}{\an{\tilde{\omega}}}{\Z}\to\ho{k}{\spl{2}{\Z/4}}{\Z}$ for all $k\geq 0$.
In particular $\ho{k}{\spl{2}{\Z/4}}{\Z}$ has a cyclic direct factor of order $4$ for all odd $k$.

\begin{lem}\label{lem:d1}  The differential $d^1:E^1_{1,3}\to E^1_{1,2}$ is the zero map.
\end{lem}
\begin{proof} Let $C_\bullet(G)$ denote the standard homogeneous resolution of $\Z$ as a left module over a group $G$, and let $C_\bullet:=C_\bullet(\spl{2}{\Z/4})$.

Since $\tilde{\omega}$ generates $\gab{\spl{2}{\Z/4}}\cong\ho{1}{\spl{2}{\Z/4}}{\Z}$, the class $1\otimes (1,\tilde{\omega})\in \Z\otimes_{\Z[\spl{2}{\Z/4}]}C_1$ represents a generator of $\ho{1}{\spl{2}{\Z/4}}{\Z}$. Now the class $\scyc{\oct}\in Z_2(L_\bullet)$ maps to a generator of the trivial $\spl{2}{\Z/4}$-module $H_2(|Y|)\cong \Z$ under the map 
$Z_2(L_\bullet)\to H_2(|Y|)$.  Thus $\scyc{\oct}\otimes (1,\tilde{\omega})\in Z_2(L_\bullet)\otimes C_1$ maps to a generator of $\ho{1}{\spl{2}{\Z/4}}{H_2(|Y|)}\cong\ho{1}{\spl{2}{\Z/4}}{\Z}$.

Now $E^1_{1,3}$ is naturally a quotient of $\ker{\mathrm{Id}\otimes d_1: Z_2(L_\bullet)\otimes C_1\to Z_2(L_\bullet)\otimes C_0=Z_2(L_\bullet)}$. We have
\[
(\mathrm{id}\otimes d_1)(\scyc{\oct}\otimes (1,\tilde{\omega})=\scyc{\oct}\cdot(\tilde{\omega}-1).
\]
But 
\begin{eqnarray*}
\scyc{\oct}\cdot\tilde{\omega}&=&
\left[
\begin{array}{ccc}
\infty&0&1\\
\frac{1}{2}&2&-1\\
\end{array}
\right]\cdot\tilde{\omega}=\left[
\begin{array}{ccc}
0&\infty&-1\\
{2}&\frac12&1\\
\end{array}
\right]\\
&=&-\left[
\begin{array}{ccc}
0&\infty&1\\
{2}&\frac12&-1\\
\end{array}
\right]=-(12)\cdot\scyc{\oct}\\
&=&\scyc{\oct}-D_1(\oct)\\
\end{eqnarray*}
by Lemma \ref{lem:oct2} (1), and thus $\scyc{\oct}\cdot(\tilde{\omega}-1)=-D_1(\oct)\in Z_2(\sigma_\bullet)$.

Let $\Psi:=D_1^{\{ \infty,0\}}(1,-1)\in Z_2(\Sigma_\bullet)$. Let 
\[
A_1:=
\left[
\begin{array}{cc}
0&1\\
-1&2
\end{array}
\right],\quad A_2:=
\left[
\begin{array}{cc}
2&-1\\
1&0\\
\end{array}
\right],
\quad
A:=
\left[
\begin{array}{cc}
1&2\\
2&1\\
\end{array}
\right]\in\spl{2}{\Z/4}.
\]
Then
\[
D_1(\oct)=\Psi\cdot(I-A_1+A-A_2)=(\mathrm{Id}\otimes d_1)(\Psi\otimes [(A_2,A)-(I,A_1)]\in Z_2(\Sigma_\bullet)\otimes C_1\subset Z_2(L_\bullet)\otimes C_1
\]
and hence 
\[
z:=\scyc{\oct}\otimes (1,\tilde{\omega})+\Psi\otimes [(A_2,A)-(I,A_1)]\in Z_2(L_\bullet)\otimes C_1
\]
represents an element of $E^1_{1,3}$ mapping to a generator of $\ho{1}{\spl{2}{\Z/4}}{\Z}$.

Now, from the definitions, we calculate
\begin{eqnarray*}
\scyc{\oct}&=&(\infty,0,1)\cdot (M_{-1}-I)\cdot X\mbox{ and }\\
\Psi&=&(\infty,0,1)\cdot (M_{-1}-I)(I-\tilde{\omega})\\
\end{eqnarray*}
where 
\[
X:=I+\left[
\begin{array}{cc}
1&0\\
2&1\\
\end{array}
\right]+\left[
\begin{array}{cc}
1&2\\
0&1\\
\end{array}
\right]+\left[
\begin{array}{cc}
1&2\\
2&1\\
\end{array}
\right]\in \Z[\spl{2}{\Z/4}].
\]
We thus find $z=(\infty,0,1)\cdot (M_{-1}-I)\otimes Y$ where 
\begin{eqnarray*}
Y&=&\left[ (I,\tilde{\omega})+\left(\left[
\begin{array}{cc}
1&0\\
2&1\\
\end{array}
\right],\left[
\begin{array}{cc}
0&1\\
-1&2\\
\end{array}
\right]\right)+
\left(\left[
\begin{array}{cc}
1&2\\
0&1\\
\end{array}
\right],\left[
\begin{array}{cc}
2&1\\
-1&0\\
\end{array}
\right]\right)
+\left(\left[
\begin{array}{cc}
1&2\\
2&1\\
\end{array}
\right],\left[
\begin{array}{cc}
2&1\\
-1&2\\
\end{array}
\right]\right)\right]\\
&&+\left[
\left(\left[
\begin{array}{cc}
2&-1\\
1&0\\
\end{array}
\right],\left[
\begin{array}{cc}
1&2\\
2&1\\
\end{array}
\right]\right)
-
\left(I,\left[
\begin{array}{cc}
0&1\\
-1&2\\
\end{array}
\right]\right)
-
\left(\left[
\begin{array}{cc}
1&0\\
2&1\\
\end{array}
\right],\left[
\begin{array}{cc}
2&1\\
-1&2\\
\end{array}
\right]\right)
+
\left(\tilde{\omega},\left[
\begin{array}{cc}
-1&2\\
0&-1\\
\end{array}
\right]\right)
\right]\\
\end{eqnarray*}
in $C_1$.  

To calculate $d^1(z)$ we first consider it as an element of $L_2\otimes C_1$, via the inclusion $Z_2(L_\bullet)\to L_2$.  Its image under the sequence of isomorphisms
\begin{eqnarray*}
L_2\otimes_{\Z[\spl{2}{\Z/4}]}C_1&=&(\infty,0,1)\cdot\Z[\gl{2}{\Z/4}/\mu_2]\otimes_{\Z[\spl{2}{\Z/4}]} C_1\\
&\cong&(\Ind{\mu_2}{\gl{2}{\Z/4}}\Z)\otimes_{\Z[\spl{2}{\Z/4}]}C_1\\
&\cong& \left(\mathcal{R}\otimes_{\Z}\Ind{\mu_2}{\spl{2}{\Z/4}}\Z\right)\otimes_{\Z[\spl{2}{\Z/4}]}C_1\\
&\cong&\mathcal{R}\otimes_{\Z}\left(\Ind{\mu_2}{\spl{2}{\Z/4}}\Z \otimes_{\Z[\spl{2}{\Z/4}]}C_1\right)\\
&\cong &\mathcal{R}\otimes_{\Z[\mu_2]}C_1
\end{eqnarray*}
is $\pf{-1}\otimes Y$ (which lies in $\mathcal{I}\otimes C_1\cong \Z\otimes C_1$, as expected).

Now $1\otimes Y\in \Z\otimes_{\Z[\mu_2]}C_1$ represents an element of $\ho{1}{\mu_2}{\Z}\cong\mu_2$. To determine this element, we construct an augmentation-preserving map of 
left $\Z[\mu_2]$-complexes $\gamma:C_\bullet\to C_\bullet(\mu_2)$ as follows:
Fix a section $S:\mu_2\backslash \spl{2}{\Z/4}=\pspl{2}{\Z/4}\to\spl{2}{\Z/4}$ and let $\rho:\spl{2}{\Z/4}\to\mu_2$ be the map of left $\mu_2$-sets $g\mapsto g\cdot s(\bar{g})^{-1}$.
Then we can define $\gamma_n:C_n\to C_n(\mu_2)$ by $\gamma(g_0,\ldots,g_n):=(\rho(g_0),\ldots,\rho(g_n))$. Then $1\otimes Y$ corresponds to the cycle $1\otimes \gamma_1(Y)\in \Z\otimes_{\Z[\mu_2]} C_1(\mu_2)$. If $\gamma(Y)=\sum_i n_i(\epsilon_{0,i},\epsilon_{1,i})$, then $1\otimes\gamma(Y)$ represents $\prod_i(\epsilon_{0,i}\epsilon_{1,i})^{n_i}\in \ho{1}{\mu_2}{\Z}=\mu_2$.

We choose a section $s$ with the following values: $s(I)=I$, $s(\tilde{\omega})=\tilde{\omega}$ and $s(\bar{g})=g$ for 
\[
g=\left[
\begin{array}{cc}
-1&0\\
2&-1\\
\end{array}
\right],\quad 
\left[
\begin{array}{cc}
1&2\\
0&1\\
\end{array}
\right],\quad
\left[
\begin{array}{cc}
1&2\\
2&1\\
\end{array}
\right],\quad
\left[
\begin{array}{cc}
0&-1\\
1&2\\
\end{array}
\right],\quad
\left[
\begin{array}{cc}
2&1\\
-1&0\\
\end{array}
\right],\mbox{ and }
\left[
\begin{array}{cc}
2&-1\\
1&2\\
\end{array}
\right].
\]
Then $\gamma(Y)$ is equal to
\[
(1,1)+(-1,-1)+(1,1)+(1,-1)+(-1,1)-(1,-1)-(1,-1)+(1,-1)
\]
which maps to the trivial element $1\in \mu_2$.
\end{proof}

% It follows also that the maps  $\ho{k}{\mu_2}{\Z}\to\ho{k}{\spl{2}{\Z/4}}{\Z}$ are injective for all $k$. Hence $E^2_{3,0}=\ho{3}{\mu_2}{\Z}=E^\infty_{3,0}$ and the differential
%$d^3:E^3_{1,3}\to E^3_{3,0}$ is the zero map. 
It follows that $E^2_{1,2}=\mathcal{I}\otimes\mu_2\cong\Z/2$ and $E^2_{1,3}\cong\ho{1}{\spl{2}{\Z/4}}{\Z}\cong\Z/4$.
The differential $d^2:E^2_{1,2}\to E^2_{2,0}=\ho{2}{B}{\Z}$ is the zero map (\cite[Corollary A.9]{hut:bloch1}).
% while

The $E^3$-page thus  looks like:
\[
\xymatrix{
0&0&0&\cdots\\
\rpbker{\Z/4}\ar^-{d^3}[rrddd]&E^3_{1,3}\cong\Z/4\ar^-{d^3}[rrddd]&E^\infty_{2,3}&\cdots\\
2\mathcal{I}/4\mathcal{I}\cong\Z/2&E^3_{1,2}=\mathcal{I}\otimes\mu_2\cong\Z/2&0&\cdots\\
0&0&0&\cdots\\
\Z&\Z/4&\Z/2&\ho{3}{B}{\Z}\cong(\Z/2)^2\oplus\Z/4\\
}
\]

The differential $d^3:E_{0,3}=\rpbker{\Z/4}\to E^3_{2,0}\cong\ho{2}{B}{\Z}\cong\mu_2\otimes\an{\rho}\cong\Z/2$ is surjective (as noted above, in the calculation of $\rbl{\Z/4}$). 
Thus $E^4_{2,0}= E^\infty_{2,0}=0$ and $\rbl{\Z/4}:=E^\infty_{0,3}=\Z/3\cdot\cconst{\Z/4}$. 

It follows that the edge homomorphism $\ho{2}{\spl{2}{\Z/4}}{\Z}\to E^\infty_{0,2}=
2\mathcal{I}/4\mathcal{I}\cong \Z/2$ is an isomorphism.

%Finally, we wil calculate the differential $d^3:E^3_{1,3}\cong \Z/4\to E^3_{3,0}\cong\ho{3}{B}{\Z}$  (by an indirect argument):
\begin{prop}\label{prop:z4}  There is a commutative diagram of abelian groups with split exact rows and columns
\[
\xymatrix{
&0\ar[d]&0\ar[d]&&\\
0\ar[r]&\ho{3}{\an{\tilde{\omega}}}{\Z}\ar[r]\ar[d]&\ho{3}{\spl{2}{\Z}}{\Z}\ar[r]\ar[d]&\rbl{\Z}\ar[r]\ar^-{\cong}[d]&0\\
0\ar[r]&\mathcal{K}\ar[r]\ar[d]&\ho{3}{\spl{2}{\Z/4}}{\Z}\ar[r]\ar[d]&\rbl{\Z/4}\ar[r]&0\\
&C\ar^-{\cong}[r]\ar[d]&C \ar[d]&&\\
&0&0&&\\
}
\]
where $C\cong (\Z/2)^2$, $\Z/4$ or $\Z/2 \ \oplus\  \Z/4$.

In particular, as an abelian group $\ho{3}{\spl{2}{\Z/4}}{\Z}\cong C\oplus \Z/12$.
\end{prop}

\begin{proof}
Let $\mathcal{K}$ denote the kernel of the edge homomorphism $\ho{3}{\spl{2}{\Z/4}}{\Z}\to E^\infty_{0,3}=\rbl{\Z/4}\cong\Z/3$. Thus $\mathcal{K}$ is the 
$2$-Sylow subgroup of $\ho{3}{\spl{2}{\Z/4}}{\Z}$.  Then the calculations above, together with the map of spectral sequences\\
 $E(\spl{2}{\Z},\LL)\to E(\spl{2}{\Z/4},\LL)$,  gives a commutative diagram with exact rows
\[
\xymatrix{0\ar[r]&\ho{3}{\mu_2}{\Z}\ar[r]\ar[d]&\ho{3}{\an{\tilde{\omega}}}{\Z}\cong\Z/4\ar[r]\ar[d]&E^\infty_{1,2}\cong\Z/2\ar[r]\ar^-{\cong}[d]&0\\
\Z/4\cong E^3_{1,3}\ar^-{d^3}[r]&\ho{3}{B}{\Z}\ar[r]&\mathcal{K}\ar[r]&E^\infty_{1,2}\ar[r]&0
}
\]
where all vertical arrows are injective.  

If $d^3$ is the zero map, it then follows that $\mathcal{K}\cong \ho{3}{\an{\tilde{\omega}}}{\Z}\oplus \ho{3}{\an{\rho}}{\Z}\oplus \Z/2$ and the statement of the proposition follows 
with $C= \ho{3}{\an{\rho}}{\Z}\oplus \Z/2\cong \Z/4\ \oplus \ \Z/2$.

Observe that the natural map $\an{\rho}\to\gab{\spl{2}{\Z/4}}$ is an isomorphism, so that the inclusion $\an{\rho}\to\spl{2}{\Z/4}$ is split and hence the map
$\Z/4\cong\ho{3}{\an{\rho}}{\Z}\to\ho{3}{\spl{2}{\Z/4}}{\Z}$ is injective. It follows that the map $d^3$ in our diagram cannot be injective. Thus, if it is non-zero, it's image must have order $2$. 

Likewise, we note that the map $\an{\tilde{\omega}}\to \gab{\spl{2}{\Z/4}}$ is an isomorphism, and hence $\an{\tilde{\omega}}\to\spl{2}{\Z/4}$ is also a split inclusion, so that 
$\ho{3}{\an{\tilde{\omega}}}{\Z}\to\ho{3}{\spl{2}{\Z/4}}{\Z}$ is also a split inclusion. It follows that if $d^3\not=0$ then $\mathcal{K}\cong C\oplus\ho{3}{\an{\tilde{\omega}}}{\Z}$
where $C\cong \Z/4$ or $(\Z/2)^2$.
\end{proof}

\begin{rem}
 Note also that it follows that $E^\infty_{1,3}=\ker{d^3}\not=0$ and, since $E^\infty_{0,4}=0$ there is an induced surjective edge homomorphism $\ho{4}{\spl{2}{\Z/4}}{\Z}\to E^\infty_{1,3}$.
\end{rem}

%%%%%%%%%%%%%%%%%%%%%%%%%%%%%%%%%%%%%%%%%%%%%%%%%%%%%%%%%%%%%%%%%%
\subsection{The ring $\Z/8$}\label{sec:z8}
%%%%%%%%%%%%%%%%%%%%%%%%%%%%%%%%%%%%%%%%%%%%%%%%%%%%%%%%%%%%%%%%%%%%%%%%%%%%%%%

We recall that, by Corollary \ref{cor:2dim}, for any ring local ring $A$ with residue field $\F{2}$ there is an exact sequence of $\pgl{2}{A}$-modules 
\[
0\to Z_2(\Sigma_\bullet)\to Z_2(L_\bullet)\to H_2(|Y(A)|)\to 0
\]
which, on taking $\pspl{2}{A}$-coinvariants gives a long exact sequence terminating with
\[
\xymatrix{
\ho{1}{\pspl{2}{A}}{H_2(|Y(A)|)}\ar^-{\delta}[r]&\crpb{A}\ar[r]&\rpb{A}\ar[r]&\ho{0}{\pspl{2}{A}}{H_2(|Y(A)|)}\ar[r]&0.\\
}
\]
By Proposition \ref{prop:localhom} $Y(\Z/8)$ is homotopy equivalent to a pointed sum of $27=3^3$ $2$-spheres and the cycles $\scyc{\frac{1}{m}\ n\ u}$, $m,n\in \{2,4,6=-2\}$, $u\in \{ 3, 5=-3,7=-1\}$ 
represent a basis of $H_2(|Y(\Z/8)|)\cong \Z^{27}$. For convenience, we enumerate these $27$ basic cycles as follows: For $n=0,1,\ldots,26$ we let
\[
\octt{n}:=\scyc{\frac{1}{2i+2}\ 2j+2\ 2k+3}
\]
where $n=ijk=9\cdot i+3\cdot j+k$ is the base $3$ expansion of $n$, with $i,j,k\in\{ 0,1,2\}$.

By Lemma \ref{lem:action} (e) and (f), we can calculate the action of $\pspl{2}{\Z/8}$ on this basis of $H_2(|Y(\Z/8|)$. The group $\pspl{2}{\Z/8}$ is generated by the elementary matrices $E_{12}(1)$ and $E_{21}(1)$ and we have: For $i,j\in \{ 0,1,2\}$
\begin{eqnarray*}
\octt{9i+3j+2}\cdot E_{12}(1)&=&\octt{3(2-j)+i}\\
\octt{9i+3j+k}\cdot E_{12}(1)&=&\octt{3(2-j)+i}-\octt{9(k+1)+3(2-j)+i}\mbox{ for }k=0,1\\
\octt{9i+3j+2}\cdot E_{21}(1)&=&\octt{9(2-i)+j}\\
\octt{9i+3j+k}\cdot E_{21}(1)&=&\octt{9(2-i)+j}-\octt{9(2-i)+3(k+1)+j}\mbox{ for }k=0,1.\\
\end{eqnarray*}
From these equations we easily calculate that $\ho{0}{\pspl{2}{\Z/8}}{H_2(|Y(\Z/8)|)}$ is a free abelian group of rank $2$ with basis $\{ \octt{1},\octt{2}\}$ and such that 

\begin{eqnarray*}
\octt{3}=\octt{5}=\octt{7}=\octt{9}=\octt{11}=\octt{17}=\octt{19}=\octt{23}&=&\octt{1}\\
\octt{8}=\octt{20}&=&\octt{0}\\
\octt{4}=\octt{10}=\octt{13}=\octt{14}&=&2\octt{1}\\
\octt{2}=\octt{6}=\octt{18}=\octt{26}&=&\octt{0}-\octt{1}\\
\octt{15}=\octt{21}=\octt{25}&=&-\octt{1}\\
 \octt{12}=\octt{16}=\octt{22}&=&0\\
\mbox{ and } \octt{24}&=&\octt{0}-2\octt{1}.\\
\end{eqnarray*}

The $\sgr{\Z/8}$-module structure on  $\ho{0}{\pspl{2}{\Z/8}}{H_2(|Y(\Z/8)|)}$  is given by:
\[
\an{-1}\octt{0}=-\octt{0},\quad \an{-1}\octt{1}=-\octt{1},\quad \an{3}\octt{0}=2\octt{1}-\octt{0},\quad \an{3}\octt{1}=\octt{1}.
\]

Thus we have an exact sequence of $\sgr{\Z/8}$-modules
\[
\xymatrix{
\ho{1}{\pspl{2}{A}}{H_2(|Y(A)|)}\ar^-{\delta}[r]&\crpb{A}\ar[r]&\rpb{A}\ar[r]&\Z\cdot\octt{0}\oplus\Z\cdot\octt{1}\ar[r]&0.\\
}
\]

By Proposition \ref{prop:crpbstruct}, $\crpb{\Z/8}=\ccconstmod{\Z/8}\oplus \cks{i}{\Z/8}$ for $i=1$ or $2$. Since $\mu_2(\Z/8)=(\Z/8)^\times$ it follows that $\tsgr{\Z/8}=\Z/2$ and hence there is an exact sequence 
\[
0\to \Z/2\cdot \suss{i}{-1}\to\cks{i}{\Z/8}\to (1+\an{-1})\aug{\Z/8}\to 0.
\]
Since $\suss{i}{-1}\not=0$ in $\rpb{\Z/8}$ (since the its image in $\rpb{\Z/4}$ does not vanish) it follows that the natural map $\cks{i}{\Z/8}\to\ks{i}{\Z/8}\subset\rpb{\Z/8}$ is injective and thus $\ker{\crpb{\Z/8}\to\rpb{\Z/8}}=\ker{\ccconstmod{\Z/8}\to\cconstmod{\Z/8}}$. 

By Proposition \ref{prop:crpbstruct} again, $\ccconstmod{\Z/8}=\Z/3\cdot\cconst{\Z/8}\oplus \Z/3\cdot\an{3}\cconst{\Z/8}\cong \Z/3\oplus\Z/3$. Since $\cconst{\Z/8}\not=0$ ($\cconst{\Z/4}\not=0$) it follows in turn that $\ker{\crpb{\Z/8}\to\rpb{\Z/8}}$ is contained in the cyclic subgroup of $\ccconstmod{\Z/8}$ of order $3$ generated by the element
$\pf{3}\cconst{\Z/8}=\suss{1}{3}-\suss{2}{3}$.

By Lemma \ref{lem:action}, (e)and (f) -- and using Proposition \ref{prop:susscalc} -- we have
\begin{eqnarray*}
\octt{26}\cdot E_{12}(1)&=& \octt{2}+4\suss{1}{3}-2\suss{2}{3}\\
\octt{26}\cdot E_{21}(1)&=& \octt{2}+2\suss{2}{3}\\
\end{eqnarray*}
in $Z_2(L_\bullet)$. Thus in $\rpb{\Z/8}=Z_2(L_\bullet)_{\pspl{2}{\Z/8}}$ we have 
\begin{eqnarray*}
0=\octt{26}\cdot\left[ (E_{12}(1)-I)-(E_{21}(1)-I)\right] =4(\suss{1}{3}-\suss{2}{3})=4\pf{3}\cconst{\Z/8}=\pf{3}\cconst{\Z/8}
\end{eqnarray*}
and hence $\ker{\crpb{\Z/8}\to\rpb{\Z/8}}$ is cyclic of order $3$ with generator $\pf{3}\cconst{\Z/8}$. Observe that it follows that the homology class in $\ho{1}{\pspl{2}{\Z/8}}{H_2(|Y(\Z/8)|)}$ represented by the cycle 
\[
\octt{26}\otimes \left((I,E_{12}(1))-(I,E_{21}(1)\right)\in H_2(|Y(\Z/8)|)\otimes_{\Z[\pspl{2}{\Z/8}]}C_1(\pspl{2}{\Z/8})
\]
maps, under the connecting homomorphism $\delta$,  to the nonzero element $\pf{3}\cconst{\Z/8}\in \crpb{\Z/8}$.

We have $\lambda_1(\suss{i}{3})=\lambda_1(\suss{i}{-3})=(1+\an{-1})\pf{3}\in\aug{\Z/8}$ for $i=1$ or $2$. Thus $\lambda_1(\ks{i}{\Z/8})\subset\aug{\Z/8}$ is free of rank $1$ with generator $(1+\an{-1})\pf{3}=\pf{3}+\pf{-3}-\pf{-1})$. 
Furthermore $H_2(|Y(\Z/8)|)_{\pspl{2}{\Z/8}}\cong\Z$ is generated by the cycles $\octt{0}=\scyc{\frac{1}{2}\ 2\ 3}$ and 
$\octt{1}=\scyc{\frac{1}{2}\ 2\ 5}$. We have  $\lambda_1(\octt{0})=3\pf{-1}+\pf{3}+\pf{-3}$ and $\lambda_1(\octt{1})=2(\pf{-1}-\pf{3})\in\aug{\Z/8}$. It follows that $\rpbker{\Z/8}=\ker{\lambda_1}=\Z/2\cdot\suss{1}{-1}\oplus \Z/3\cdot\cconst{\Z/8}\cong \Z/6$. Thus the natural map $\rpbker{\Z/8}\to\rpbker{\Z/4}$ is an isomorphism. It  follows that $\rbl{\Z/8}=\Z/3\cdot\cconst{\Z/8}$ and thus $\rbl{\Z/8}\cong\rbl{\Z/4}\cong\rbl{\/2}\rbl{\Z}\cong \Z/3$.
%%%%%%%%%%%%%%%%%%%%%%%%%%%%%%%%%%%%%%%%%%%%%%%%%%%%%%%%%%%%%%%%%%%%%%
\subsection{The ring $\Z_{\an{2}}$}
%%%%%%%%%%%%%%%%%%%%%%%%%%%%%%%
  We consider the localization of $\Z$ at the prime ideal $\an{2}=2\Z$. The resulting ring $\Z_{\an{2}}\subset \Q$ is a local PID with maximal ideal 
$2\Z_{\an{2}}$ and residue field $\F{2}=\Z/2$. Thus $Y:=Y(\Z_{\an{2}})$ is a $2$-dimensional simplicial complex which is homotopy equivalent to a pointed sum of $2$-spheres and there is a natural exact sequence
\[
\crpb{\Z_{\an{2}}}\to\rpb{\Z_{\an{2}}}\to H_2(|Y|)_{\pspl{2}{\Z_{\an{2}}}}\to 0.
\]

We will show that the term $ H_2(|Y|)_{\pspl{2}{\Z_{\an{2}}}}$ contributes non-trivially to $\rbl{\Z_{\an{2}}}$ (unlike the cases $A=\Z/4$ and $A=\Z/8$). More precisely:   Recall that for a ring $A$, $\rbl{A}=\ker{d^3:\rpbker{A}\to E^3_{2,0}}$ where the map $d^3$ is the appropriate differential in the spectral sequence $E(\spl{2}{A}, L^\tau)$.  We also let $\crbl{A}$ denote the kernel of the associated composite homomorphism $\crpbker{A}\to \rpbker{A}\to E^3_{2,2}$.
%For a ring $A$, let \[
%\crbl{A}:=\rbl{A}\cap\image{\crpb{A}\to\rpb{A}}.
%\]
  In this section we will exhibit a non-zero element of the quotient group $\rbl{\Z_{\an{2}}}/\crbl{\Z_{\an{2}}}$. Note that this quotient  is naturally contained in $\rpb{\Z_{\an{2}}}/\crpb{\Z_{\an{2}}}\cong H_2(|Y|)_{\pspl{2}{\Z_{\an{2}}}}$.

We begin by finding a nonzero element of $H_2(|Y|)_{\pspl{2}{\Z_{\an{2}}}}$. Let $M=2\Z_{\an{2}}$ and let 
\[
x=(x_0,x_1,x_2)=(0,1,1/2), s=(s_0,s_1,s_2)=(2,-1,-1/2)\in M\times (1+M)\times 1/M.
\]
Then,  by Proposition \ref{prop:localhom},
\begin{eqnarray*}
\gamma&=&\sum_{y\in\prod_{i=0}^2\{ s_i,x_i\} } (-1)^{\delta(y,s)}y\\
&=&-(0,1,1/2)+(2,1,1/2)+(0,-1,1/2)+(0,1,-1/2)\\
&&-(2,-1,1/2)-(0,-1,-1/2)-(2,1,-1/2)+(2,-1,-1/2)\in Z_2(L_\bullet)\subset C_2(Y)
\end{eqnarray*}
represents an element of infinite order in $H_2(|Y|)$. It also represents an element in $\rpb{\Z_{\an{2}}}=(Z_2(L_\bullet))_{\pspl{2}{\Z_{\an{2}}}}$. We now calculate the image of $\gamma$  under the map $\rpb{\Z_{\an{2}}}\to\rpb{\Q}$ induced by the inclusion $\Z_{\an{2}}\to \Q$. Let
\begin{eqnarray*}
\tilde{\gamma}:=(\gamma,\infty)&=&-(0,1,1/2,\infty)+(2,1,1/2,\infty)+(0,-1,1/2,\infty)+(0,1,-1/2,\infty)\\
&&-(2,-1,1/2,\infty)-(0,-1,-1/2,\infty)-(2,1,-1/2,\infty)+(2,-1,-1/2,\infty)\in L_3(\Q).\\
\end{eqnarray*}
Then $d(\tilde{\gamma})=\gamma\in Z_2(\Q)\subset L_2(\Q)$. Since 
\[
\rpb{\Q}=Z_2(\Q)_{\pspl{2}{\Q}}\cong\coker{L_4(\Q)_{\pspl{2}{\Q}}\to L_3(\Q)_{\pspl{2}{\Q}}}
\]
is a quotient of $L_3(\Q)$, $\tilde{\gamma}$ represents the image of $\gamma$ in $\rpb{\Q}$. Now recall that,  by the results of section \ref{sec:yabasic} above, for any ring $A$, 
if $(a,b,c,d)\in L_3(A)$ then 
\[
(a,b,c,d)=[u(a,b,c,d]\cdot M_{a,b,c}\mapsto \an{\det{M_{a,b,c}}}[u(a,b,c,d)]:=D(a,b,c,d)\in L_3(A)_{\pspl{2}{A}}
\] 
where 
\[
u(a,b,c,d)=d\cdot M_{a,b,c}^{-1}=\frac{d(a,c)d(b,d)}{d(a,d)d(b,c)}\in \wn{A}.
\]
In particular, if $a,b,c\in \Q^\times$ are distinct then
 \[
D(a,b,c,\infty)=\an{(a-b)(b-c)(c-a)}\gpb{\frac{a-c}{b-c}}\in \rpb{\Q}.
\]
Thus $\gamma\in \rpb{\Z_{\an{2}}}$ maps to 
\begin{eqnarray*}
D(\tilde{\gamma})&=&-\an{-1}\gpb{-1}+\an{-3}\gpb{3}+\an{-3}\gpb{1/3}+\an{3}\gpb{1/3}\\
&&-\an{3}\gpb{-1}-\gpb{-1}-\an{-15}\gpb{5/3}+\an{15}\gpb{-5}\in \rpb{\Q}.\\
\end{eqnarray*}

We recall (see, for example, \cite[Section 3.2]{hut:sl2Q}) that for any prime $p$ there is a natural $\sgr{\Q}$-module homomorphism $S_p:\rpb{\Q}\to\qpb{\F{p}}$. Here, for a ring $A$, $\qpb{A}$ denotes the quotient
of $\pb{A}$ by the subgroup generated by the elements $\sus{u}$, $u\in A^\times$. Furthermore, $\qpb{\F{p}}$ has the $\sgr{\Q}$-module structure $\an{a}\cdot x:= (-1)^{v_p(a)}x$ for 
$a\in \Q^\times$, $x\in \qpb{A}$. The homomorphism $S_p$ is determined   by 
\[
S_p(\gpb{a})=
\left\{
\begin{array}{ll}
\gpb{\bar{a}}, & v_p(a)=0,\\
\bconst{\F{p}},& v_p(a)>0,\\
-\bconst{\F{p}},& v_p(a)<0.\\
\end{array}
\right.
\]
Note further that $\bconst{\F{p}}=0$ in $\pb{\F{p}}$ when $p\equiv 1\pmod{12}$ by \cite[Lemma 7.11]{hut:cplx13}. It follows also that $\crpb{\Q}\subset \ker{S_p}$ when $p\equiv 1\pmod{12}$ since $\crpb{\Q}$ is generated by the elements $\cconst{\Q}=4\bconst{\Q}$ and $\suss{i}{u}$ which map to $4\cconst{\F{p}}$ and $\sus{\bar{u}}$. It follows that there is an induced homomorphism $S_p:\rpb{\Z_{\an{2}}}/\crpb{\Z_{\an{2}}}\to \qpb{\F{p}}$ when $p\equiv 1\pmod{12}$.

We consider the case $p=13$:  We begin by reviewing some of the relevant facts about $\pb{\F{13}}$. Recall that for any field $F$ 
\begin{enumerate}
\item There is a homomorphism
$\psi:F^\times/(F^\times)^2\to \pb{F}$, $u\mapsto \sus{u}:=\gpb{u}+\gpb{u^{-1}}$ 
\item There is a homomorphism 
\begin{eqnarray*}
\lambda:\pb{F}&\to &\asym{2}{\Z}{F^\times}:=\frac{F^\times\otimes_{\Z}F^\times}{\an{ x\otimes y+y\otimes x\ |\ x,y\in F^\times}},\\
 \gpb{a}&\mapsto & a\circ (1-a)\\
\end{eqnarray*}
whose kernel is $\bl{F}$.
\end{enumerate}

\begin{prop}\label{prop:f13}  $\pb{\F{13}}=\bl{\F{13}}\times \psi(\F{13}^\times)$ is cyclic of order $14$. $\bl{\F{13}}$ is cyclic of order $7$ with generator $\gpb{4}$. The second factor 
is cyclic of order $2$ with generator $\gpb{7}=\sus{2}$.
\end{prop}
\begin{proof} $\pb{\F{13}}$ is cyclic of order $13+1=14$ by \cite[Corollary 7.6]{hut:cplx13}. As noted above, since $13\equiv 1\pmod{13}$ we have
\begin{eqnarray}\label{cf}
0=\bconst{\F{13}}=\gpb{x}+\gpb{1-x} \mbox{ in } \pb{\F{13}} \mbox{ for all } 1\not=x\in \F{13}^\times.
\end{eqnarray}
Since $7=1-7$ in $\F{13}$ it follows that $2\gpb{7}=0$. Since $7$ is a nonsquare in $\F{13}$, $\lambda(\gpb{7})=7\circ 7$ is a generator of the cyclic group (of order $2$) $\asym{2}{\Z}{\F{13}^\times}$. Thus $\gpb{7}$ is the unique element of $\pb{\F{13}}$ of order $2$. 

Since $-1\in (\F{13}^\times)^2$ we have $0=\sus{-1}=2\gpb{-1}$. However, $\lambda(\gpb{-1})=0$ and hence $\gpb{-1}\not=\gpb{7}$. It follows that $\gpb{-1}=0$.
 Furthermore,  by (\ref{cf}), $0=\gpb{-1}+\gpb{2}$ and hence $\gpb{2}=0$ also. It follows  that $\sus{2}=\gpb{2}+\gpb{7}=\gpb{7}$. Thus, in $\pb{\F{13}}$ we have
\begin{eqnarray}\label{sus}
\sus{u}=\gpb{u}+\gpb{u^{-1}}=
\left\{
\begin{array}{ll}
0,& u\in (\F{13}^\times)^2\\
\gpb{7},&  u\not\in (\F{13}^\times)^2.\\
\end{array}
\right.
\end{eqnarray}
Let 
\[
t=
\left[
\begin{array}{cc}
4&-2\\
-1&4\\
\end{array}
\right]
\in \pspl{2}{\F{13}},
\]
of order $7$. By \cite[p. 64]{hut:cplx13} a generator of $\bl{\F{13}}$ is 
\begin{eqnarray*}
\sum_{i=1}^4\gpb{\frac{\infty\cdot t-\infty\cdot t^{i+1}}{\infty\cdot t-\infty\cdot t^{i+2}}}&=&\gpb{\frac{-4-1}{-4-5}}+\gpb{\frac{-4-5}{-4+5}}+\gpb{\frac{-4+5}{-4+1}}+
\gpb{\frac{-4+1}{-4-4}}\\
&=&\gpb{2}+\gpb{4}+\gpb{4}+\gpb{2}=2\gpb{4}.\\
\end{eqnarray*}
Since $\gpb{4}\in\bl{\F{13}}=\ker{\lambda}$ it follows that $\gpb{4}$ has order $7$ and is a generator. 
\end{proof}

\begin{cor}  The composite of natural homomorphisms $\bl{\F{13}}\to\pb{\F{13}}\to\qpb{\F{13}}$ is an isomorphism.
\end{cor}

\begin{cor}  Let $\gamma\in \rpb{\Z_{\an{2}}}$ as above. Then $S_{13}(\gamma)=3\gpb{4}$ in $\qpb{\F{13}}$ has order $7$. In particular, $\gamma$ is a nonzero element of
$\rpb{\Z_{\an{2}}}/\crpb{\Z_{\an{2}}}\cong H_2(|Y|)_{\pspl{2}{\Z_{\an{2}}}}$.
\end{cor}

\begin{proof} We have 
\begin{eqnarray*}
S_{13}(\gamma)=S_{13}(D(\tilde{\gamma}))&=&-\gpb{-1}+\gpb{3}+\gpb{9}+\gpb{9}-\gpb{-1}-\gpb{-1}-\gpb{6}+\gpb{8}\\
&=& \gpb{9}-\gpb{6}+\gpb{8}\\
\end{eqnarray*}
since $\gpb{-1}=0$ and $\gpb{3}+\gpb{9}=0$ by (\ref{sus}). Now $\gpb{8}+\gpb{6}=0$  and $0=\gpb{9}+\gpb{5}$ by (\ref{cf}) and $0=\gpb{5}+\gpb{8}$ by (\ref{sus}) (since $\gpb{7}=\sus{2}=0$ in $\qpb{\F{13}}$). Thus $S_{13}(\gamma)=3\gpb{8}$. However, the $5$-term relation
\[
0=\gpb{4}-\gpb{8}+\gpb{\frac{8}{4}}-\gpb{\frac{1-\frac{1}{4}}{1-\frac{1}{8}}}+\gpb{\frac{1-4}{1-8}}=\gpb{4}-\gpb{8}+\gpb{4}-\gpb{-1}+\gpb{6}
\]
gives $0=2(\gpb{4}-\gpb{8})$  and hence $\gpb{4}=\gpb{8}$ in $\qpb{\F{13}}$.
\end{proof}

Now recall that $\rpbker{A}=\ker{\lambda_1:\rpb{A}\to \aug{A}\subset\sgr{A}}$ where $\lambda_1(\gpb{u})=\pf{u}\pf{1-u}$. Note that $\lambda_1(\suss{1}{u})=-(1+\an{-1})\pf{u}$ for all 
$u\in A^\times$. Suppose now that $\lambda_1(\gamma)=\sum_im_i\pf{u_i}\in\aug{A}$. Then let $\rho=\rho_\gamma:=\sum_im_i\suss{1}{u_i}\in \crpb{\Z_{\an{2}}}$ and let\\ 
$\gamma_1:=(1+\an{-1})\gamma+\rho$. Observe that $\lambda_1(\gamma_1)=0$ so that $\gamma_1\in \rpbker{\Z_{\an{2}}}$. 

Furthermore, $S_{13}(\gamma_1)=2S_{13}(\gamma)=6\gpb{4}=-\gpb{4}$ in $\qpb{\F{13}}$.

Finally, recall again that $\rbl{\Z_{\an{2}}}:=\ker{d^3: \rpbker{\Z_{\an{2}}}\to E^3_{2,0}}$. Here, the term $E^3_{2,0}$ is a quotient of $E^1_{2,0}=H_2(B,\Z)$ where 
$B$ is the stabilizer of $\infty$ in $\spl{2}{\Z_{\an{2}}}$; i.e., $B$ is the group of upper triangular matrices. We let $T\cong \Z_{\an{2}}^\times$ be the subgroup of 
triangular matrices. We have a natural short exact sequence of groups $0\to \Z_{\an{2}}\to B\to T\to 1$ where the inclusion $T\to B$ gives a splitting. Note that $T$ acts -- by conjugation -- on $\Z_{\an{2}}$ via $u\cdot a=u^2a$.

\begin{lem} \label{lem:tb}
This group extension  induces a decomposition
 \[
H_2(B,\Z)\cong H_2(T,\Z)\oplus H_1(T,\Z_{\an{2}})
\]
 where the second factor is annihilated by $8$.
\end{lem}

\begin{proof} The Hochscilld-Serre spectral sequence associated to this extension has the form $E^2_{p,q}=H_p(T,H_2(\Z_{\an{2}},\Z))\Longrightarrow H_{p+q}(B,\Z)$.
 Observe that edge homomorphisms 
$H_p(B,\Z)\to H_p(T,\Z)=E^2_{p,0}$ are split surjective  for all $p$ since the map $B\to T$ is split by the inclusion map $T\to B$. In particular, $E^2_{p,0}=E^\infty_{p,0}$ and the higher differentials from $E^r_{p,0}$ are all zero for $r\geq 2$. 

It follows at once that $H_1(T,\Z_{\an{2}})=E^2_{1,1}=E^\infty_{1,1}$. Let $g:=\mathrm{diag}(3,1/3)\in T$. On the one hand, the map of pairs $(T,\Z_{\an{2}})\to(T,\Z_{\an{2}})$, 
$(t,m)\mapsto (gtg^{-1},g\cdot  m)=(t,3^2m)$ induces the identity map on the homology groups $H_i(T,\Z_{\an{2}})$. On the other hand, this map is clearly multiplication by $3^2=9$. Thus
$9=1$ in $\mathrm{End}(H_i(T,\Z_{\an{2}}))$ and hence $8=0$. In particular, $8$ annihilates $H_1(T,\Z_{\an{2}})=E^\infty_{1,1}$.

Finally note that $H_2(\Z_{\an{2}},\Z)=\Z_{\an{2}}\wedge \Z_{\an{2}}=0$ and hence\\
 $E^\infty_{0,2}=E^2_{0,2}=H_0(T,H_2(\Z_{\an{2}},\Z))=0$.
\end{proof}

\begin{cor}\label{cor:8x}  Consider the spectral sequence $E(\spl{2}{\Z_{\an{2}}}, L^\tau)$. If $x\in E^3_{2,0}$ then $\pf{u}8x=0$ for all $u\in \Z_{\an{2}}^\times$. 
\end{cor}
\begin{proof}
Since the action of $\sgr{A}$ on $H_i(\spl{2}{A},\Z)$ is induced from conjugation by matrices $\mathrm{diag}(u,1)\in \gl{2}{A}$ and since such conjugation acts trivially on the group $T$ of diagonal matrices, we see that $\sgr{\Z_{\an{2}}}$ acts trivially on $H_2(T,\Z)$. Now $E^3_{2,0}$ is a quotient of $E^1_{2,0}=H_2(B,\Z)$. By Lemma \ref{lem:tb} it follows that  if $x\in E^3_{2,0}$, then $\an{u}\cdot 8x=8x$ for all $u\in \Z_{\an{2}}^\times$; i.e., $\pf{u}8x=0$ for all $u$.
\end{proof}

\begin{prop}\label{prop:gamma2}
 Let $\gamma_2:=\pf{13}8\gamma_1$. Then $\gamma_2$ is a nonzero element of ${\rbl{\Z_{\an{2}}}}/{\crbl{\Z_{\an{2}}}}$.
\end{prop}
\begin{proof} Certainly $\gamma_2\in \rpbker{\Z_{\an{2}}}$. Now $d^3(\gamma_2)=\pf{13}8d_3(\gamma_1)=0$ by Corolllary \ref{cor:8x}, and thus $\gamma_2\in \rbl{\Z_{\an{2}}}$.

 On the other hand, 
since $\an{13}$ acts as multiplication by $-1$ on $\qpb{\F{13}}$, $\pf{13}$ is multiplication by $-2$ and we have $S_{13}(\gamma_2)=-16S_{13}(\gamma_1)=16\gpb{4}=2\gpb{4}$ is a generator (of order $7$) of $\qpb{\F{13}}$.
Since $S_{13}$ annihilates $\crpb{\Z_{\an{2}}}$, the result follows.
\end{proof}

%%%%%%%%%%%%%%%%%%%%%%%%%%%%%%%%%%%%%%%%%%%%%%%%%%%%%%%%%%%%%%%%%%%%%%%%%%%%%%%%
%%%%%%%%%%%%%%%%%%%%%%%%%%%%%%%%%%%%%%%%%%%%%%%%%%%%%%%%%%%%%%%%%%%%%%%%%%%%%%%

%%%%%%%%%%%%%%%%%%%%%%%%%%%%%%%%%%%%%%%%%%%%%%%%%%%%%%%%%%%%%%%%%%%%%
\appendix
%%%%%%%%%%%%%%%%%%%%%%%%%%%%%%%%%%%%%%%%%%%%%%%%%%%%%%%%%%%%%%%%%%%%%%%%%%%%

%%%%%%%%%%%%%%%%%%%%%%%%%%%%%%%%%%%%%%%%%%%%%%%%%%%%%%%%%%%%%%%%%%%%%%%
\section{The clique complex of a complete multipartite graph}\label{app:clique}
%%%%%%%%%%%%%%%%%%%%%%%%%%%%%%%%%%%%%%%%%%%%%%%%%%%%%%%%%%%%%%%%%%%%%%%%%%%%%%%%%%%%%

\subsection{Pointed sums}   Let $\{ (X(i),p_i)\}_{i\in I}$ be a pairwise disjoint family of pointed simplicial complexes ($p_i\in X(i)_0$ for all $i$). Then the simplicial complex $\bigvee_i X(i)$  is the (disjoint) union of the $X(i)$ with the $0$-simplices $p_i$ identified. We will regard $\bigvee_i X(i)$ as a pointed complex with point $p=p_i$ for any $i$.
(If the $X(i)$ are not disjoint, when constructing $\bigvee_i X(i)$, we will tacitly replace them by isomorphic disjoint complexes $X(i)\times \{ i\}$, $i\in I$).

The natural maps $X(i)\to \bigvee_i X(i)$ induce maps of geometric realizations $|X(i)|\to |\bigvee_i X(i)|$ and hence a map of pointed topological spaces $\bigvee_i|X(i)|\to |\bigvee_i X(i)|$ 
which is a homeomorphism (i.e., geometric realization commutes with sums).

\subsection{Multicones}  Let $Y$ be a simplicial complex and let $T$ be a set disjoint from $X$. We let $\mc{T}{Y}$ denote the simplicial complex given by
\[
\mc{T}{Y}_n:= Y_n\cup \{ \sigma\cup \{ t\}\  |\ \sigma\in Y_{n-1}, t\in T\}
\]
(where $Y_{-1}:=\emptyset$).

%\begin{rem}\label{rem:mc}  We observe for future use  that if $T$ $S$ are disjoint from the complex $Y$, and if $S\cap T=\emptyset$, then $\mc{S}{\mc{T}{Y}}=\mc{T}{\mc{S}{Y}}$ since %they have the same sets of $r$-simplices for all $r$.
%\end{rem}

If $T=\{ t\}$ is a singleton we will simply write $\mc{t}{Y}$.  Then $|\mc{t}{Y}|$ is just the cone on $Y$ with vertex $t$. We spell this out in detail:

\begin{lem}\label{lem:cone}  For any simplicial complex $Y$ there is a homeomorphism of spaces 
\[
\frac{|Y|\times [0,1]}{|Y|\times \{ 1\}}\cong |\mc{t}{Y}|.
\]
\end{lem}
\begin{proof}
The inclusion $Y\to \mc{t}{Y}$ induces an embedding of spaces $|Y|\to |\mc{t}{Y}|$. For each $x\in |Y|$, there is a unique line segment $[x,t]$ in $|\mc{t}{Y}|$ and 
$|\mc{t}{Y}|=\cup_{x\in |Y|}[x,t]$. We have $[x,t]\cap [x',t]=\{ t\}$ unless $x=x'$. We can thus assign a unit length to each of these line segments and label the points of $|\mc{t}{Y}|\setminus\{ t\}$ by the points $(x,\lambda)\in |Y|\times [0,1)$, with $(x,0)=x$. We also of course set $(x,1):=t$ for all $x\in |Y|$.  
\end{proof}

\begin{cor}\label{lem:contract}   The space  $|\mc{t}{Y}|$ deforms to the point $t$ by a contraction fixing $t$.
\end{cor}

\begin{proof}
The required contraction is
\begin{eqnarray*}
\mc{t}{Y}\times [0,1]&\to& \mc{t}{Y}\\
((x,\lambda),s)&\mapsto &\left\{
\begin{array}{ll}
(x,\lambda+(1-\lambda)s),& \lambda <1\\
t=(x,1),&\lambda=1.\\
\end{array}
\right.
\end{eqnarray*}

\end{proof}

Note that for any set $T$ disjoint from $Y$, we have $\mc{T}{Y}:=\bigcup_{t\in T} \mc{t}{Y}$. Observe that for disjoint sets $S$ and $T$, $\mc{S}{\mc{T}{Y}}\cong\mc{T}{\mc{S}{Y}}$
canonically.

When $t_1\not=t_2\in T$, then clearly 
 $\mc{t_1}{Y}\cap \mc{t_2}{Y}=Y\subset \mc{T}{Y}$.
Thus to any point $x$ in  $|\mc{T}{Y}|$, if $x\in |\mc{t}{Y}|$ for some $t\in T$,  as in the proof of Lemma \ref{lem:cone}, we may assign a label $(y,\lambda)_t$ to $x$ where $\lambda\in [0,1]$, $y\in |Y|$. We have $(y,1)_t=t$  and
 $(y,0)_t=y$ for all 
$y\in |Y|$, $t\in T$.

\begin{lem}\label{lem:susp} For any simplicial complex $Y$ and any two-element set $T=\{ t_1,t_2\}$, there is a homeomorphism of spaces
\[
\left| \mc{T}{Y}\right| \approx \Sigma |Y|
\]
where $\Sigma (X)$ is the (free) suspension of the space $X$.
\end{lem}

\begin{proof}
Recall that $\Sigma|Y|=|Y|\times [0,1]/(|Y|\times \{ 0\}\cup |Y|\times \{ 1\})$. 

%Now $\mc{\{ t_1,t_2\}}{Y}=\mc{t_1}{Y}\cup\mc{t_2}{Y}$ where $\mc{t_1}{Y}\cap\mc{t_2}{Y}=Y$. The corresponding description then holds for geometric realizations:
%$|\mc{\{ t_1,t_2\}}{Y}|=|\mc{t_1}{Y}|\cup|\mc{t_2}{Y}|$ where $|\mc{t_1}{Y}|\cap|\mc{t_2}{Y}|=|Y|$.
 %For $i=1,2$, we label the points of $\mc{t_i}{Y}|$ with coordinates
% $(x,\lambda)_i$, $x\in |Y|$, $\lambda\in [0,1]$ as in the proof of \ref{lem:contract}. 
For $i=1,2$, let $(x,\lambda)_i:=(x,\lambda)_{t_i}$ in $\mc{T}{Y}|$. Thus that $(x,0)_1=(x,0)_2=x$, $(x,1)_1=t_1$, $(x,1)_2=t_2$ for all $x\in |Y|$. Then the required homeomorphism 
is given by
\begin{eqnarray*}
(x,\lambda)_1&\mapsto & \left( x, (1+\lambda)/2)\right)\\
(x,\lambda)_2&\mapsto & \left(x,(1-\lambda)/2\right).\\
\end{eqnarray*}
\end{proof}

 Let $\tilde{C}_\bullet(Y)$ denote the augmented chain complex
\[
\xymatrix{
\ar[r]&C_1(Y)\ar[r]&C_0(Y)\ar^-{\epsilon}[r]& \Z\ar[r]&0\\
}
\]
which computes the reduced homology groups $\tilde{H}_\bullet(|Y|)$. 

\begin{cor}\label{cor:susp} Let $Y$ be a simplicial complex and let $\gamma\in \tilde{C}_n(Y)$ be a cycle. Let $T=\{ t_1,t_2\}$. Then $(t_1,\gamma)-(t_2,\gamma)\in \tilde{C}_{n+1}(\mc{\{ t_1,t_2\}}{Y})$ 
is also a cycle. (Here, for a chain $w=\sum n_i(x^i_0,\ldots,x^i_n)\in \tilde{C}_n(Y)$, $(t,w)$ denotes the chain $\sum n_i(t,x^i_0,\ldots,x^i_n)\in C_{n+1}(\mc{t}{Y})$. When $w=1\in \tilde{C}_{-1}(Y)=\Z$, $(t,w)$ denotes $t\in \tilde{C}_0(\mc{t}{Y}$.)

Furthermore, for $n\geq 0$, the map on reduced homology $\tilde{H}_n(|Y|)\cong \tilde{H}_{n+1}(|\mc{\{ t_1,t_2\}}{Y}|)$  which sends the homology class $[\gamma]$ to the class
$[(t_1,\gamma)-(t_2,\gamma)]$ is an isomorphism.
\end{cor}

\begin{proof}
Given $w\in \tilde{C}_n(Y)$, the element $(t,w)\in \tilde{C}_{n+1}(\mc{t}{Y})$ satisfies $d(t,w)=w-(t,dw)$.  It folloows that the map $f:\tilde{C}_\bullet(Y)\to \tilde{C}_{n+1}(\mc{\{ t_1,t_2\}}{Y}), w\mapsto (t_1,w)-(t_2,w)$ is a maps of complexes of degree 1. This establishes the first assertion.

There is a short exact sequence of complexes
\[
0\to\tilde{C}_\bullet(Y)\to\tilde{C}_\bullet(\mc{t_1}{Y})\oplus\tilde{C}_\bullet(\mc{t_2}{Y})\to\tilde{C}_\bullet(\mc{\{ t_1,t_2\}}{Y})\to 0
\]
where the second map sends $w$ to $(w,-w)$ and the third map sends $(w_1,w_2)$ to $w_1+w_2$ (where we omit the obvious functorial inclusion maps). There is an induced long exact sequence of homology groups
\[
\xymatrix{
\cdots\ar[r]&\tilde{H}_{n+1}(|\mc{t_1}{Y}|)\oplus\tilde{H}_{n+1}(|\mc{t_2}{Y}|)\ar[r]&\tilde{H}_{n+1}(|\mc{\{ t_1,t_2\}}{Y}|)\ar^-{\delta}[r]&\tilde{H}_n(|Y|)\\
\ar[r]
&\tilde{H}_{n}(|\mc{t_1}{Y}|)\oplus\tilde{H}_{n}(|\mc{t_2}{Y}|)\ar[r]&\cdots&\\
}
\]
Since the groups $\tilde{H}_n(\mc{t_i}{Y})$ vanish for all $n$ by Corollary \ref{lem:contract}, it follows that the connecting homomorphisms $\delta$ induce isomorphisms 
$\tilde{H}_{n+1}(|\mc{\{ t_1,t_2\}}{Y}|)\cong \tilde{H}_n(|Y|)$ for all $n$. By the definition of the homomorphism $\delta$, for any cycle $\gamma\in \tilde{C}_n(Y)$, it sends the class of $(t_1,\gamma)-(t_2,\gamma)$ to the class of $\gamma$. It follows that the chain map $f$ induces the inverse isomorphism.
\end{proof}
%In general $|\mc{\{ t_1,\ldots,t_n\}}{Y}|=\cup_{i=1}^n|\mc{t_i}{Y}|$ where $|\mc{t_i}{Y}|\cap|\mc{t_j}{Y}|=|Y|$ for $i\not=j$. As above, we can label the points of 
%the space $|\mc{\{ t_1,\ldots,t_n\}}{Y}|$ by $(x,\lambda)_i$, with $x\in |Y|$, $\lambda\in [0,1]$, $i=1,\ldots,n$ where $(x,1)_i=t_i$ and $(x,0)_i=x$ for all $x\in |Y|$ and for all $i$.

Now let $\pmc{t}{Y}$ denote the space $|\mc{t}{Y}|/|Y|$, which we will regard as a pointed space with point $\star:=$ the class $(x,0)$ for any $x\in |Y|$.

\begin{lem}\label{lem:homotopy} Let $Y$ be a simplicial complex and let $T$ be a discrete set disjoint from $Y$, with $|T|\geq 2$. Fix $t_0\in T$.
 Let $f=f_T:|\mc{T}{Y}|\to \bigvee_{t\not=t_0}\pmc{t}{Y}$ be the map of spaces defined by
\[
f((x,\lambda)_t):=\left\{
\begin{array}{ll}
\star,& t=t_0\\
(x,\lambda)_t, & t\not= t_0\\
\end{array}
\right.
\]
for $x\in Y$, $\lambda\in [0,1]$.
Then $f$ is a homotopy equivalence.
\end{lem}

\begin{proof} Define $g:\bigvee_{t\not=t_0}\pmc{t}{Y}\to |\mc{T}{Y}|$ as follows: 
\[
g((x,\lambda)_t):=
\left\{
\begin{array}{ll}
(x,2\lambda-1)_t,&\lambda\geq 1/2\\
(x,1-2\lambda)_{t_0}, &\lambda\leq 1/2.\\
\end{array}
\right.
\]
(In particular $g(\star)=g((x,0)_t)=(x,1)_{t_0}=t_0$.)

Thus $gf:|\mc{T}{Y}|\to|\mc{T}{Y}|$ is given by:
\[
gf((x,\lambda)_t)=
\left\{
\begin{array}{ll}
t_0,& t=t_0\\
( x,2\lambda-1)_t,& t\not= t_0, 1\geq \lambda\geq 1/2\\
(x,1-2\lambda)_{t_0},& t\not=t_0,0\leq \lambda\leq 1/2.\\
\end{array}
\right.
\]

Now let  $H_1:|\mc{T}{Y}|\times [0,1]\to |\mc{T}{Y}|$ be given by
\[
H_1((x,\lambda)_t),s):=
\left\{
\begin{array}{ll}
(x,\lambda+s(1-\lambda))_{t_0},&t=t_0\\
\left( x,\frac{2\lambda-s}{2-s}\right)_t&t\not= t_0, \lambda\geq s/2,\\
\left(x,s-2\lambda\right)_{t_0},& t\not= t_0,0\leq  \lambda\leq s/2.\\
\end{array}
\right.
\]
Then $H_1((x,\lambda)_t,0)=(x,\lambda)_t$ for all $t$, while $H_1((x,\lambda)_t,1)=gf((x,\lambda)_t)$  for all $x\in |Y|$, $\lambda\in [0,1]$. Thus $H_1$ is a homotopy from 
$\id{|\mc{T}{Y}|}$ to $gf$.

The map $fg: \bigvee_{t\not=t_0}\pmc{t}{Y}\to \bigvee_{t\not=t_0}\pmc{t}{Y}$ is given by
\[
fg((x,\lambda)_t)=
\left\{
\begin{array}{ll}
(x,2\lambda-1)_t,& \lambda\geq 1/2,\\
\star,& \lambda\leq 1/2.\\
\end{array}
\right.
\]

Let $H_2:\bigvee_{t\not=t_0}\pmc{t}{Y}\times [0,1]\to \bigvee_{t\not=t_0}\pmc{t}{Y}$ be the map
\[
H_2((x,\lambda)_t,s)=
\left\{
\begin{array}{ll}
\left( x,\frac{2\lambda-s}{2-s}\right)_t,&\lambda\geq s/2,\\
\star,&\lambda\leq s/2.\\
\end{array}
\right.
\]
Then $H_2$ is a homotopy from $\id{\bigvee_{t\not=t_0}\pmc{t}{Y}}$ to $fg$ as required.
\end{proof}

\begin{prop}\label{prop:homotopy} Let $Y$ be a simplicial complex and let $T$ ($|T|\geq 2$) be disjoint from $Y$. Fix $p\in Y_0$ and $t_0\in T$. The inclusions of complexes
$\mc{\{ t_0,t\}}{Y}\to \mc{T}{Y}$, $t\not= t_0$ induce a map, $G=G_T$, of simplicial complexes 
\[
G:\bigvee_{t\not=t_0}\mc{\{ t_0,t\} }{Y}\longrightarrow \mc{T}{Y}.
\]
(On the left, the complexes have basepoint the chosen vertex $p\in Y_0$.) Then $G$ induces a homotopy equivalence of geometric realizations, and hence a homotopy equivalence of spaces
\[
\xymatrix{
\bigvee_{t\not=t_0}\left| \mc{\{ t_0,t\}}{Y}\right|\ar^-{\cong}[r]&\left| \bigvee_{t\not=t_0}\mc{\{ t_0,t\} }{Y}\right| \ar^-{\simeq}[r]&|\mc{T}{Y}|.\\
}
\]
\end{prop}
\begin{proof}  For $t\not=t_0$ let $f_t:|\mc{\{ t_0,t\}}{Y}|\to \pmc{t}{Y}$ be the map of Lemma \ref{lem:homotopy}.
We have a commutative diagram of spaces:
\[
\xymatrix{
\left|\bigvee_{t\not=t_0}\mc{\{ t_0,t\}}{Y} \right|\ar^-{|G|}[rr]\ar^-{\cong}[d]& &\left| \mc{T}{Y}\right|\ar^-{f_n}[d]\\
\bigvee_{t\not=t_0}\left| \mc{\{ t_0,t\}}{Y}\right|\ar^-{\vee_t f_t}[rr]&&\bigvee_{t\not=t_0}\pmc{t}{Y}\\
}
\]
Here the left vertical arrow is a homeomorphism. The right vertical arrow is a homotopy equivalence by Lemma \ref{lem:homotopy}. The maps $f_t$ are also a homotopy equivalences by 
Lemma \ref{lem:homotopy}, and hence so is $\vee_t f_t$. It follows, by viewing this as a diagram in the homotopy category, that $|G|$ is also a homotopy equivalence.
\end{proof}

%%%%%%%%%%%%%%%%%%%%%%%%%%%%%%%%%%%%%%%%%%%%%%%%%%%%%%%%%%%%%%%%
\subsection{Complete multipartite graphs}
%%%%%%%%%%%%%%%%%%%%%%%%%%%%%%%%%%%%%%%%%%%%%%%%%%%%%%%%%%%%%%%%%

Let $T_0,T_1,\ldots, T_n$ be a collection of disjoint  sets. We let $K(T_0;\ldots;T_n)$ denote the complete multipartite graph on the sets $T_i$, $i=0,\ldots,n$. Thus the vertices are 
$\cup_{i=0}^n T_i$ and the edges are all pairs $\{ t,s\}$ with $t\in T_i$ and $s\in T_j$ for some $i\not=j$. If $|T_i|=m_i<\infty$ for $0\leq i\leq n$, we also denote (the isomorphism type of) this graph by $K_{m_0,\ldots,m_n}$.

We will let $Y(T_0;\ldots;T_n)$ -- or $Y_{m_0,\ldots,m_n}$ --  denote the clique complex of $K(T_0;\ldots;T_n)$. Thus, by definition, for $r\geq 0$, 
\[
Y(T_0;\ldots;T_n)_r=\{ \{ t_{i_1},\ldots,t_{i_r}\}\ |\ t_{i_j}\in T_{i_j}\mbox{ where }0\leq i_1<\cdots <i_r\leq n\}.
\]
In particular, $Y(T_0;\ldots;T_n)_r=\emptyset$ if $r>n$; it is an $n$-dimensional complex.

Note that $K(T_0)=Y(T_0)$ is just the discrete graph on the set $T_0$.  We will simply denote this by $T_0$.

From the definitions it also follows that, for $n\geq 1$,  $Y(T_0;\ldots;T_n)=\mc{T_n}{Y(T_0;\ldots;T_{n-1})}$ and hence
\[
Y(T_0;\ldots;Y_n)=\mc{T_n}{\mc{T_{n-1}}{\ldots (\mc{T_1}{T_0}})}.
\]

Observe that if $|T_i|=m_i=1$ for any $i\geq 0$, we may suppose without loss that $i=n$ so that $T_n=\{ t\}$ and hence $Y(T_0;\ldots;T_n)=\mc{t}{Y(T_0;\ldots;T_{n-1})}$ is a contractible 
space. 

We consider next the case that $|T_i|=2$ for all $i$.

\begin{prop}\label{prop:mc2}  Suppose that $|T_i|=2$ for $0\leq i\leq n$. Then $Y=Y(T_0;\ldots;Y_n)$ is homeomorphic to an $n$-sphere $S^n$.

 If $T_i=\{ t^i_0,t^i_1\}$ for 
$0\leq i\leq n$ then the cycle \[
\gamma_n:=\sum_{j=0}^n\sum_{i_j\in \{ 0,1\}}(-1)^{i_0+\cdots+i_n}(t^n_{i_n},t^{n-1}_{i_{n-1}},\ldots,t^0_{i_0})\in L_n(Y)\subset C_n(Y)
\]
 represents a generator of the  group $H_n(|Y|)\cong H_n(S^n)\cong \Z$.
\end{prop}

\begin{proof} We have $Y= \mc{T_n}{\mc{T_{n-1}}{\ldots (\mc{T_1}{T_0}})}$ and hence $Y\cong \Sigma^n(T_0)\cong \Sigma_n(S^0)\cong S^n$, using Lemma \ref{lem:susp}.
Certainly $\gamma_0=t^0_0-t^0_1$ is a generator of the reduced homology $\tilde{H}_0(T_0)=\tilde{H}_0(S^0)$. Using Corollary \ref{cor:susp}, whenever $\gamma_{n-1}$ is a generator of $\tilde{H}_{n-1}(Y(T_0;\ldots;T_{n-1}))\cong \tilde{H}_{n-1}(S^{n-1})$, then $\gamma_n=(t^n_0,\gamma_{n-1})-(t^n_1,\gamma_{n-1})$ is a generator of $\tilde{H}_n(Y)$. So the statement follows by induction on $n$.
\end{proof}

\begin{thm}\label{thm:spheresum} Let $T_i$, $i=0,\ldots,n$ be disjoint sets with $|T_i|\geq 2$ for all $i$. Then $Y(T_0;\ldots ;T_n)$ is homotopy equivalent to a pointed sum 
of $n$-spheres.

More precisely:  
%Suppose $T_i=\{ t^i_0,\ldots, t^i_{m_i}\}$ for $0\leq i\leq n$. 

Fix $t_i\in T_i$ for $0\leq i\leq n$.  Let $S:=\prod_{i=0}^n(T_i\setminus\{ t_i\})$.
For $s=(s_0,\ldots,s_n)\in S$, let $Y_s$
%$\mathbf{j}:=(j_0,\ldots,j_n)\in \mathbf{J}_n:=[1,m_0]\times\cdots\times [1,m_n]$ let $Y_\mathbf{j}$ 
be the $n$-sphere 
$Y(\{ t_0,s_0\};\ldots ;\{ t_n,s_n\})$, regarded as a pointed simplicial complex with point $t_0$. Then the map of complexes
\[
\bigvee_{s\in S}Y_s \longrightarrow Y(T_0;\ldots;T_n)
\]
is a homotopy equivalence of spaces.

 In particular, $H_n(Y(T_0;\ldots;T_n))\cong \oplus_{s\in S}\Z$ with basis represented by the family of cycles 
$\{ \gamma^{s}_n\}_{s\in S}$ as in Proposition \ref{prop:mc2}, and furthermore $H_r((Y(T_0;\ldots;T_n))=0$ when $r\not=0,n$.
\end{thm}

\begin{proof} We use Proposition \ref{prop:homotopy} and induction on $n\geq 0$. When $n=0$,  then $S=T_0\setminus\{ t_0\}$ and the map 
\[
\bigvee_{s\in S}\{ t_0,s\}\to T_0
\]
is a bijection of discrete spaces (where the sets on the left are pointed, with point $t_0$).

Now suppose that $n\geq 0$ and that we have $n+2$ sets $T_0,\ldots, T_{n+1}$, with chosen elements $t_0,\ldots,t_{n+1}$. 
Let $S_n:= \prod_{i=0}^n(T_i\setminus{t_i})$ and suppose that the map $\bigvee_{s\in S_n}Y_{s}\to Y(T_0;\ldots;Y_n)$ is a homotopy equivalence.
% Let $T_{n+1}=\{ t^{n+1}_0,\ldots,t^{n+1}_{m_{n+1}}\}$. 
Applying $\mc{T_{n+1}}{-}$ to both sides gives a homotopy equivalence
\[
\bigvee_{s\in S_n}\mc{T_{n+1}}{Y_{s}}\longrightarrow \mc{T_{n+1}}{Y(T_0;\ldots;Y_n)}=Y(T_0;\ldots;T_n;T_{n+1}).
\]

However, Proposition \ref{prop:homotopy}  implies that for each $s\in S_n$ the map
\[
\bigvee_{t\in T_{n+1}\setminus\{ t_{n+1}\}}\mc{\{ t_{n+1},t\}}{Y_s}\longrightarrow \mc{T_{n+1}}{Y_{s}}
\]
is a homotopy equivalence. Since $\mc{\{ t_{n+1},t\}}{Y_{s}}=Y(\{ t_0,s_0\};\ldots ;\{ t_n,s_n\} ;\{ t_{n+1},t\})=Y_{s'}$
where $s'\in S_n\times (T_{n+1}\setminus\{t_{n+1}\}):=S_{n+1}$,
%$\mathbf{j}'\in \mathbf{J}_{n+1}:=\mathbf{J}_n\times [1,m_{n+1}]$,
 it follows that the resulting  composite map
\[
\bigvee_{s'\in S_{n+1}}Y_{s'}\longrightarrow 
\bigvee_{s\in S_n}\mc{T_{n+1}}{Y_{s}}\longrightarrow Y(T_0;\ldots;T_n;T_{n+1})
\]
is a homotopy equivalence as required.
\end{proof}
%\subsection{Infinite complete multipartite graphs} The above results extend easily to the case where some or all of the $T_i$ are infinite sets. 
%%%%%%%%%%%%%%%%%%%%%%%%%%%%%%%%%%%%%%%%%%%%%%%%%%%%%

%%%%%%%%%%%%%%%%%%%%%%%%%%%%%%%%%%%%%%%%%%%%%%%%%%%%%%%%%%%%%%%%%%%%%%%%%%%%%%%%%%

\end{document}